\documentclass{amsart}
 \pdfoutput=1   
\usepackage{amssymb, amsmath, amscd, latexsym, mathrsfs, appendix, mathdots}
\usepackage{graphics}
\usepackage{amsmath}
\usepackage{amsxtra}
\usepackage{amsfonts}
\usepackage{appendix}
\usepackage[all]{xy}
\usepackage{tikz}
\usetikzlibrary{matrix,arrows}

\usepackage{hyperref}
\hypersetup{colorlinks, linkcolor=blue, citecolor=red}

\usepackage{color}

\numberwithin{equation}{section}
\newtheorem{thm}{Theorem}[section]
\newtheorem{cor}[thm]{Corollary}
\newtheorem{lem}[thm]{Lemma}
\newtheorem{prop}[thm]{Proposition}

\theoremstyle{definition}
\newtheorem{defn}[thm]{Definition}

\newtheorem{rem}[thm]{Remark}
\newtheorem{expl}[thm]{Example}

\newtheorem{notn}[thm]{Notation}

\newcommand{\lra}{\longrightarrow}

\newcommand{\co}{\colon\!}

\newcommand{\smin}{\smallsetminus}

\newcommand{\id}{\textup{id}}
\newcommand{\im}{\textup{im}}

\newcommand{\holim}{\textup{holim}}
\newcommand{\hocolim}{\textup{hocolim}}
\newcommand{\colim}{\textup{colim}}

\newcommand{\mor}{\textup{mor}}
\newcommand{\ob}{\textup{ob}}
\newcommand{\map}{\textup{map}}
\newcommand{\mapt}{\textup{mapt}}
\newcommand{\grmap}{\textup{\underline{ma}p}}  
\newcommand{\rgrmap}{\RR \grmap}

\newcommand{\rmap}{\mathbb R\textup{map}}

\newcommand{\lad}{\textup{la}}
\newcommand{\fre}{\textup{cla}}
\newcommand{\prt}{{p}}
\newcommand{\bund}{\textup{bun}}
\newcommand{\fat}{\textup{fat}}

\newcommand{\vep}{\varepsilon}
\newcommand{\lean}{\textup{lean}}
\newcommand{\match}{\textup{match}}
\newcommand{\latch}{\textup{latch}}

\newcommand{\cvar}{\curvearrowright}

\newcommand{\TOP}{\textup{TOP}}

\newcommand{\G}{\textup{G}}

\newcommand{\emb}{\textup{emb}}
\newcommand{\config}{\mathsf{con}} 
\newcommand{\cconfig}{\mathsf{ucon}}
\newcommand{\grcconfig}{\mathsf{grucon}}

\newcommand{\loc}{\textup{loc}}
\newcommand{\imm}{\textup{imm}}

\newcommand{\fin}{\mathsf{Fin}}
\newcommand{\cfin}{\mathsf{uFin}}
\newcommand{\epifin}{\mathsf{biFin}}
\newcommand{\cepifin}{\mathsf{ubiFin}}

\newcommand{\topcat}{\mathsf{Top}}

\newcommand{\sol}{\sharp}
\newcommand{\inj}{\textup{inj}}

\newcommand{\pack}{\textup{pack}}
\newcommand{\Pack}{\textup{Pack}}
\newcommand{\skel}{\textup{sk}}

\newcommand{\sA}{\mathcal A}
\newcommand{\sB}{\mathcal B}

\newcommand{\sC}{\mathcal C}
\newcommand{\sD}{\mathcal D}
\newcommand{\sE}{\mathcal E}

\newcommand{\sK}{\mathcal K}
\newcommand{\sL}{\mathcal L}

\newcommand{\sN}{\mathcal N}
\newcommand{\sP}{\mathcal P}
\newcommand{\sQ}{\mathcal Q}

\newcommand{\sT}{\mathcal T}

\newcommand{\sU}{\mathcal U}
\newcommand{\sV}{\mathcal V}
\newcommand{\sW}{\mathcal W}
\newcommand{\op}{\textup{op}}

\newcommand{\RR}{\mathbb R}

\newcommand{\ZZ}{\mathbb Z}

\newcommand{\twosub}[2]{\begin{array}{cc}
\scriptstyle{#1} \\  [-1mm] \scriptstyle{#2}  \end{array}}

\newcommand{\colimsub}[1]{\begin{array}[t]{cc} \!\!\textup{colim} \\
[-1mm] \scriptstyle{#1} \end{array}}

\newcommand{\limsub}[1]{\begin{array}[t]{cc} \textup{lim} \\
[-1.6mm] \scriptstyle{#1} \end{array}}

\newcommand{\holimsub}[1]{\begin{array}[t]{cc} \textup{holim} \\ [-1mm]
\scriptstyle{#1} \end{array}}

\newcommand{\hocolimsub}[1]{\begin{array}[t]{cc}\!\! \textup{hocolim} \\
[-1.2mm] \scriptstyle{#1} \end{array}}

\newcommand{\uli}{\underline}

\begin{document}

\title{The torus trick for configuration categories}
\author{Pedro Boavida de Brito and Michael S. Weiss}%

\address{Dept. of Mathematics, Instituto Superior Tecnico, Univ.~of Lisbon, Av.~Rovisco Pais, Lisboa, Portugal}%
\email{pedrobbrito@tecnico.ulisboa.pt}

\address{Math. Institut, Universit\"at M\"{u}nster, 48149 M\"{u}nster, Einsteinstrasse 62, Germany}%
 \email{m.weiss@uni-muenster.de}

\thanks{The project was funded by the Deutsche Forschungsgemeinschaft (DFG, German Research Foundation) – Project-ID 427320536 – SFB 1442,
as well as under Germany’s Excellence Strategy EXC 2044 390685587, Mathematics Münster: Dynamics–Geometry–Structure.
P.B. was supported by FCT 2021.01497.CEECIND and grant SFRH/BPD/99841/2014.}

\date{\today}%
\begin{abstract} We show that in codimension at least 3, spaces of locally flat topological embeddings of manifolds
are correctly modeled by derived spaces of maps between their configuration categories
(under mild smoothability conditions). That general claim was reduced in an earlier paper to the special cases where the manifolds in question are euclidean spaces.
We deal with these special cases by comparing to other special cases where the manifolds have the form ``torus'' and
``torus times euclidean space'', respectively, and by setting up a torus trick for
configuration categories.
\end{abstract}
\maketitle

\tableofcontents

\section{Introduction} \label{sec-goal}
Our main theorem relates the space $\emb^t(\RR^\ell,\RR^m)$ of locally flat topological embeddings $\RR^\ell\to \RR^m$
to the space of derived maps (over the nerve of $\fin$) between the configuration categories of $\RR^\ell$ and
$\RR^m$, respectively.
\begin{thm} \label{thm-maineuc} The map $\emb^t(\RR^\ell,\RR^m)\to \rmap_\fin(\config(\RR^\ell),\config(\RR^m))$ determined
by the functoriality of configuration categories is a weak equivalence if $m-\ell\ge 3$ and $m\ge 5$.
\end{thm}
We point out that $\rmap_\fin(\config(\RR^\ell),\config(\RR^m))$ is weakly equivalent to
the space of derived maps from $E_\ell$ to $E_m$, where $E_m$ denotes the little $m$-disk operad.
This was shown in \cite{BoavidaWeiss2018}. Therefore the theorem implies
\[  \emb^t(\RR^\ell,\RR^m) \simeq \rmap(E_\ell,E_m) \]
for $m-\ell\ge 3$ and $m\ge 5$. In this formulation it is more difficult to name the map
which gives the equivalence.

\begin{rem}
Horel shows \cite[Thm. 8.5]{Horel17} that theorem~\ref{thm-maineuc} is valid for $\ell=m=2$. Krannich and Kupers show \cite[Cor E]{KraKu} that theorem~\ref{thm-maineuc} does not extend to the cases $\ell=m$, except for $m\le 2$. They also have counterexamples in codimensions 1 and 2; see \cite[Cor 8.19]{KraKu}.

Salvatore and Turchin \cite[Thm. 2.3]{SalTur}, reacting to an earlier version of this paper, have pointed out that our theorem \ref{thm-maineuc} is also valid in the case $(\ell, m) = (1,4)$ if topological embeddings $\emb^t$ are replaced by PL embeddings $\emb^{pl}$. On the other hand, we can also replace $\emb^t$ by $\emb^{pl}$ throughout in theorem \ref{thm-maineuc}, using \cite{Lashof76}.
\end{rem}

\begin{cor} \label{cor-mainman}
Let $L$ and $M$ be smooth manifolds of dimension $\ell$ and $m$ respectively, with $m-\ell\ge 3$ and $m\ge 5$. Then
\[
 \emb^t(L,M) \to \rmap_\fin(\config(L),\config(M))
\]
has contractible homotopy fibers over each component which comes from a smooth embedding $L \to M$. 
\end{cor}

Something analogous holds for manifolds with boundary, where a smooth embedding $\partial L \to \partial M$ is prescribed.

\smallskip
\emph{Proof, conditional on theorem~\ref{thm-maineuc}.}
Let $\imm^s(L,M)$ and
$\imm^t(L,M)$ be the spaces of smooth resp.~locally flat topological immersions from $L$ to $M$.
In the commutative
diagram
\[
\xymatrix@R=16pt{
\emb^s(L,M) \ar[r] \ar[d] & \emb^t(L,M) \ar[d] \ar[r] & \rmap_\fin(\config(L),\config(M)) \ar[d] \\
\imm^s(L,M) \ar[r] & \imm^t(L,M) \ar[r]^-v & \rmap_\fin(\config^\loc(L),\config^\loc(M))
}
\]
the left-hand square is a homotopy pullback square by \cite[Thm.A]{Lashof76} and
the outer square is a homotopy pullback square by \cite[Cor.~5.2]{BoavidaWeiss2018}.
The arrow labeled $v$ is a weak equivalence by theorem~\ref{thm-maineuc}. (Source and target
of $v$ can be described as section spaces of certain fibrations over $L$. In the case of the source, the fiber over $x\in L$
is, informally speaking, the space of pairs $(y,f)$ where $y\in M$ and $f$ is a locally flat embedding $T_xL\to T_yM$. In the case of the
target, the fiber over $x\in L$ is the space of pairs $(y,f)$ where $y\in M$ and $f$ is a derived map from
$\config(T_xL)$ to $\config(T_xM)$.) The claim follows with some diagram diagnostics. \qed

\medskip
\begin{rem} Krannich-Kupers \cite[Thm. B]{KraKu2} have recently announced a convergence result for the topological embedding calculus tower in the case where the target manifold is smoothable and of dimension $\geq 5$, and the codimension is at least three. That implies the map in corollary \ref{cor-mainman} is a weak equivalence and the assumption that $L$ is smoothable is not necessary. (Convergence is essentially equivalent to the statement that the right-hand square in the diagram of the proof is homotopy cartesian.)
\end{rem}

\smallskip
\emph{Remark.} If $L$ is compact, and $m-\ell\ge 3$, $m\ge 5$,
then the space of locally flat embeddings $\emb^t(L,M)$ is weakly equivalent to the space of
all injective continuous maps from $L$ to $M$. See \cite[Appendix]{Lashof76}. Note in passing that we normally think of
$\emb^t(L,M)$ as a simplicial set, while the space of all injective continuous maps from $L$ to $M$ could also be regarded as an
honest space with the compact-open topology; more about that in \cite{Lashof76}. We do not know whether $\emb^t(\RR^\ell,\RR^m)$ is
weakly equivalent to the space of all injective continuous maps from $\RR^\ell$ to $\RR^m$. Theorem~\ref{thm-maineuc} implies that one is a homotopy
retract of the other since the ``functoriality of configuration categories'' is perfectly valid for continuous injective maps.

\medskip
We end this introduction with another corollary of theorem \ref{thm-maineuc}. Let $G(d)$ denote the space of homotopy automorphisms of $S^{d-1}$ and $G = \colim_d G(d)$.
 
\begin{cor} For $d \geq 3$, we have
$
\colim_{\ell} ~\rmap(E_\ell, E_{\ell + d}) \simeq \G / \G(d) \; .
$
\end{cor}
\begin{proof}
It is known \cite[p.147]{Lashof76} that the map
\[
\colimsub{\ell} \emb^t(\RR^\ell, \RR^{\ell + d}) \to \colimsub{\ell} \G(\ell+d)/\G(d)
\]
which classifies the normal spherical fibration is a weak equivalence if $d \geq 3$.
\end{proof}

\smallskip
\emph{Guide to the paper.} The torus trick is explained in section~\ref{sec-totri}. It is easy to understand as a method, but the
preparations required to make it work are sobering. The most important of these help us to translate metric conditions on maps
between configuration categories into homotopical conditions. For that we have mainly
sections~\ref{sec-pseudomet} and~\ref{sec-bolift}. Sections~\ref{sec-smallbig} and~\ref{sec-grobo} support
sections~\ref{sec-pseudomet} and~\ref{sec-bolift} respectively by dropping some useful anchors.
In section~\ref{sec-pusimp} we develop a homotopical decomposition method
for configuration categories, based on partitions of unity. This is needed in sections~\ref{sec-grobo} and~\ref{sec-bolift}.
In the appendix sections~\ref{sec-models},~\ref{sec-homodels},~\ref{sec-joins} and~\ref{sec-fibreal}
we introduce language, models and points of view. Section~\ref{sec-fibreal} introduces a
new procedure for conservatization. (This is new only in relation to the conservatization procedure used in \cite[\S8]{BoavidaWeiss2018}.)
\section{The torus trick} \label{sec-totri}
\subsection{A form of torus trick for embeddings.} Let $T^\ell$ be the $\ell$-dimensional torus.
Let $f\co T^\ell \to T^\ell \times \RR^d$ be a locally flat topological embedding, and let $H$ be a homotopy from $f$ to
the standard inclusion. For any self-covering map $\pi: T^\ell \to T^\ell$, path lifting determines a unique lift of the homotopy $H$
to a homotopy $H^\prime$ satisfying $(\pi \times {\id}) H^\prime = H \pi$ with $H^\prime_1$
the standard inclusion. Then $H^\prime_0$ is an embedding that lifts $f$:
\[
\xymatrix@M=5pt@R=18pt{
T^\ell \ar[d]^{\pi} \ar[r]^-{H'_0} &T^\ell \times \RR^d \ar[d]^{\pi \times \id} \\
T^\ell \ar[r]^-{f} & T^\ell \times \RR^d
}
\]

\medskip
The following instance of a geometric torus trick has guided this paper. Think of $T^\ell$ as a topological abelian group
(hence a $\ZZ$-module) and suppose $d\ge 3$, $\ell+d\ge 5$.

\begin{lem}\label{lem:thetorustrick}
Let $f^{(k)}$ be the lift of $f$ across the covering map $\pi\co T^\ell\to T^\ell$ where $\pi(x)=kx$ for some positive
integer $k$. If $k$ is sufficiently large, then the lift $f^{(k)}$ is isotopic to the standard inclusion.
\end{lem}
(The meaning of \emph{sufficiently large} depends on $f$ and $H$.)
The argument for this can be broken up into several steps.

\begin{enumerate}
\item $f$ is $\varepsilon$-\emph{bounded} for some $\varepsilon$ (which can be large)
\item if $f$ is $\varepsilon$-bounded then $f^{(k)}$ is $\varepsilon/k$-bounded
\item therefore $f^{(k)}$ is \emph{isotopic} to the standard inclusion if $k$ is large enough.
\end{enumerate}

By an $\varepsilon$-\emph{bounded} map $f: T^\ell \to T^\ell \times \RR^d$ (equipped with a homotopy to the standard inclusion)
we mean a map whose lift to the universal cover $${f}^\infty : \RR^\ell \to \RR^\ell \times \RR^d$$ is $\varepsilon$-close
to the standard inclusion $e$ in the sense that the distance between $p_1e(x)$ and $p_1f^\infty(x)$ is $<\varepsilon$ for all $x\in \RR^\ell$,
where $p_1\co \RR^\ell\times\RR^d\to \RR^\ell$
is the first projection. We are not very interested in the distance between $p_2e(x)$ and $p_2f^\infty(x)$ since we can easily make it
as small as we like (for all $x$ simultaneously) by applying a suitable isotopy to $f$ (shrinking in the $\RR^d$ factor).

Points (1) and (2) are clear. For (3) we must refer to the appendix of \cite{Lashof76}, and specifically to theorems 1 and 2 in there.
Together these imply a complicated local contractibility property for spaces of locally flat embeddings of topological manifolds in codimension $\ge 3$.
Lashof attributes this to \cite{Cernavskii}.
Here the consequence is that $f^{(k)}$ for sufficiently large $k$ admits a ``small'' isotopy to the standard inclusion, but that isotopy might not be
locally flat. On the other hand there is a theorem saying that, if there is such an isotopy, then there is also one which is locally flat,
although that one might not be small.

Lemma \ref{lem:thetorustrick} has a family version where the embedding $f$  is replaced by a compact family of embeddings. We restate
this in more homotopical terms, after fixing some language and notation.

\begin{defn}
A topological embedding $T^\ell \to T^\ell \times \RR^d$ together with a homotopy to the standard inclusion will be called a \emph{grounded embedding}. The space of such is denoted $\underline{\emb}^{t}(T^\ell, T^\ell \times \RR^d )$. It is the homotopy fiber of
\[
\emb^{t}(T^\ell, T^\ell \times \RR^d ) \to \map(T^\ell, T^\ell \times \RR^d)
\]
over the standard inclusion.
\end{defn}

As we noted already, grounded embeddings lift uniquely across self-coverings of the torus. For a covering $\pi : T^\ell \to T^\ell$, this gives a map
\[
\underline{\emb}^{t}(T^\ell, T^\ell \times \RR^d) \xrightarrow{\pi^*} \underline{\emb}^{t}(T^\ell, T^\ell \times \RR^d) \; .
\]
The family version of Lemma \ref{lem:thetorustrick} is then:

\begin{prop}\label{prop:torustrickT} Let $\pi : T^\ell \to T^\ell$ be the covering map given by $\pi(x)=2x$. The homotopy colimit of the tower
\[
\underline{\emb}^{t}(T^\ell, T^\ell \times \RR^d)  \xrightarrow{\pi^*} \underline{\emb}^{t}(T^\ell, T^\ell \times \RR^d)  \xrightarrow{\pi^*} \dots
\]
is contractible.
\end{prop}

\subsection{The torus trick for configuration categories}
For the purposes of this section, we make the following abbreviations. For topological manifolds $L$ and $M$
let $\bund^s(TL,TM)$ be the space of pairs $(f,g)$ where $f\co L\to M$
is any map and $g\co TL\to TM$ covers $f$, and is fiberwise linear and injective.
The definition of $\bund^t(TL,TM)$ is similar, but we replace the linearity condition by local flatness. We write
\[
\emb^\prt(L,M) := \RR \map_{\fin}(\config(L), \config(M)),
\]
\[
\bund^\prt(TL,TM) := \RR \map_{\fin}(\config^\loc(L), \config^\loc(M)).
\]
The superscript $\prt$ is a reference to particles or patches.
Another justification (for the $\bund^\prt$ notation) comes from noting that $\config^\loc(L)$ is in some sense
sense a bundle of configuration categories $\config(T_xL)$ where $x\in L$.

\begin{defn}Let $\underline{\emb}^\prt(T^\ell,T^\ell \times \RR^d)$ be the homotopy fiber, over the standard inclusion, of the map
\[
\emb^\prt(T^\ell, T^\ell \times \RR^d ) \to \map(T^\ell, T^\ell \times \RR^d)
\]
given by restriction to the space of objects of $\config(T^\ell)$ over $\uli 1$ in $\fin$. We again refer to
the homotopy as a \emph{grounding} and the elements as \emph{grounded} things.
\end{defn}

\begin{rem}\label{rem:groundeddown}
The forgetful map
\[
\underline{\emb}^\prt(T^\ell, T^\ell \times \RR^d) \to \emb^\prt(T^\ell, T^\ell \times \RR^d)
\]
is a fibration with homotopically discrete fibers. An alternative description of the elements of $\underline{\emb}^\prt(T^\ell, T^\ell \times \RR^d)$
is as pairs $(f, f^\infty)$ where $f \in \emb^\prt(T^\ell, T^\ell \times \RR^d)$ and
\[ f^\infty : \RR^\ell \to \RR^\ell \times \RR^d \]
is a lift of the map $T^\ell\to T^\ell\times\RR^d$ determined by $f$ to the universal covers.
This variant will be used later in the text.
\end{rem}

In \cite{confcover}, we have shown that derived maps of configuration categories can be lifted to maps
of configuration categories of covering spaces under mild assumptions.
Consequently, given a covering map $\pi : T^\ell \to T^\ell$, there is a dashed arrow
\[
\xymatrix@M=5pt@R=20pt{
\underline{\emb}^t(T^\ell,T^\ell \times \RR^d) \ar[d] \ar[r]^{\pi^*} & \underline{\emb}^t(T^\ell,T^\ell \times \RR^d) \ar[d] \\
\underline{\emb}^\prt(T^\ell,T^\ell \times \RR^d) \ar@{..>}[r]^{\pi^*} &  \underline{\emb}^\prt(T^\ell,T^\ell \times \RR^d)
}
\]
making the diagram homotopy commutative.

The main technical result of this paper is then:
\begin{thm}[Torus trick for configuration categories]\label{thm:torustrickE} Let $\pi : T^\ell \to T^\ell$ be the degree $2$ covering of the torus. The homotopy colimit of the tower
\[
\underline{\emb}^\prt(T^\ell, T^\ell \times \RR^d) \xrightarrow{\pi^*} \underline{\emb}^\prt(T^\ell, T^\ell \times \RR^d) \xrightarrow{\pi^*} \dots
\]
is contractible.
\end{thm}

The remaining sections in the paper taken together prove this theorem. We will indicate how in a moment.
Before we do so, let us explain how theorem \ref{thm:torustrickE} implies theorem \ref{thm-maineuc}.

\subsection{Proof of theorem \ref{thm-maineuc}} \label{subsec-oldtrick}
To shorten notation, we abbreviate
\[
L = T^\ell \quad \textup{ and } \quad M = T^\ell \times \RR^d
\]
(and $m:=\ell+d$).
In analogy with the definitions and notation from before, let
\[
\underline{\bund}^?(TL,TM)
\]
be the homotopy fiber of $\bund^?(TL,TM) \to \map(L,M)$ where $?$ can be $s$, $t$ or $\prt$.

The master diagram, relating the different types of spaces, is:
\[
\xymatrix@M=5pt@R=20pt{
\underline{\emb}^s(L,M) \ar[d] \ar[r] & \underline{\emb}^t(L,M) \ar[d] \ar[r] & \underline{\emb}^\prt(L, M) \ar[d] \\
\underline{\bund}^s(TL,TM) \ar[r] & \underline{\bund}^t(TL,TM) \ar[r] &  \underline{\bund}^\prt(TL,TM)
}
\]
This has a forgetful map to a similar diagram where no grounding is imposed (the underlining in the terms disappears).
The left-hand square is homotopy cartesian by \cite[Thm.A]{Lashof76} and the outer rectangle is homotopy cartesian
by \cite[Thm.5.1]{BoavidaWeiss2018}. Therefore the right-hand square is almost homotopy cartesian
(i.e., the vertical homotopy fiber over the base point in the middle column maps by a weak equivalence to the vertical homotopy fiber over
the base point in the right-hand column).

These diagrams are the layers in a tower which is obtained by invoking the self-covering map $\pi\co L\to L$ of proposition~\ref{prop:torustrickT}
and letting it act repeatedly. This extends the towers for $\underline{\emb}^s(L,M)$, $\underline{\emb}^t(L,M)$ and $\underline{\emb}^\prt(L,M)$.
We have explained this for the arrows in the top row (and it is clear for the left-hand square).
For the remaining arrows, it is a consequence of two observations (for more details, see \cite{BoavidaWeiss2018}).
First, the lifting-across-$\pi$ map $$\pi^* : \underline{\emb}^\prt(L,M) \to  \underline{\emb}^\prt(L,M)$$ can be made
natural with respect to inclusions of open subsets of $L$. Second, the lower row is obtained functorially from the upper row by
homotopy sheafification wrt the ordinary
notion of open cover.

Passing to the top of the tower, i.e., taking the sequential homotopy colimit in each column,
associated to the right-hand square in the diagram above we obtain a square
\[
\xymatrix@M=5pt@R=20pt{
{\colimsub{\textup{iterated } \pi^*}} \underline{\emb}^t(L,M) \ar[d] \ar[r] & {\colimsub{\textup{iterated } \pi^*}}  \underline{\emb}^\prt(L, M) \ar[d] \\
{\colimsub{\textup{iterated } \pi^*}} \underline{\bund}^t(TL,TM) \ar[r] & {\colimsub{\textup{iterated } \pi^*}}  \underline{\bund}^\prt(TL,TM) \; .
}
\]
Since homotopy pullbacks commute with directed homotopy colimits, this is still almost homotopy cartesian; i.e., the induced map on vertical homotopy fibers over the respective
base points is a weak equivalence. Since the terms in the top row are contractible, it follows that the lower horizontal map is a weak equivalence on
base point components.

The last step is easy. Since $L$ and $M$ are parallelized, the map
\[
{\colimsub{\textup{iterated } \pi^*}} \underline{\bund}^t(TL,TM) \to {\colimsub{\textup{iterated } \pi^*}}  \underline{\bund}^\prt(TL,TM)
\]
is weakly equivalent to
\[
{\colimsub{\textup{iterated } \pi^*}} \map(L, \emb^t(\RR^\ell, \RR^m)) \to {\colimsub{\textup{iterated } \pi^*}}
\map(L, \emb^\prt(\RR^\ell, \RR^m).
\]
(The maps in these towers are given by precomposition with $\pi$.) The last map has the standard comparison map
\[
\emb^t(\RR^\ell, \RR^m) \to \emb^\prt(\RR^\ell, \RR^m)
\]
as a retract (use restriction of maps out of $L$ to the base point, and inclusion of constant maps). Therefore we may conclude
that the said comparison map is a weak equivalence on base point components. But these spaces happen to be connected.
For the target space, this follows from \cite{Goeppl}, and here again the condition $m-\ell\ge 3$ is important.
For the source, our local flatness assumption implies that the map $\TOP(m)\to \emb^t(\RR^\ell, \RR^m)$ given by restriction
induces a surjection of path components. Famously $\TOP(m)$ has only two path components, and clearly these determine the same path component of
$\emb^t(\RR^\ell, \RR^m)$.

\subsection{Proof of Theorem \ref{thm:torustrickE}}

We keep the abbreviations of section~\ref{subsec-oldtrick}. For reasons given in \cite{confpres}, we prefer to work with the Rezk completions of
the Segal spaces $\config(L)$ and $\config(M)$. They are denoted $\cconfig(L)$ and $\cconfig(M)$, respectively, where the ``u'' is for \emph{unordered}
as in unordered configuration.
The manifolds $L$ and $M$ will be regarded as Riemannian manifolds (with the standard flat Riemannian metrics).
This allows us to use the Riemannian multipatch models for $\cconfig(L)$ and $\cconfig(M)$, which are
more useful here than the particle models. The multipatch models are described in \cite{BoavidaWeiss2018} and/or in appendix~\ref{sec-models} below.
In the Riemannian multipatch model, $\cconfig(L)$ and $\cconfig(M)$ are nerves of topological posets (whose elements
are the multipatches). Then we can provisionally re-define
\[  \emb^\prt(L,M) := \rmap_\cfin(\cconfig(L),\cconfig(M))= \map_\cfin(\cconfig(L),\varphi\cconfig(M)) \]
where $\varphi$ is for a Reedy fibrant replacement (whereas $\cconfig(L)$ is already Reedy cofibrant). This is in agreement with
the previous definition of $\emb^\prt(L,M)$ up to weak equivalence.

As a first step, we restrict cardinalities. That is, we replace the configuration category $\config(L)$ by the
subcategory $\config(L; \alpha)$. For each $\alpha$, we have a tower
\[
\underline{\emb}^\prt(L,M; \alpha) \xrightarrow{\pi^*} \underline{\emb}^\prt(L,M; \alpha) \xrightarrow{\pi^*} \dots
\]
(in self-explanatory notation). It suffices to show that the colimit of each of these towers is contractible. With that in mind, we fix an $\alpha$ throughout.

In a similar vein, we often find it helpful to restrict the size of patches in a multipatch
in $L$. This leads to notation like $\cconfig_\delta(L;\alpha)$, where $\delta$ is an upper bound on the radius of patches.
The inclusion of $\cconfig_\delta(L;\alpha)$ in $\cconfig(L;\alpha)$ is a (degreewise) weak equivalence. We seize the
opportunity to re-define once more
\[ \emb^\prt(L,M) := \colimsub{\delta\to 0} \map_\cfin(\cconfig_\delta(L),\varphi\cconfig(M)). \]
(If we have to be precise, the colimit is taken in the category of simplicial sets.
There will be no further re-definitions of $\emb^\prt(L,M)$ in this section.)

Next, there is a notion of $\vep$-boundedness for elements (more correctly, simplices) of $\underline{\emb}^\prt(L,M)$. This
is very similar to the notion of $\vep$-boundedness for elements of $\underline{\emb}^t(L,M)$. It is slightly more complicated though;
see definitions~\ref{defn-eprcontrol} and~\ref{defn-swepcontrol}.

\medskip
Now the necessary definitions are in place and the strategy can be outlined. It is close in spirit to the topological torus trick above.
Let $K$ be a finitely generated simplicial set and let $f\co K\to \emb^\prt(L,M)$ be a map of simplicial sets.
We like to think of $f$ as a family $(f_u)_{u\in K}$.
\begin{enumerate}
\item After a homotopy applied to the family,
each $f_u$ is $\varepsilon$-bounded for some
$\vep$ independent of $u\in K$.
\item If a simplex in $\underline{\emb}^\prt(L,M; \alpha)$ is $\varepsilon$-bounded, and $k$ is any positive integer,
then for sufficiently large $r$ (depending on $k$) the image of that simplex under the $r$-fold iteration of
\[ \pi^*\co \underline{\emb}^\prt(L, M; \alpha)\to \underline{\emb}^\prt(L, M; \alpha) \]
is $(\vep/k)$-bounded.
\item If $\vep$ in (1) is sufficiently small, then the family $(f_u)$ is nullhomotopic.
(There is a preferred  base point in $\underline{\emb}^\prt(L, M; \alpha)$.)
\end{enumerate}
Each of these steps is a major undertaking. Step (1) is achieved in theorem~\ref{thm-grbound}. Step (2) is a
consequence of theorem~\ref{thm-step2} and the observation, justified in the last section of~\cite{confcover},
that the $r$-fold iteration of
$\pi^*\co \underline{\emb}^\prt(L, M; \alpha)\to \underline{\emb}^\prt(L, M; \alpha)$ is homotopic to the map
determined by lifting across the covering map
\[ \underbrace{\pi\circ\pi\circ\cdots\circ\pi}_r\co T^\ell\lra T^\ell\,. \]
Step (3) is theorem~\ref{thm-shrink}. Strictly speaking, theorem~\ref{thm-shrink} makes a claim which
looks slightly weaker than (3). It says that under conditions as in (3), the composition
\begin{equation}
  K \xrightarrow{\,u\mapsto f_u\,} \underline{\emb}^\prt(L, M; \alpha) \xrightarrow{\textup{forgetful}}  \emb^\prt(L, M; \alpha)
  \end{equation}
is nullhomotopic. In order to show that this does not make any difference, we expand the map
$\underline{\emb}^\prt(L, M; \alpha) \to  \emb^\prt(L, M; \alpha)$ into a Barratt-Puppe sequence. This gives
\[  \cdots\to \Omega\emb^\prt(L, M; \alpha) \to \Omega\map(L,M) \to \underline{\emb}^\prt(L, M; \alpha) \to  \emb^\prt(L, M; \alpha). \]
The map $\Omega\emb^\prt(L, M; \alpha) \to \Omega\map(L,M)$ in the sequence has a right homotopy inverse, since the base point
component of $\map(L,M)$ is weakly equivalent to $L$
and since $L$, being a topological abelian group, is a retract of $\emb^\prt(L,M)$. Therefore, if a map
from a simplicial set $K$ to $\underline{\emb}^\prt(L, M; \alpha)$ is nullhomotopic as a map to
$\emb^\prt(L, M; \alpha)$, then it was already nullhomotopic to begin with.

\section{Lean and fat multipatches}  \label{sec-smallbig}
\emph{Vocabulary and notation}: see section~\ref{sec-models}.

\begin{lem} \label{lem-upperb} \emph{i)} Let $\sW$ be an open $\alpha$-cover of the closed Riemannian manifold $L$. The inclusion
$\cconfig_\sW(L;\alpha)\to \cconfig(L;\alpha)$ is a weak equivalence. \emph{ii)} Let $\delta$ be a positive real number.
The inclusion $\cconfig_\delta(L;\alpha)\to \cconfig(L;\alpha)$ is a weak equivalence.
\qed
\end{lem}

\begin{lem} \label{lem-alphacover} Let $L$ be a closed smooth Riemannian manifold
and let $\sV$ be an open $\alpha$-cover of $L$. Then there exists $\delta>0$ such that
for every finite subset $S\subset L$ where $|S|\le \alpha$, the open $\delta$-neighborhood
of $S$ in $L$ is contained in some $U\in \sV$. \qed
\end{lem}
(In the case $\alpha=1$, this lemma is a special case
of Lebesgue's covering lemma.) The important implication of lemma~\ref{lem-alphacover} for us is that
if we wish to make a sensible selection of objects of $\config(L;\alpha)$ by imposing upper bounds on size,
as in lemma~\ref{lem-upperb}, then we can often do so by just imposing a uniform upper bound on
the radii of the (multi-)patches.

\medskip
In the remainder of the section we mainly ask how far we can go in imposing \emph{lower bounds}
on the size of objects of $\cconfig(L;\alpha)$ or $\config(L;\alpha)$, if we still wish to have a sensible selection.
This turns out to be a much more difficult topic.

Again let $L$ be a closed smooth Riemannian manifold. Let $k$ be a positive integer.
For $c>0$ let $\pack(k,c)\subset \emb(\uli k,L)$ be the open subset consisting of the
embeddings $f\co \uli k\to L$ such that $d_L(f(x),f(y))>2c$
whenever $x,y\in \uli k$ are distinct (where $d_L$ is the geodesic distance). If $c$ is less than
the global injectivity radius of $L$, then it is alright to think of $\pack(k,c)$ as the
space of multipatches in $L$ (as in definition~\ref{defn-models1}) with exactly $k$ components, where each component has radius exactly $c$.

\begin{prop} \label{prop-elbow} If $c$ is sufficiently small, then the inclusion
\[ \pack(k,c)\to \emb(\uli k,L) \]
is a homotopy equivalence.
\end{prop}
\proof Note first that $\emb(\uli k\,,L)$ is a smooth manifold and $\pack(k,c)$ is an open subset of it.
We may assume from the outset that $6c$ is less than the global injectivity radius of the Riemannian manifold $L$.
Call a tangent vector $v$ to $f\in \pack(k,c)$ \emph{admissible} (in this proof) if, for every
curve $\gamma\co J\to \pack(k,c)$ (where $J\subset\RR$ is an open interval containing $0$) having $\gamma(0)=f$
and $\gamma'(0)=v$, and every choice of distinct $x,y\in \uli k$ such that $d_L(f(x),f(y))\le 6c$, the distance in question
increases along $\gamma$. More precisely, we wish to have
\[   \frac{d}{dt}\Big|_{t=0} ~d_L(\gamma(t)(x),\gamma(t)(y)) > 0 \]
for such $x$ and $y$.
(The differential quotient exists because of the condition on $c$.) Clearly if elements $v,w$ of the
tangent space  $T_f\pack(k,c)$ are admissible, then $v+w$ is also admissible, and $sv$ is admissible for any
positive $s\in \RR$. The admissibility condition is an \emph{open condition}, i.e., the set of admissible elements
in $T\pack(k,c)$ (total space of the tangent bundle of $\pack(k,c)$) is an open set. If
$f\in \pack(k,3c)\subset \pack(k,c)$, then \emph{every} $v\in T_f\pack(k,c)$ is admissible. --- The next thing to observe
is that if $c$ is sufficiently small, then for every $f\in \pack(k,c)$ there exists an admissible tangent vector
in $T_f \pack(k,c)$\,. To see this, let $R$ be the smallest equivalence relation on $\uli k$ which contains all
pairs $(x,y)\in \uli k\times\uli k$ such
that $d_L(f(x),f(y))\le 6c$. We may assume from now on that $6kc$ is smaller than the global
injectivity radius of $L$. Each
equivalence class $A$ of $R$ is contained in an open metric ball $B_A$ of radius $6kc$ about some point $z_A\in L$.
These balls need not be disjoint,
but it is clear that if $A_1$ and $A_2$ are distinct equivalence classes and $x\in A_1$, $y\in A_2$, then $d_L(x,y)>6c$.
Therefore in our effort to find an admissible $v\in T_f(c,k)$ we may treat each equivalence class $A$ of $R$
separately. By applying a ``radial expansion'' of $B_A$ with center point $z_A$, and restricting that to $f(A)$,
we obtain a curve $\gamma_A$ in $\emb(A,L)$ such that $\gamma(0)=f|_A$
and such that $\gamma_A'(0)$ has the required positivity properties, as far as elements $x,y$ of $A$ are concerned,
provided $6kc$ is sufficiently small. (Here the smallness requirement could be quantified in terms of the curvature
properties of the Riemannian metric,
in addition to the global injectivity radius.)
Together the vectors $\gamma_A'(0)$
make up a tangent vector $v\in T_f \pack(k,c)$ which is admissible. ---
Using all that and partitions of unity, one can easily construct a smooth vector field $\xi$ on $\pack(k,c)$ which is everywhere
admissible. The flow of $\xi$ is ``forward complete''. More precisely, if $f\in \pack(k,c)$, then by definition
of $\pack(k,c)$ there exists a positive $\delta$ such that $f\in \pack(c+\delta,k)$ and $\delta<2c$.
The integral curve of $\xi$ passing through $f\in \pack(k,c)$ at time $t=0$ can never escape from the closure
of $\pack(c+\delta,k)$ in $\pack(k,c)$, which is compact. More to the point, if $C$ is any compact subset
of $\pack(k,c)$, then the flow of $\xi$ will move it into the open subset $\pack(k,3c)$ in finite time.
--- Now let $W\subset \pack(k,c)$ be defined as follows: $f\in \pack(k,c)$ belongs to
$W$ if and only if $d_L(x,y)>6c$ for all distinct $x,y\in \{1,2,3,\dots,k-1\}$. Then we have a diagram of inclusion maps
\[  \pack(k,3c) \xrightarrow{g_1} W \xrightarrow{g_2} \pack(k,c) \xrightarrow{g_3} \emb(\uli k\,,L). \]
We still want to show that $g_3$ is a homotopy equivalence.
The forgetful projection $W\to \pack(k-1,3c)$ is a fiber bundle. (Each fiber is
the complement in $L$ of $k-1$ pairwise disjoint metric closed balls of radius $2c$; the balls are disjoint
because their center points have distance $>6c$.)
By induction, and making $c$ smaller if necessary, we may assume that the inclusion of $\pack(k-1,3c)$ in
$\emb(\underline{k-1}\,,L)$ is a
homotopy equivalence. It follows that $g_3g_2$ in the above diagram
is a homotopy equivalence. (We can describe $g_3g_2$ as a map between the total spaces of two fiber bundles
which respects the bundle projections, induces a homotopy equivalence of the base spaces, and restricts to
homotopy equivalences between corresponding fibers.) Now we can conclude with a formal argument. For every \emph{compact}
CW-space $Z$, the map $g_2$ induces a surjection, thanks to the vector field $\xi$, from the set of
homotopy classes $[Z,W]$ to $[Z,\pack(k,c)]$. As we have seen, $g_3g_2$ induces a bijection from
$[Z,W]$ to $[Z,\emb(\uli k\,,L)]$. Therefore $g_3$ induces a bijection
from $[Z,\pack(k,c)]$ to $[Z,\emb(\uli k\,,L)]$. By itself that does not allow us to conclude that $g_3$ is a
weak homotopy equivalence. But we know also that the target of $g_3$ is homotopy equivalent to a compact
CW-space, so that $g_3$ has a homotopy right inverse. With that it is an exercise to show that $g_3$ is a weak equivalence.
Then it is also a homotopy equivalence, since source and target are finite dimensional smooth manifolds. \qed

\smallskip
There is a mild generalization of proposition~\ref{prop-elbow} in which $L$ is allowed to have a boundary.
In the generalization, suppose that $L$ is a compact smooth Riemannian manifold with boundary and let $k$ be a positive integer.
For $c>0$ let $\pack(k,c)$ be the open subset of $\emb(\uli k,L\smin\partial L)$ consisting of the
embeddings $f\co \uli k\to L\smin\partial L$ such that $d_L(f(x),f(y))>2c$
whenever $x,y\in \uli k$ are distinct (where $d_L$ is the geodesic distance), and also $d_L(f(x),z)>c$
for any $z\in \partial L$ and $x\in \uli k$.

\begin{prop} \label{prop-elbowbdry}
For sufficiently small $c$, the inclusion
\[ \pack(k,c)\to \emb(\uli k,L\smin\partial L) \]
is a homotopy equivalence.
\end{prop}
The proof is left to the reader. It can be modeled on the proof of proposition~\ref{prop-elbow}.
We will only need this in the cases where $L$ is a disk, though not necessarily a disk with
the standard Riemannian metric. \qed

\medskip
It is a easy to think of generalizations and variants of proposition~\ref{prop-elbow} and~\ref{prop-elbowbdry}
where the multipatches are replaced by nested systems of multipatches. We will also need such
variants (in the proof of theorem~\ref{thm-shrink} below). To specify the type of nested
system, we fix a diagram in $\fin$,
\begin{equation} \label{eqn-systype}
  \sD=\left(\uli k_0 \leftarrow \uli k_1\leftarrow \cdots\leftarrow \uli k_{p-1} \leftarrow \uli k_p\right)
\end{equation}
(where $k_0,k_1,\dots,k_p$ are positive).

Let $c_0,c_1,\dots,c_p$ be a string of positive real numbers. Suppose that $c_0$ is
less than the global injectivity radius of $L$. Let $\Pack(\sD,c_0,c_1,\dots,c_p)$
be the space of systems of multipatches
\begin{equation} \label{eqn-syspatch} U_0\supset U_1 \supset \dots \supset U_{p-1} \supset U_p \end{equation}
in $L$ where the string of finite sets and maps obtained by applying $\pi_0$
to~\eqref{eqn-syspatch} is identified with $\sD$ in~\eqref{eqn-systype},
and moreover, for $j\in\{0,1,\dots,p\}$, the individual patches in $U_j$ have radii
$\ge c_0c_1\cdots c_{j-1}c_j$. (In the case $p=0$, we can write $\uli k_0$ instead of $\sD$.
There is a small distinction between
$\Pack(\uli k_0,c_0)$ and $\pack(k_0,c_0)$ because in one of them we allow patches whose radii
are at least $c_0$, whereas in the other one the patches must have radii equal to $c_0$. But the
inclusion $\pack(k_0,c_0)\hookrightarrow \Pack(\uli k_0,c_0)$ is clearly a homotopy equivalence.)

\begin{thm} \label{thm-elbow} If $c_0,c_1,\dots,c_p$ are sufficiently small, then the inclusion
\[ \Pack(\sD,c_0,c_1,\dots,c_p) \hookrightarrow \config(L)_\sD \]
is a homotopy equivalence.
\end{thm}
\proof (See remark~\ref{rem-elbow} for an overview.)
Let
\[ \Pack_1(\sD,c_0,c_1,\dots,c_p)\subset \Pack(\sD,c_0,c_1,\dots,c_p) \]
be the subspace
determined by the additional condition that the innermost multipatch, whose name is $U_p$ in~\eqref{eqn-syspatch},
must have all patch radii \emph{equal} to $c_0c_1\cdots c_p$.
It is clear that the inclusion of $\Pack_1(\sD,c_0,c_1,\dots,c_p)$ in $\Pack(\sD,c_0,c_1,\dots,c_p)$
is a homotopy equivalence. Therefore it is enough to investigate $\Pack_1(\sD,c_0,c_1,\dots,c_p)$.

Let $d_p\sD$ be the diagram obtained from $\sD$ by deleting the object $\uli k_p$. There is a forgetful map
\begin{equation} \label{eqn-funnyind1} \Pack_1(\sD,c_0,c_1,\dots,c_{p-1},c_p) \lra
\Pack(d_p\sD,c_0,c_1,\dots,c_{p-1}). \end{equation}
If we can show that this is a \emph{weak Serre fibration}, see definition~\ref{defn-wkSfib} below, then we can proceed
by induction on $p$. Indeed we have a good understanding of the fibers of that forgetful map
(by dint of proposition~\ref{prop-elbowbdry}), and the case $p=0$, the induction beginning, is taken care of
by proposition~\ref{prop-elbow}. Therefore we concentrate on showing that~\eqref{eqn-funnyind1} is a weak
fibration. (This appears to be nontrivial even in the case where the Riemannian metric on $L$ is flat.)
Here we can again proceed by induction on the positive integer $k_p$. The case where the
map from $\uli k_p$ to $\uli k_{p-1}$ in $\sD$ is injective is trivial, i.e., in that case~\eqref{eqn-funnyind1} is
clearly a weak fibration. In particular that takes care of the case $k_p=1$. Suppose now that
the map from $\uli k_p$ to $\uli k_{p-1}$ is not
injective. Then $k_p\ge 2$ and without loss of generality, the elements $k_p$ and $k_p-1$ of $\uli k_p$
both map to the same element of $\uli k_{p-1}$. Let $\sD'$ be the diagram obtained from $\sD$ by deleting
the element $k_p$ from the set $\uli k_p$, and restricting the map $\uli k_p\to \uli k_{p-1}$
accordingly. By inductive assumption, the forgetful map
\begin{equation} \label{eqn-funnyind2} \Pack_1(\sD',c_0,c_1,\dots,c_{p-1},3c_p) \lra
\Pack(d_p\sD,c_0,c_1,\dots,c_{p-1}) \end{equation}
is a weak Serre fibration. (Note that $d_p\sD=d_p\sD'$.) Let
\[  W \subset \Pack_1(\sD,c_0,c_1,\dots,c_{p-1},c_p) \]
be the preimage of $\Pack_1(\sD',c_0,c_1,\dots,c_{p-1},3c_p)$ under the forgetful map
\[  \Pack_1(\sD,c_0,c_1,\dots,c_{p-1},c_p) \to \Pack_1(\sD',c_0,c_1,\dots,c_{p-1},c_p). \]
(Strictly speaking it is incorrect to say that $\Pack_1(\sD',c_0,c_1,\dots,c_{p-1},3c_p)$ is a subspace of
$\Pack_1(\sD',c_0,c_1,\dots,c_{p-1},c_p)$, but there is a preferred embedding of one
into the other given by fattening the innermost multipatches.)
By analogy with the proof of proposition~\ref{prop-elbow}, the map
\begin{equation} \label{eqn-funnyind3}  W \lra \Pack(d_p\sD,c_0,c_1,\dots,c_{p-1})  \end{equation}
obtained from~\eqref{eqn-funnyind1} by restriction is also a weak Serre fibration, since it is the
composition of the fiber bundle projection $W\to \Pack_1(\sD',c_0,c_1,\dots,c_{p-1},3c_p)$
and the weak Serre fibration~\eqref{eqn-funnyind2}. Now it is almost
true, though perhaps not quite true, that~\eqref{eqn-funnyind1} is a fiberwise
homotopy retract of~\eqref{eqn-funnyind3}. More precisely: the vector field argument from the proof
of proposition~\ref{prop-elbow} can be employed fiberwise to show that
if $K\subset \Pack_1(\sD,c_0,c_1,\dots,c_{p-1},c_p)$
is any compact subset, then there exists a vertical homotopy
\[  \big(h_t\co K \to \Pack_1(\sD,c_0,c_1,\dots,c_{p-1},c_p)\big)_{t\in [0,1]} \]
such that $h_0$ is the inclusion and $h_1(K)$ is contained in $W$. (\emph{Vertical} means that each
$h_t$ is a map over $\Pack(d_p\sD,c_0,c_1,\dots,c_{p-1})$.) Using that and the information that~\eqref{eqn-funnyind3}
is a weak Serre fibration, it is straightforward to deduce that~\eqref{eqn-funnyind1} is also a weak Serre
fibration. \qed

\begin{defn} \label{defn-wkSfib} A map $p\co E\to B$ of spaces is a \emph{weak Serre fibration} if the following holds.
For every compact CW-space $X$ with a map $g\co X\to E$ and a homotopy $(h_t\co X\to B)$ such that $h_0=pg$ and $h_t=h_0$
for $t$ in a neighborhood of $0\in [0,1]$, there exists a homotopy $(H_t\co X\to E)_{t\in [0,1]}$
such that $H_0=g$ and $pH_t=h_t$ for all $t\in [0,1]$. (This is a variant of the notion of \emph{weak fibration};
see \cite{DoldPofU}.)
\end{defn}

\begin{rem} \label{rem-elbow}
The formulation of theorem~\ref{thm-elbow} reflects the proof. An important guiding principle in that proof
is that the face operator $d_1\co \config(L)_1\to \config(L)_0$\,, also known as \emph{target}, has some good properties, such as being a Serre fibration.
The face operator $d_0$ from $\config(L)_1$ to $\config(L)_0$ cannot compete with that. Therefore,
to find positive numbers $c_0,c_1,\dots,c_p$ small enough in theorem~\ref{thm-elbow}, we determine first $c_0$ which sets the lower bound for patch radii in
$U_0$ of~\eqref{eqn-syspatch}, the ultimate target. The upper bound for allowed choices of $c_0$ will depend strongly on the metric properties of $L$, such as diameter,
curvature and the like. Then we select $c_1$. Together, $c_0$ and $c_1$ determine the lower bound $c_0c_1$ for patch radii in
$U_1$ of~\eqref{eqn-syspatch}. (Therefore $c_1$ has the purpose of a ratio.) The upper bound for the choice of $c_1$
depends less strongly on the metric properties of $L$, especially in the limit $c_0\to 0$, because it
only reflects the Riemannian metric properties of disks or open balls in $L$ of radius $\le c_0$ (which we may rescale to radius $1$ without losing
essential information). Then we select $c_2,c_3,\dots,c_p$ in the same manner. The upper bounds
for $c_1,c_2,\dots,c_p$ may depend on $c_0$ but they converge to the same positive real number for $c_0\to 0$, and that positive
real number does not even depend on $L$.
\end{rem}

\begin{cor} \label{cor-elbow} For every choice of positive integers $r$ and $\alpha$
there exist a positive real number $\rho$ and a simplicial subspace $X^\fat$ of the $r$-skeleton $\skel_r\cconfig(L;\alpha)$
such that every patch in a multipatch $V\in X^\fat_0$ has radius $\ge\rho$, and the inclusion of $X^\fat$ in   $\skel_r\cconfig(L;\alpha)$ is a weak equivalence.
\end{cor}
\proof By theorem~\ref{thm-elbow}, for every diagram
\[
\sD=\big(\uli k_0 \xleftarrow{g_1} \uli k_1\xleftarrow{g_2} \cdots\xleftarrow{g_{p-1}} \uli k_{p-1}
\xleftarrow{g_p} \uli k_p\big)~~\in (N\fin)_p~,
\]
positive real numbers $b_0,...,b_p$ can be selected, all $\le 1$ with the possible exception of $b_0$\,, such that
the inclusion of $\Pack(\sD,c_0,\dots,c_p)$ in $\config(L)_\sD$ is a homotopy equivalence
whenever $0<c_j\le b_j$ for all $j\in\{0,1,\dots,p\}$. The selection can be made in such
a way that $b_j=1$ whenever $g_j$ is injective.
It is a trivial matter to translate this into a statement about the unordered configuration category,
but it is convenient to keep the diagram $\sD$. Therefore we write $\Pack([\sD],c_0,\dots,c_p)$
and $\cconfig(L)_{[\sD]}$. Let $b$ be the minimum of all the $b_j$ obtained in this way,
for all $\sD\in (N\fin)_p$ involving only objects of $\fin$ of cardinality $\le \alpha$,
and all $p\le r$. This is still a positive real number. By definition, the
preimage of $[\sD]$ in $X^\fat_p$ is $\Pack([\sD],c_0,\dots,c_p)$,
where $c_0,c_1,\dots,c_p$ depend on $\sD$ and $r$ as follows:
\begin{itemize}
\item[] $c_0c_1\cdots c_p=b^{1+r}$;
\item[] $c_j=1$ if $j>0$ and $g_j$ is injective,
\item[] $c_j=b$ if $j>0$ and $g_j$ is not injective.
\end{itemize}
This defines $X^\fat$ in degrees $p\le r$. By inspection,
the simplicial operators in $\skel_r\cconfig(L;\alpha)$ respect $X^\fat$ as long as they are induced by
morphisms $[p]\to [q]$ in $\Delta$ where $p,q\le r$. There is exactly one way to finish the
construction of $X^\fat$ in such a way that it is a simplicial subspace of $\skel_r\cconfig(L;\alpha)$.
The inclusion of $X^\fat$ in $\skel_r\cconfig(L;\alpha)$ is a degreewise weak equivalence. This is
true by theorem~\ref{thm-elbow} in degrees $\le r$. In degrees $>r$ it is an easy consequence
of the statement for degrees $\le r$ and the observation that both $X^\fat$ and $\skel_r\cconfig(L;\alpha)$
are Reedy cofibrant simplicial spaces. --- Therefore $\rho:=b^{1+r}$ is a correct decision. \qed

\section{Perturbation and tolerance} \label{sec-pseudomet}
\emph{Vocabulary and notation}: see section~\ref{sec-models}.

\smallskip
In this section we assume that $M$ is a Riemannian manifold and
$L$ is a compact smooth submanifold of $M$ (both without boundary). Then $L$ inherits a Riemannian metric from $M$, and we assume
that $L$ is totally geodesic in $M$ (all geodesics in $L$ are also geodesics in $M$) and that $M$ is geodesically
complete. For simplicity we also assume that $M$ is connected. Then the Riemannian metric on $M$ induces an honest
metric $d_M$ on $M$, the geodesic distance.

We also want to speak of a preferred inclusion $\iota$
of $\cconfig(L)$ in $\cconfig(M)$. This should be ``induced'' by the inclusion $L\to M$.
On objects, the idea is therefore that we take a multipatch in $L$
with center points $z_1,z_2,\dots,z_k$ and radii $c_1,c_2,\dots,c_k$ to the multipatch in $M$ with the
same center points $z_1,z_2,\dots,z_k$ and the same radii. But this may not work. We are on the safe
side if we assume that for every $z\in L$ the injectivity radius of $M$ at $z$ agrees with the injectivity
radius of $L$ at $z$. This is quite a strong condition on the inclusion $L\hookrightarrow M$.

\smallskip
\begin{defn} \label{defn-epthick} Let $U$ be an object in $\cconfig(L)$, a multipatch in $L$. For $\vep>0$ we denote by $\Theta_\vep U\subset M$
the $\vep$-thickening of $\iota(U)$ in $M$, in other words, the union of all metric open balls in $M$
with center in $\iota(U)$ and radius $\vep$. Nota bene: we assumed $U\subset L$, but $\Theta_\vep U$
is an open subset of $M$.

More generally, for an open subset $V$ of $L$, we sometimes write $\Theta_\vep V$ for the union of all $\Theta_\vep U$
where $U$ is an object of $\cconfig(V)$. This is also an open subset of $M$.
\end{defn}

\begin{defn} \label{defn-eprcontrol} Let $\vep$ and $\delta$ be positive real numbers.
A map from $\cconfig_\delta(L;\alpha)$ to $\varphi\cconfig(M)$
over $\cfin$ is \emph{$\vep$-bounded} if
for every object $U$ in $\cconfig_\delta(L;\alpha)$, the map takes $\cconfig_\delta(U;\alpha)$ to $\varphi\cconfig(\Theta_\vep U)$.

More generally, let $f$ be a $k$-simplex in
$\map(\cconfig_\delta(L;\alpha),\varphi\cconfig(M))$. This can be thought of as a family
of maps $f_z\co \cconfig_\delta(L;\alpha)\to \varphi\cconfig(M)$ depending continuously on
$z\in \Delta^k$. Consequently we say that $f$ is \emph{$\vep$-bounded} if each $f_z$
is $\vep$-bounded.
\end{defn}

Let $Z_\vep\subset \map(\cconfig_\delta(L;\alpha),\varphi\cconfig(M))$
be the simplicial subset consisting of the $\vep$-bounded simplices.
\begin{thm} \label{thm-shrink} For every choice of $\alpha$ there exist $\vep,\delta>0$
such that the inclusion $Z_\vep\to \map(\cconfig_\delta(L;\alpha),\varphi\cconfig(M))$
is based nullhomotopic.
\end{thm}
\emph{Clarification.} The real number $\vep$ depends on $L$ and $M$ with their Riemannian metrics,
and on $\alpha$.
Meaning of \emph{the inclusion is based nullhomotopic}: the inclusion is based homotopic to the constant map
taking everything to the base point $\iota$.

\bigskip
\emph{Proof of theorem~\ref{thm-shrink}.} The finite presentation result
of \cite{confpres} says that for large enough $r$, the $r$-skeleton $\skel_r\config(L;\alpha)$ of $\cconfig(L;\alpha)$
contains a finite presentation of $\cconfig(L;\alpha)$. Therefore it suffices to show: for every $r\ge 1$ there exist
$\vep,\delta>0$ such that the composition
\[ Z_\vep\hookrightarrow  \map_\cfin(\cconfig_\delta(L;\alpha),\varphi\cconfig(M)) \xrightarrow{\textup{res.}}
\map_\cfin(\skel_r\cconfig_\delta(L;\alpha),\varphi\cconfig(M)) \]
is nullhomotopic.

Here is an overview of the rest of the proof. Let us write $X^\lean$ for $\skel_r\cconfig_\delta(L;\alpha)$,
simplicial subspace of $\skel_r\cconfig(L;\alpha)$. By imposing certain
\emph{lower} bounds on the radii of patches in multipatches, we define another
simplicial subspace $X^\fat$ of $\skel_r\cconfig(L;\alpha)$; this was done already in corollary~\ref{cor-elbow}.
The inclusion of $X^\fat$ in $\skel_r\cconfig(L;\alpha)$ is a degreewise weak equivalence by construction.
Also by construction, $X^\lean$ and $X^\fat$ are degreewise disjoint simplicial subspaces of
$\skel_r\cconfig(L;\alpha)$. Another simplicial subspace $X$ of $\skel_r\cconfig(L;\alpha)$ will be
defined which contains the (disjoint) union of $X^\lean$ and $X^\fat$, and then
some more to ensure that the inclusion $X\to \skel_r\cconfig(L;\alpha)$
is a good approximation. It turns out to be a conservatization map.
Then we construct the broken arrow in a strictly commutative diagram
\begin{equation} \label{eqn-guidethm2}
\begin{aligned}
\xymatrix@R=30pt@M=8pt@C=50pt{
&  {\map_\cfin(X^\fat,\varphi(\cconfig(M)))} \\
Z_{\vep,r} \ar@{..>}[r] \ar@/^1pc/[dr]_-{\textup{forgetful }\,} \ar@/_1pc/[ur]^-{\textup{trivial map\quad}}
&  {\map_\cfin(X,\varphi(\cconfig(M)))} \ar[d]^-{\textup{restriction}} \ar[u]_-{\textup{restriction}} \\
&  {\map_\cfin(X^\lean,\varphi(\cconfig(M)))}
}
\end{aligned}
\end{equation}
The two vertical arrows are weak equivalences since the inclusion maps
$X^\fat\to X$ and $X^\lean\to X$ are weak equivalences modulo conservatization,
and the common target $\varphi(\cconfig(M))$ is conservative. This completes the argument. It
does not construct an explicit nullhomotopy. Instead it proceeds by
making a few useful homotopy equivalences. (\emph{End of overview.})

First order of business: definition of $X^\fat$ and the choice of $\vep,\delta>0$. For $X^\fat$ we can take the
construction in the proof of corollary~\ref{cor-elbow}.
The only condition on $\vep$ is $\vep<\rho/2$, where $\rho$ is the lower bound
in corollary~\ref{cor-elbow}, also known as $b^{1+r}$. We choose $\delta$ so that $\delta< \vep/2$.

The next order of business is the definition of $X$, again a simplicial subspace of
$\skel_r\cconfig(L;\alpha)$. For motivation and warm-up, let us take
the view that $\skel_r\cconfig(L;\alpha)$ is a Segal space \emph{in degrees $\le r$}, fiberwise
complete over the nerve of $\cfin$. (The prefix $\skel_r$ has done some damage to
the Segal property.)
Consequently $X^\fat$ and $X^\lean$ are both
Segal spaces \emph{in degrees $\le r$}, fiberwise
complete over the nerve of $\cfin$. (The inclusions $X^\fat\to \skel_r\cconfig(L;\alpha)$
and $X^\lean\to \skel_r\cconfig(L;\alpha)$ are degreewise weak equivalences.) As noted
in the overview, $X^\fat$ and $X^\lean$ are degreewise disjoint as simplicial subspaces of
$\skel_r\cconfig(L;\alpha)$ (because multipatches in $L$ which are objects of $X^\fat$ have all patch radii
$>2\vep$, whereas multipatches which are objects of $X^\lean$ have all patch radii $<\vep/2$).
We wish to define $X$ as a simplicial subspace of $\skel_r\cconfig(L;\alpha)$,
in such a way that it is also a Segal
space in degrees $\le r$. It must contain $X^\fat$ and $X^\lean$. It should not contain objects other
than those in $X^\fat$ and $X^\lean$.
It should however contain ``most'' of the morphisms in $\cconfig(L;\alpha)$ from objects in $X^\lean$ to objects in $X^\fat$.
(There are no morphisms in $\cconfig(L;\alpha)$  from objects in $X^\fat$ to objects in $X^\lean$.) To be precise: a morphism
in $\cconfig(L;\alpha)$ from an object in $X^\lean$ to an object in $X^\fat$ (which is an inclusion
$U\to V$ of multipatches) qualifies as a
morphism in $X$ if the closure of $U$ in $V$ has distance $>2\vep$ from the boundary
of $V$ in $L$. Therefore $X_n$ for $n\le r$ is the space of systems of multipatches
\begin{equation} \label{eqn-seq0}  U_0\supset U_1 \supset \dots \supset U_{n-1} \supset U_n \end{equation}
in $L$ satisfying one of the following (mutually exclusive) conditions (i),(ii),(iii).
\begin{itemize}
\item[(i)] The entire system is an element in degree $n$ of $X^\fat$.
\item[(ii)] The entire system is an element in degree $n$ of $X^\lean$.
\item[(iii)] There exists $p\in \{0,1,2,\dots,n\}$ such that
$U_{p+1}\supset U_{p+2} \supset \dots \supset U_{n-1} \supset U_n$
is an element of $X^\lean$ in degree $n-p-1$ and
$U_0\supset U_1 \supset \dots \supset U_{p-1} \supset U_p$
is an element of $X^\fat$ in degree $p$, and the closure of $U_{p+1}$ in $U_p$ has
distance $>2\vep$ from the boundary
of $U_p$ in $L$.
\end{itemize}
We give almost the same definition for $p>r$, but then we require in addition
that $U_0\supset U_1 \supset \dots \supset U_{n-1} \supset U_n$ is an element of
$\skel_r\cconfig(L;\alpha)$ in degree $p$. Equivalently, at least $n-r$ of the
containments $U_{j-1}\supset U_j$ are identities.

By inspection, $X$ is a Segal space in degree $\le r$. As such it is also fiberwise complete 
over the nerve of $\cfin$. But it is not conservative over the nerve of $\cfin$ (except perhaps
for some extreme choices of $L$, such as $L=\emptyset$). Indeed, there are elements in $X_1$
of the form $U_0\supset U_1$
where $U_0$ is an element of $X^\fat$ in degree $0$ and $U_1$ is an element of $X^\lean$ in degree $0$,
and the inclusion induces a bijection $\pi_0(U_1)\to \pi_0(U_0)$. These do not qualify as
homotopy invertible elements in $X_1$.

The Segal space $X$ comes with an important simplicial map $w$ to $\Delta[1]$, the nerve of $[1]^\op$. This expresses the
fact that objects of $X$ belong either
to $X^\fat$, in which case $w$ takes them to $0$, or to $X^\lean$, in which case $w$ takes them to $1$.

Next order of business: construction of the broken arrow in diagram~\eqref{eqn-guidethm2}. Readers are
advised to try a construction of their own making, but if none is forthcoming then they should read on.
Let $f$ be a $0$-simplex of $Z_{\vep,r}$. (We ought to begin with: \emph{let $f$ be a $j$-simplex
of $Z_{\vep,r}$}, but we take $j=0$ and leave the
generalization to arbitrary $j\ge 0$ to the reader. Remember that a $j$-simplex of
$Z_{\vep,r}$ is nothing but a family of $0$-simplices of $Z_{\vep,r}$
parametrized continuously by $\Delta^k$.) Then $f$ restricts to a map from $X^\lean$ to $\varphi\cconfig(M)$, which we still denote by $f$.
We must extend this to a map $f^!$ defined on $X$ in such a way that $f^!$ also extends the
standard inclusion
$\iota\co \skel_r\cconfig(L;\alpha)\to \cconfig(M)$ restricted to $X^\fat$.
In the language and notation of section~\ref{sec-joins}, the construction of $f^!$ can be indicated
by the following diagram. Write $N\sP$ for
$\cconfig(L;\alpha)$ and $N\sQ$ for $\cconfig(M)$. These are nerves of topological posets $\sP$ resp.~$\sQ$,
because we are using the Riemannian patch models for the configuration categories. The goal
is now to define $f^!\co X\to \varphi N\sQ$.

\begin{equation} \label{eqn-joinshow}
\begin{aligned}
\xymatrix@R=20pt@C=30pt@M=5pt{
X \ar[d]_-{\textup{(incl.,$w$)}} \ar[d] \ar@{..>}@/_1pc/[drr] && \varphi N\sQ  \\
N\sP\times\Delta[1] \ar[d]_-{\textup{sep.\,diag., defn.~\ref{defn-fracdi}}} &&
\varphi(N\sQ\times\Delta[1])  \ar[d]^-{\varphi(\textup{sep.\,diag.})}
\ar[u]_-{\varphi(\textup{proj.})}   \\
N\sP*N\sP  \ar[r]^{\iota * f} & N\sQ*\varphi N\sQ  \ar[r]^-{\textup{defn.~\ref{defn-phijoin}}}  & \varphi(N\sQ* N\sQ)
}
\end{aligned}
\end{equation}
The vertical arrow in the right-hand column with the label $\varphi(\textup{sep.\,diag.})$ is an embedding. It follows that the broken arrow
is unique if it exists. We have to show that it exists, i.e., that the image
of the map $X\to \varphi(N\sQ*N\sQ)$ obtained by going along the lower circuit is contained in the
image of the map with the label $\varphi(\textup{sep.\,diag.})$.
We begin with an $n\in\{0,1,\dots,r\}$ and an element~\eqref{eqn-seq0} of $X_n$. In the cases (i) and (ii) of~\eqref{eqn-seq0},
there is nothing to verify, so we can assume (iii). Moving forward to $N\sQ*\varphi N\sQ$ in~\eqref{eqn-joinshow},
we obtain the element
\[   \big(\iota(U_0\supset ... \supset U_p),\,f(U_{p+1}\supset \dots\supset U_n)\big) \in N\sQ_p\times (\varphi N\sQ)_{n-p-1}
\subset (N\sQ* \varphi N\sQ)_n\,. \]
The important observation to make now is that
\begin{equation} \label{eqn-impobs} f(U_{p+1}\supset \dots\supset U_n) \in (\varphi\cconfig(\Theta_\varepsilon U_{p+1}))_{n-p-1}\subset
(\varphi N\sQ)_{n-p-1} \end{equation}
since $f$ is $\varepsilon$-bounded, and by (iii), we have
\[
\Theta_\varepsilon U_{p+1} \subset \iota(U_p).
\]
(Indeed, if $z \in \Theta_\varepsilon U_{p+1}$ then $z$ has distance $< \varepsilon + \delta$ from some center point of $U_{p+1}$,
but that center point has distance $\geq 2 \varepsilon + \delta$ from the complement of $\iota(U_p)$.
Therefore $z$ is not in the complement of $\iota(U_p)$.) 
Now we can apply lemma~\ref{lem-asinus} and it follows that
the broken arrow in~\eqref{eqn-joinshow} exists.

There is something tedious left to do under the same heading. The extension
$f^!$ of $f$ (restricted to $\skel_r\cconfig_\delta(L;\alpha)$) has been defined, but it remains to be shown that it is a map over $\cfin$.
We adopt the point of view developed at the end of section~\ref{sec-models}, around lemma~\ref{lem-classy}.
Therefore, to turn $f^!$ into a map over $\cfin$ we require, for every
element $V$ of $X_0$ (which is a multipatch in $L$), a preferred bijection
\[  \omega_V\co  \pi_0(V) \to \pi_0(f^!(V)). \]
We have this already because $V$ belongs either to $X^\lean$, in which case $\omega_V$
comes with the package $f$, or $V$ belongs to $X^\fat$, in which case $f^!(V)=\iota(V)$ and we define $\omega_V$ to be the
bijection induced by the inclusion of $V$ in $\iota(V)$. It remains only to test for naturality.
Suppose that $U$ and $V$ are elements of $X_0$ (multipatches in $L$) such that $U\subset V$. Then we must show that the diagram
of finite sets
\begin{equation} \label{eqn-ufincheck}
\begin{aligned}
\xymatrix{
\pi_0(U) \ar[d]_-{\textup{ind. by incl.}} \ar[r]^-{\omega_U} & \pi_0(f^!(U)) \ar[d]^-{\textup{determined by $f^!(V\supset U)$}} \\
\pi_0(V) \ar[r]^-{\omega_V} & \pi_0(f^!(V))
}
\end{aligned}
\end{equation}
commutes. Since $f^!$ is a derived thing, we must not take for granted that
$f^!(U)$ is contained in $f^!(V)$. The cases where $U$ and $V$ are both in $X^\lean$ or both in $X^\fat$ are
nevertheless trivi\-al.
Therefore we may assume that $U$ belongs to $X^\lean_0$ and $V$ belongs to $X^\fat_0$. Then~\eqref{eqn-ufincheck} simplifies to
\begin{equation} \label{eqn-ufincheck2}
\begin{aligned}
\xymatrix@C=50pt{
\pi_0(U) \ar[d]_-{\textup{ind. by incl.}} \ar[r]^-{\omega_U} & \pi_0(f(U)) \ar[d]^-{\textup{det. by $f^!(V\supset U)$}} \\
\pi_0(V) \ar[r]^-{\textup{ind. by incl.}} & \pi_0(\iota(V))
}
\end{aligned}
\end{equation}
As stated in section~\ref{sec-models}, after lemma~\ref{lem-itsenough}, we can reduce to the situation where $\pi_0(U)$ has only
one element. Assuming this, let $V_0$ be the connected component of $V$ containing $U$.
Then $(V_0\supset U)$ is an element of $X$ in degree $1$, and $(V\supset V_0\supset U)$
is an element of $X$ in degree $2$, and so diagram~\eqref{eqn-ufincheck2} can be expanded
as follows:
\begin{equation} \label{eqn-ufincheck4}
\begin{aligned}
\xymatrix@R=30pt@C=60pt{
\pi_0(U) \ar[d]_-{\textup{ind. by incl.}} \ar[r]^-{\omega_U} & \pi_0(f(U)) \ar[d]_-{\textup{det.~by $f^!(V_0\supset U)$}}
\ar@/^4pc/[dd]^-{\textup{det.~by $f^!(V\supset U)$}} \\
\pi_0(V_0) \ar[d]_-{\textup{ind. by incl.}} \ar[r]^-{\textup{ind. by incl.}} & \pi_0(\iota(V_0)) \ar[d]_-{\textup{ind. by incl.}} \\
\pi_0(V) \ar[r]^-{\textup{ind. by incl.}} & \pi_0(\iota(V))
}
\end{aligned}
\end{equation}
Here the upper square commutes because $\pi_0(\iota(V_0))$ has only one element, and the lower square commutes
because $\pi_0$ is a functor.
The deformed triangle commutes because we have the element $f^!(V\supset V_0\supset U)$ in degree 2 of $\cconfig(M)$.
So~\eqref{eqn-ufincheck4} commutes, and therefore~\eqref{eqn-ufincheck} commutes.

Next and last order of business in this proof: showing that the inclusion of $X$ in $Y:=\skel_r\cconfig(L;\alpha)$
is a conservatization map. There is a good reason why this should be so: the inclusion $X\hookrightarrow Y$
has a factorization
\begin{equation} \label{eqn-factoconserv}  X \xrightarrow{(\textup{incl.},w)}  Y\times \Delta[1] \xrightarrow{\textup{proj.}} Y  \end{equation}
where the first arrow is a weak equivalence (by inspection). Here we view $Y$ as a
simplicial space with reference map $Y\to \cfin$, as usual, and we view $\Delta[1]$ as a simplicial space with reference map
$v\co \Delta[1]\to \Delta[0]$. Therefore $Y\times \Delta[1]$ is a simplicial
space with reference map to $\cfin\times\Delta[0]\cong \cfin$. In this situation proposition~\ref{prop-prodconserv} states that
the conservatization $\Lambda$ respects products. (We use a new incarnation $\Lambda^\fre$ of the conservation procedure here,
but we write $\Lambda$ for short. See section~\ref{sec-fibreal}.) Therefore the map
\[ \id_Y\times v \co Y\times \Delta[1] \to Y \times \Delta[0]  \]
is a conservatization map over $\cfin\times\Delta[0]\cong \cfin$ because $\id\co Y\to Y$ is a conservatization map over $\cfin$
and $\Delta[1]\to \Delta[0]$ is a conservatization map over $\Delta[0]$. This completes the verification that
$X\to Y$ is a conservatization map, and thereby the proof as a whole. \qed

\section{Partitions of unity in a simplicial setting} \label{sec-pusimp}
Let $Y$ be a space and let $\sU=(U_i)_{i\in J}$ be an open cover of $Y$.
For every finite nonempty subset $S$ of $J$, write $U_S=\bigcap_{i\in S} U_i$.
Then $S\mapsto U_S$ is a contravariant functor from the poset of nonempty subsets of $J$ to
spaces. There is a standard projection map
\begin{equation} \label{eqn-Segalproj}
\hocolimsub{\textup{fte nonempty }S\subset J} U_S\quad  \lra\quad Y~=~\colimsub{\textup{fte nonempty }S\subset J} U_S\,.
\end{equation}

\begin{lem} \label{lem-Segalproj} \emph{(Segal \cite[\S4]{Segal1968}.)} If $Y$ is paracompact, then~\eqref{eqn-Segalproj} is a homotopy
equivalence. \qed
\end{lem}
We repeat some of Segal's arguments because we want to make a statement which is stronger in some respects.
The official definition of the hocolim in~\eqref{eqn-Segalproj} is
\[   {\coprod_{S} A_S \times U_S}\Big/\textup{relations}  \]
where $S$ is a nonempty finite subset of $J$ and $A_S$ is the classifying space of the
poset of nonempty subsets of $S$. (The relations in the denominator of the ``fraction'' are of the coend type; indeed,
$S\mapsto A_S$ is a covariant functor and $S\mapsto U_S$ is contravariant.) There is a natural identification
of $A_S$ with $\Delta(S)$, the simplex spanned by $S$. Therefore the source in~\eqref{eqn-Segalproj} can be
written as
\[   {\coprod_{S} \Delta(S) \times U_S}\Big/\textup{relations}.  \]
A partition of unity $(\psi_i\co Y\to [0,1])_{i\in J}$ subordinate to the
cover $\sU$ gives rise to a map from $Y$ to $\hocolim_S~U_S$ which takes $x\in Y$ to the element
represented by
\[ ((\psi_i(x))_{i\in T},x)\in \Delta(T)\times U_T  \]
where $T$ is the finite nonempty set $\{i\in J~|~\psi_i(x)>0\}$. We are using barycentric coordinates in
$\Delta(T)$. It is easily seen that this map from
$Y$ to $\hocolim_S~U_S$ is
continuous, and that it defines a section of~\eqref{eqn-Segalproj}. This is our cue for making a
stronger statement.

\begin{lem} \label{lem-SegalPUext} The space of sections of \eqref{eqn-Segalproj}
can be identified with the space $\sP(Y,\sU)$ of partitions of unity subordinate to $\sU$.
\end{lem}
Here it is necessary to explain what is meant by \emph{space of sections} and the like.
In the absence of an obvious or otherwise preferred topology on the set of such sections, we define
the space of sections as a simplicial set. A $k$-simplex is a map from
$\Delta^k\times Y$ to $\hocolim_S U_S$ such that (post-)composition with~\eqref{eqn-Segalproj} gives the
projection from $\Delta^k\times Y$ to $Y$. Similarly, $\sP(Y,\sU)$ is a simplicial set. A $k$-simplex
is a partition of unity subordinate to the open cover $(\Delta^k\times U_i)_{i\in J}$ of $\Delta^k\times Y$.

As regards the proof of lemma~\ref{lem-SegalPUext}, the arguments already given can be extended mechanically to give us a map
from $\sP(Y,\sU)$ to the space (simplicial set) of sections of~\eqref{eqn-Segalproj}. This is
injective by inspection. Surprisingly, it is in fact bijective. See \cite[Lem. 2.3]{WeissWhatclass}.
(Keep in mind that the hocolim in~\eqref{eqn-Segalproj} was \emph{not} defined to be a subspace of
$\Delta(J)\times Y$, not even in the cases where $J$ is finite. It has a continuous and injective map
to $\Delta(J)\times Y$, but often this is not a topological embedding. Example:
take $Y=[0,1]$, let $J=\{0,1\}$ and let $\sU$ be the open cover consisting of $U_0=Y\smin \{0\}$
and $U_1=Y\smin\{1\}$. Then the $\hocolim$ is the homotopy pushout of the diagram
$U_0\leftarrow U_0\cap U_1\rightarrow U_1$.
It is not metrizable and is consequently not homeomorphic to any subspace of $\Delta^1\times [0,1]$
whatsoever.)

The following two lemmas are in the nature of observations.

\begin{lem} \label{lem-PU1} If $Y$ is a paracompact space and $\sU=(U_i)_{i\in J}$ is an open cover of $Y$,
then $\sP(Y,\sU)$ is a contractible simplicial set. It is also fibrant, i.e., it has the Kan filling property. \qed
\end{lem}

\begin{lem} \label{lem-PU2} If $Y$ is a paracompact space, $\sU=(U_i)_{i\in J}$ is an open cover of $Y$
and $C$ is a closed subspace of $Y$, then the forgetful map of simplicial sets
\[  \sP(Y,\sU) \lra \sP(C,\sU_C)  \]
\emph{(where $\sU_C:=(U_i\cap C)_{i\in J}$)} is a Kan fibration. \qed
\end{lem}

\bigskip
Next we will develop simplicial variants of lemmas~\ref{lem-Segalproj} and~\ref{lem-SegalPUext}.
So let $Y$ be a \emph{simplicial} space. Suppose that $Y_0$ is equipped with an
open cover $\sU=(U_i)_{i\in J}$. For $j\in [k]$ let $e_{k,j}\co Y_k\to Y_0$ be the face operator
determined by the morphism $[0]\to [k]$ in $\Delta$ taking $0$ to $j$.  We write
\begin{itemize}
\item[] $U(k,j)_i$ for the preimage of $U_i$ under $e_{k,j}$
(an open subset of $Y_k$);
\item[] similarly $U(k,j)_S$ for the preimage of $U_S$ (where $S\subset J$ is finite, nonempty);
\item[] $\sU(k,j)$ for the open cover of $Y_k$ determined by $\sU$, open cover of $Y_0$, and pullback along $e_{k,j}$.
\end{itemize}
We make the following assumptions.
\begin{itemize}
\item[-] $Y$ is degreewise paracompact Hausdorff and Reedy cofibrant.
\item[-] For each $i\in J$, we have the following inclusions of open subsets of $Y_k$:
\[ U(k,0)_i\subset U(k,1)_i\subset \dots \subset U(k,k)_i\,. \]
\end{itemize}
Let $S$ be a finite nonempty
subset of $J$. Then $[k] \mapsto U(k,0)_S$
defines a simplicial subspace of $Y$.

\begin{prop} \label{prop-dersecabstract} There exists a factorization as in
\begin{equation} \label{eqn-factosol}
\begin{aligned}
\xymatrix@R=24pt@C=0pt@M=6pt{
 &   {\hocolimsub{\textup{fte nonempty }S\subset J} (\,[k]\mapsto U(k,0)_S\,)} \ar[d]^-\simeq \\
{\varphi^\lad Y} \ar@{..>}[ur] \ar[r]^-\simeq &  Y
}
\end{aligned}
\end{equation}
where the horizontal arrow is adjoint to the inclusion of $Y$ in $\varphi Y$.
\end{prop}
(The hocolim is a degreewise hocolim of simplicial spaces. The vertical arrow is a
levelwise homotopy equivalence by lemma~\ref{lem-Segalproj}. The factorization is strict.)

\proof For this undertaking we need:
\begin{itemize}
\item[(i)] for every $p\ge 0$ and every string of $p$ composable morphisms
\[  \sD: \qquad [k_0] \xleftarrow{g_1} [k_1] \xleftarrow{g_2} [k_2]  \xleftarrow{g_3} \cdots \xleftarrow{g_p} [k_p] \]
in $\Delta_\inj$, a $p$-simplex $\sigma(\sD)$ in $\sP(Y_{k_0},\sU(k_0,a))$ where $a=g_1g_2\cdots g_p(0)$.
\end{itemize}
The simplices $\sigma(\sD)$ are jointly subject to the following condition:
\begin{itemize}
\item[(ii)] $f^*(\sigma(\sD))$ is in agreement with $\sigma(f^*\sD)$
whenever $f\co [p_1]\to [p_0]$ is a morphism in $\Delta$ and $\sD$ is a string of $p_0$ composable morphisms
in $\Delta_\inj$.
\end{itemize}
The agreement asked for is an equality of $p_1$-simplices in the simplicial set
which is the home of $\sigma(\sD)$ and $f^*(\sigma(\sD))$, but $\sigma(f^*\sD)$ is a $p_1$-simplex
in
\[ \sP(Y_{k_{f(0)}},\,\sU(k_{f(0)},b)) \]
where $b=g_{f(0)+1}g_{f(0)+2}\cdots g_{f(p_1)}(0)$.
The request is meaningful because there is a simplicial map
\[  \sP(Y_{k_{f(0)}},\,\sU(k_{f(0)},b)) \lra \sP(Y_{k_0},\,\sU(k_0,a)) \]
by pullback along the face operator $Y_{k_0}\to Y_{k_{f(0)}}$
determined by the morphism
\[ g_1g_2\dots g_{f(0)}\co [k_{f(0)}]\to [k_0] \]
in $\Delta_\inj$. This uses also that $a\ge g_1g_2\dots g_{f(0)}(b)$.

We have to impose one more condition.
Let $\sD$ be a diagram in $\Delta_\inj$ as before and let $s\co [k_0]\twoheadrightarrow [\ell_0]$
be a \emph{surjective} morphism in $\Delta$. Then diagram $\sD$ can be uniquely completed to a commutative diagram
of the form
\[
\xymatrix@M=7pt@R=18pt{
\sD: & [k_0] \ar@{->>}[d]^-s & \ar[l]^{g_1} [k_1] \ar@{->>}[d] & \ar[l]^-{g_2} [k_2] \ar@{->>}[d] & \ar[l]^-{g_3} \cdots &
\ar[l]^-{g_p}  [k_p] \ar@{->>}[d] \\
\sE: & [\ell_0] & \ar[l]^{h_1} [\ell_1] & \ar[l]^-{h_2} [\ell_2] & \ar[l]^-{h_3} \cdots & \ar[l]^-{h_p}  [\ell_p]
}
\]
in $\Delta$, where the arrows in the lower row (and in the upper row) belong to $\Delta_\inj$.
Now the new condition:
\begin{itemize}
\item[(iii)] $\sigma(\sD)$ is taken to $\sigma(\sE)$ under the map
\[ \sP(Y_{k_0},\,\sU(k_0,a)) \to \sP(Y_{\ell_0},\,\sU(\ell_0,s(a))) \]
induced by the degeneracy operator $s^*\co Y_{\ell_0} \to Y_{k_0}$.
\end{itemize}
Condition (ii) refers only to the face operators in $Y$, whereas
condition (iii) refers mainly to the degeneracy operators in $Y$.

We begin with the constructing of a collection $(\sigma(\sD))$
satisfying (i), (ii) and (iii). It will proceed by induction on $(p,k_0)$, where $p$ is the length of $\sD$
and $[k_0]$ is the ultimate target. We use the lexicographic
ordering on the set of such pairs $(p,k_0)$. The induction therefore begins with $p=0$, $k_0=0$
and continues with $p=0$, $k_0=1,2,3,...$. The induction steps here (where $p=0$ is fixed)
grapple with condition (iii) only. Lemma~\ref{lem-PU2} makes a contribution.
Then we make the step from $p=0$ to $p=1$, and this is
where condition (ii) becomes important.

\smallskip
\emph{Cases $p=0$, $k:=k_0$ arbitrary.} Here we have to select a $0$-simplex in
\[ \sP(Y_k,\,\sU(k,0)), \]
in other words a partition of unity on $Y_k$ subordinate to the cover $\sU(k,0)$. We proceed by induction
on $k$. For $k=0$, there is no further condition to be satisfied and a solution exists by lemma~\ref{lem-PU1}. For $k>0$,
the partition of unity is already prescribed (by inductive assumption) on the degenerate part of $Y_k$, also
known as $\latch_kY\subset Y_k$. The induction step can be carried out by lemma~\ref{lem-PU2}.

\smallskip
\emph{Cases $p=1$, $k_0$ arbitrary.} The diagram is $\sD=[k_0]\xleftarrow{g_1}[k_1]$, and we write $a=g_1(0)$.
We have to select a $1$-simplex $\sigma(\sD)$ in
\[ \sP(Y_{k_0},\,\sU(k_0,a)), \]
in other words a partition of unity on $\Delta^1\times Y_{k_0}$ subordinate to the open
cover $\Delta^1\times \sU(k_0,a)$.
If $g_1$ is onto, then it is the identity map of $[k_0]$ and $\sigma(\sD)$ is already fully prescribed by condition (ii);
it has to be the pullback of $\sigma([k_0])$ along the projection
\[ \Delta^1\times Y_{k_0}\to Y_{k_0}. \]
In the more interesting case where $g_1$ is not onto, we proceed by
induction on $k_0$. The induction begins with $k_0=1$. In that case the restriction of $\sigma(\sD)$ to
$\partial\Delta^1\times Y_1$ is prescribed because of condition (ii). Here we are also using the
standing assumption that $\Delta^1\times \sU(k_0,a)$ is refined by $\Delta^1\times \sU(k_0,0)$.
The restriction to $\Delta^1\times \latch_1Y$
is also prescribed because of condition (iii). There is a solution by lemmas~\ref{lem-PU1} and~\ref{lem-PU2}.
For $k:=k_0>1$, the restriction of $\sigma(\sD)$ to
$\partial\Delta^1\times Y_k$ is prescribed because of condition (ii).
Here we are also using the
standing assumption that $\Delta^1\times \sU(k_0,a)$ is refined by $\Delta^1\times \sU(k_0,b)$ whenever $b<a$.
The restriction to $\Delta^1\times \latch_kY$ is also prescribed because of the inductive assumption and condition (iii).
There is a solution by lemmas~\ref{lem-PU1} and~\ref{lem-PU2}.

\emph{Cases $p>1$, $k_0$ arbitrary.} If one of the maps $g_j$ in the diagram $\sD$ is onto, hence an identity
map, then $\sigma(\sD)$ is already determined thanks to condition (ii) and previous steps (smaller $p$).
Otherwise, the reasoning is as in the case $p=1$.

\smallskip
It remains to be said how the $\sigma(\sD)$ allow us to produce a broken arrow as in~\eqref{eqn-factosol}.
By equation~\eqref{eqn-adphi}, the simplicial space $\varphi^\lad Y$ is a quotient of $\coprod_n Y_n\times \varphi^\lad(\Delta[n])$,
so that it is enough to exhibit compatible maps from
\[  Y_n\times \varphi^\lad(\Delta[n])~=~Y_n\times\hocolimsub{\twosub{[m]\to [n]}{\textup{ in }\Delta_\inj}}\Delta[m] \]
to the hocolim in~\eqref{eqn-factosol}, for $n\ge 0$. (The compatibility checks will be left to the reader, however.)
In degree $k$, the right-hand side of this last equation is
\[   Y_n\times \hocolimsub{\twosub{[m]\to [n]}{\textup{ in }\Delta_\inj}}\mor_\Delta([k],[m]). \]
This is a quotient of a disjoint union of pieces $Y_n\times \Delta^p$ corresponding to pairs $(f,\sD)$ where
$f\in \mor_\Delta([k],[m])$ and
\begin{equation} \label{eqn-beforedevissage}
\sD: \qquad [n]=[m_0] \xleftarrow{g_1} [m_1] \xleftarrow{g_2} [m_2]  \xleftarrow{g_3} \cdots \xleftarrow{g_p} [m_p]=[m] \end{equation}
is a diagram much as in item (i) at the beginning of this proof; the $g_j$ are morphisms in $\Delta_\inj$,
and we can assume that none of them are identity maps. Therefore we have to exhibit compatible maps
\begin{equation} \label{eqn-devissage} Y_n\times\Delta^p \to \hocolimsub{S}  U(k,0)_S\,, \end{equation}
one for each of these pairs $(f,\sD)$. Now $\sigma(\sD)$ is a $p$-simplex in $\sP(Y_n,\sU(n,a))$ where $a=g_1g_2\cdots g_p(0)$.
It can also be viewed as a $p$-simplex in $\sP(Y_n,\sU(n,b))$ where $b=g_1g_2\cdots g_pf(0)$.
Pulling this back along $g_1g_2\cdots g_p f$ we obtain a $p$-simplex in $\sP(Y_k,\sU(k,0))$, which gives us the
second arrow in
\begin{equation} \label{eqn-enjoydevissage}
Y_n\times \Delta^p \xrightarrow{(g_1g_2\cdots g_pf)^*\times\id}
Y_k\times \Delta^p \lra \hocolimsub{S} U(k,0)_S
\end{equation}
in the manner of lemma~\ref{lem-SegalPUext}. The composition of the two arrows is the map that we have been looking for. \qed

\bigskip
We  return to configuration categories. For an open $\alpha$-cover $\sV=(V_j)_{j\in J}$ of a (smooth, complete) Riemannian
manifold $L$ we have $\cconfig_{\sV}(L;\alpha)$ as in definition~\ref{defn-models1}.
This can also be written as the co\-limit, over all finite nonempty subsets $S \subset J$, of the $\config(V_S;\alpha)$
where $V_S=\bigcap_{j\in S} V_j$.

\begin{prop} \label{prop-dersec} There exists a factorization
\[
\xymatrix@R=20pt@C=10pt@M=6pt{
 &   {\hocolimsub{S\subset J} \displaystyle{\cconfig(V_S;\alpha)}} \ar[d] \\
{\varphi^\lad(\cconfig_{\sV}(L;\alpha))} \ar@{..>}[ur] \ar[r] & \cconfig_\sV(L;\alpha)
}
\]
\end{prop}
(The horizontal arrow is adjoint to the inclusion of $\cconfig_\sV(L;\alpha)$ in $\varphi(\cconfig_\sV(L;\alpha)$.
The vertical arrow is the canonical map from the homotopy colimit to the colimit. It is a weak equivalence. The factorization is strict, not just up to homotopy.)

\proof[Proof of proposition~\ref{prop-dersec}] This is a straightforward application of proposition~\ref{prop-dersecabstract}.
We choose $Y:= \cconfig_{\sV}(L;\alpha)$. The open $\alpha$-cover $\sV$ of $L$ determines an open cover $\sU$ of $Y_0$ with the same
indexing set $J$. Namely, a multipatch $W\in Y_0$ is said to be an element of $U_i \subset Y_0$ if, as a subset of $L$, it is contained in $V_i$
(for the same $i\in J$). An element of $U(k,j)_i$ is then a nested system of multipatches
\[
W_0 \supset W_1 \supset \dots \supset W_k
\]
such that $W_j \subset U_i$. (Confusion alert: $j\in \{0,1,2,\dots,k\}$ but $i\in J$.)
Clearly, $U(k,j)_i \subset U(k,\ell)_i$ for all $\ell \in \{j,...,k\}$. So the assumptions of proposition \ref{prop-dersecabstract} are met.
\qed

\begin{rem} \label{rem-spill} Let $W\subset L$ be an open set and let $\sV|_W$ be the open cover
of $W$ consisting of the $V_j\cap W$, where $V_j\in \sV$. Then there is a
(strictly) commutative diagram
\[
	\begin{tikzpicture}[descr/.style={fill=white}, baseline=(current bounding box.base)] ]
		\matrix(m)[matrix of math nodes, row sep=1.5em, column sep=2.5em,
	text height=1.5ex, text depth=0.25ex]
	{
	\varphi^{\lad} \cconfig_{\sV|_W}(W;\alpha) & \hocolimsub{S\subset J} \cconfig(V_S\cap W;\alpha) \\
	\varphi^{\lad} \cconfig_{\sV}(L;\alpha) & \hocolimsub{S\subset J} \cconfig(V_S;\alpha) \\
	};
	\path[->,font=\scriptsize]
		(m-1-1) edge node [auto] {} (m-1-2)
		(m-1-1) edge node [auto] {} (m-2-1)
		(m-1-2) edge node [auto] {} (m-2-2)
		(m-2-1) edge node [auto] {} (m-2-2);
	\end{tikzpicture}
\]
\vskip2mm\noindent
where the vertical arrows are inclusion maps and the lower horizontal arrow is the broken arrow from proposition~\ref{prop-dersec}.
Later we may modify or simplify this. For example, in the upper right-hand term, those $S\subset J$
for which $V_S\cap W$ is empty make no contribution to the hocolim and can be left out.
\end{rem}

\section{Grounding and bounding} \label{sec-grobo}
We return to the setting and hypotheses of section~\ref{sec-pseudomet}. Specifically we have a complete Riemannian manifold $M$ and a compact
smooth submanifold $L \subset M$ (empty boundary) with the induced Riemannian metric. As in section~\ref{sec-pseudomet} we also assume
that $M$ is connected. Let $\Pi_M\co A_M\to M$ be a universal covering space and let $\Pi_L\co A_L\to L$ be the restricted covering space, so that there
is a strict pullback square
\[
\xymatrix@R=16pt{
A_L   \ar[r] \ar[d]_-{\Pi_L}  & A_M \ar[d]^-{\Pi_M} \\
L \ar[r]  & M
}
\]
Now $A_M$ inherits a Riemannian metric from $M$ such that $\Pi_M$ is a local isometry. Similarly $A_L$ inherits
a Riemannian metric from $L$ such that $\Pi_L$ is a local isometry. Important for us: the preimage under $\Pi_M$ of a
geodesic ball in $M$ (center $z$, radius less than the injectivity radius of $M$ at $z$) is a disjoint union of geodesic
balls in $A_M$, of the same radius. Something analogous can be said about $L$ and $\Pi_L$.

\begin{defn} \label{defn-ovmapnotn} For bookkeeping purposes in connection with configuration categories and
covering spaces, we introduced in \cite[Def. 2.1]{confcover} the category $\epifin$ with two forgetful functors to $\fin$.
The objects of $\epifin$ are surjective maps of finite sets, $\uli k\to \uli\ell$, subject to a little
condition (``selfic'') which ensures that each isomorphism class in $\epifin$ has only one object. A morphism is a commutative square, i.e. a map from $\uli k\to \uli\ell$ to $\uli k^\prime\to \uli\ell^\prime$, whose induced maps on fibers are injective.
The two forgetful functors are, of course, $(\uli k\to \ell)\mapsto \uli k$ and $(\uli k\to \uli\ell)\mapsto \uli\ell$.

Using notation as in \cite{confcover} and in appendix \ref{sec-models}, we have the ``covering configuration category'' $\config(\Pi_L)$
(a Segal space over $\epifin$) and forgetful maps of Segal spaces
\[ \config(\Pi_L)\to \config(A_L),\qquad \config(\Pi_L)\to \config(L) \]
which cover the respective forgetful functors from $\epifin$ to $\fin$.

In this section and the next one, we need (Rezk) complete variants of all these. Therefore we have
$\cepifin$ with two forgetful maps $\cepifin\to \cfin$, as well as $\cconfig(\Pi_L)$
(a complete Segal space over the complete Segal space $\cepifin$) and forgetful maps
\[ \cconfig(\Pi_L)\to \cconfig(A_L),\qquad \cconfig(\Pi_L)\to \cconfig(L) \]
which cover the respective forgetful maps from $\cepifin$ to $\cfin$.
\end{defn}

\begin{defn} \label{defn-grpatches} The second forgetful functor $\epifin\to \fin$ has a
section (right inverse) $s\co \fin\to \epifin$ which on objects $\uli\ell$ of $\fin$ is specified by
$\uli\ell \mapsto (\id\co \uli\ell\to \uli\ell)$. There is a completed variant of this, $\cfin\to \cepifin$,
which we also denote by $s$ (if at all). Let
\[  \cconfig^\sol(L)   \]
be the pullback (limit) of $\cconfig(\Pi_L)\to \cepifin\xleftarrow{\,\,s\,\,} \cfin$. Then we have an inclusion
$\cconfig^\sol(L)\to \cconfig(\Pi_L)$ and a forgetful projection $\cconfig^\sol(L)\to \cconfig(L)$,
which is the restriction to $\cconfig^\sol(L)$ of the forgetful projection
$\cconfig(\Pi_L)\to \cconfig(L)$.

An element in degree $0$ of $\cconfig^\sol(L)$, a.k.a.~object, can be thought of as a multipatch $U$ in $L$
together with a section of $\Pi_L\co A_L\to L$ defined on $U$. An element in degree 1
of $\cconfig^\sol(L)$ is an inclusion of multipatches in $L$ respecting the section data.
It follows that the projection $\cconfig^\sol(L)\to \cconfig(L)$ has a few noteworthy properties.
\begin{itemize}
\item[(a)] It is a right fibration with discrete fibers. (See definition~\ref{defn-leftfib}, which is
the definition of a left fibration with discrete fibers.)
\item[(b)] It is distributive. In other words, if a multipatch $U$ is the disjoint union of multipatches $U_0$ and $U_1$, then the fiber over $U$ is the product of the fibers of $U_0$ and $U_1$.
\item[(c)] It is invariant under the action of $\Gamma$ on $\cconfig^\sol(L)$ by left translation.
\end{itemize}

We can make similar definitions for $M$ and $\Pi_M$ in place of $L$ and $\Pi_L$. A new aspect here is that
we are also interested in fibrant replacements. Therefore it is worth noting that the commutative square
\[
\xymatrix@R=20pt{
\cconfig^\sol(M) \ar[d]  \ar[r] & \varphi\cconfig^\sol(M) \ar[d] \\
\cconfig(M) \ar[r] & \varphi\cconfig(M)
}
\]
is a strict pullback square, and the two vertical arrows are both right fibrations with discrete fibers.
\end{defn}

\begin{defn} \label{defn-ovmapconf}
A \emph{grounded map} over $\cfin$ from $\cconfig(L)$ to $\varphi\cconfig(M)$ is a pair $(f,f^\infty)$ as in the following commutative square
of maps over $N\cfin$,
\[
\xymatrix@R=20pt{
\cconfig^\sol(L) \ar[d]  \ar[r]^-{f^\infty} & \varphi\cconfig^\sol(M) \ar[d] \\
\cconfig(L) \ar[r]^-f & \varphi\cconfig(M)
}
\]
where the vertical arrows are forgetful, and the upper horizontal arrow \emph{respects the left actions} of
the covering translation group $\Gamma$ of $\Pi_M$. In other words, it is a map
from $\cconfig(L)$ to $\varphi\cconfig(M)$ which is ``covered'' by a map of right fibrations
with discrete fibers, the vertical arrows in the diagram (and there is an equivariance condition).
The equivariance condition is redundant in the cases that we are most interested in. It is automatically
satisfied if the inclusion $L\to M$ induces an isomorphism of fundamental groups.  --- Write
\[
\grmap_{\cfin}(\cconfig(L), \varphi\cconfig(M))
\]
for the ``space'' (simplicial set) of these grounded maps. Embellishments can be added, as in
$\grmap_{\cfin}(\cconfig_\delta(L;\alpha), \varphi\cconfig(M))$. We may also use notation like
\[  \grmap_\cfin(\cconfig(W_0),\cconfig(W_1)) \]
where $W_0$ is an open subset in $L$ and $W_1$ is open in $M$. This should be self-explanatory. ---
There is a forgetful map
\[
\grmap_{\cfin}(\cconfig(L), \varphi\cconfig(M)) \lra \map_{\cfin}(\cconfig(L), \varphi\cconfig(M))
\]
of simplicial sets, and it is a covering space. There is another forgetful map from
$\grmap_{\cfin}(\cconfig(L), \varphi\cconfig(M))$ to $\map_{\cfin}(\cconfig(L), \varphi\cconfig(A_M))$
which is also quite important to us. In most cases this is not a covering space.
\end{defn}

\begin{rem} \label{rem-othergr} The definition of \emph{grounded map} (from $\grcconfig(L)$ to $\varphi\grcconfig(M)$) over $N\cfin$ has other equivalent
formulations. In the definition, as given above, of a grounded map as a pair $(f,f^\infty)$, the $f^\infty$ is already
fully determined by $f$ and the restriction of $f^\infty$ to the space of ``single patches in $L$ with a lift to
$A_L$''. This follows from properties (a) and (c) in definition~\ref{defn-grpatches}.
More precisely, there is a strict pullback square
\[
\xymatrix@R=16pt@C=20pt{
\grmap_\cfin(\cconfig(L),\varphi\cconfig(M)) \ar[d] \ar[r] & \map^\Gamma(A_L,A_M) \ar[d] \\
\map_\cfin(\cconfig(L),\varphi\cconfig(M)) \ar[r] & \map(L,M)
}
\]
where $\map^\Gamma(A_L,A_M)$ is the
space of $\Gamma$-maps from $A_L$ to $A_M$. This amounts to another definition of
$\grmap_\cfin(\cconfig(L),\varphi\cconfig(M))$. It is very reminiscent of Prop.~3.6. in \cite{confcover}.
(For the lower horizontal arrow in that pullback square we need to make
a choice. The forgetful map $\cconfig(L;1)_0\to L$ taking single patches in $L$ to their
center points is a fibration with contractible fibers. Choose a section $t$ for it. A map
$\cconfig(L)\to\varphi\cconfig(M)$ can be restricted to give a map of spaces from $\cconfig(L;1)_0$
to $\cconfig(M;1)_0$. This can be precomposed with $t$ and postcomposed with the projection
$\cconfig(M;1)_0\to M$. This gives the lower horizontal arrow. The upper horizontal arrow is then determined.)
This reformulation of definition~\ref{defn-ovmapconf} is very reminiscent of Prop.~3.6. in \cite{confcover}.

At the opposite extreme, a grounded map from $\cconfig(L)$ to $\varphi\cconfig(M)$ could also be defined
as a pair of maps $(f,f^\infty)$ contributing to a commutative square
of maps over $\cfin$,
\[
\xymatrix@R=20pt{
\cconfig(\Pi_L) \ar[d]  \ar[r]^-{f^\infty} & \varphi\cconfig(\Pi_M) \ar[d] \\
\cconfig(L) \ar[r]^-f & \varphi\cconfig(M)
}
\]
where the vertical arrows are forgetful, and $f^\infty$ respects the left actions of $\Gamma$.

The point of view of definition~\ref{defn-ovmapconf}
has some advantages which can be appreciated in definition~\ref{defn-swepcontrol} below.
\end{rem}

\begin{defn} \label{defn-swepcontrol}
Let $\vep$ be a positive real number. An element alias $0$-simplex
$(f,f^\infty)$ of $\grmap_{\cfin}(\cconfig_\delta(L;\alpha), \varphi\cconfig(M))$ is
\emph{$\vep$-bounded} if $f$ is $\vep$-bounded. (See definitions~\ref{defn-epthick} and~\ref{defn-eprcontrol}.)
It is \emph{strongly $\vep$-bounded} if, informally speaking, $f^\infty$ is $\vep$-bounded.

More precisely, if $V$ is an object of $\cconfig^\sol_\delta(L)$, then we can view that as an
open subset of $A_L$, and we obtain a map
\[  \cconfig_\delta(V) \hookrightarrow \cconfig^\sol_\delta(L) \xrightarrow{f^\infty}
\varphi\cconfig^\sol(M) \xrightarrow{\textup{forget}} \varphi\cconfig(A_M). \]
The condition is that this lands in $\varphi\cconfig(\Theta_\vep V)\subset \varphi\cconfig(A_M)$.

For $k$-simplices of $\grmap_{\cfin}(\cconfig_\delta(L;\alpha), \varphi\cconfig(M))$ where $k>0$, there are
similar definitions. Such a $k$-simplex is a family of pairs $(f_z,f^\infty_z)$ depending continuously on $z\in \Delta^k$,
and we say that it is $\vep$-bounded, resp.~strongly $\vep$-bounded, if every $(f_z,f^\infty_z)$ in the family
is $\vep$-bounded, resp.~strongly $\vep$-bounded. Let
\[ \underline{Z}_\vep\subset \grmap_{\cfin}(\cconfig_\delta(L;\alpha), \varphi\cconfig(M)) \]
be the simplicial subset consisting of the strongly $\vep$-bounded simplices. Let $\underline{Z}_\infty$
be the union of the $\underline{Z}_\vep$ for all $\vep>0$. We may also
write $\underline{Z}_\vep(L,M)$ resp.~$\underline{Z}_\infty(L,M)$ for these simplicial subsets if it helps to avoid confusion.
\end{defn}
\emph{Remarks.} a) If a simplex of $\grmap_{\cfin}(\cconfig_\delta(L;\alpha), \varphi\cconfig(M))$ is
strongly $\vep$-bounded, then it is
$\vep$-bounded. Reason: $\Pi_M$ is distance-non-increasing.
b) The simplicial set $\grmap_{\cfin}(\cconfig_\delta(L;\alpha), \varphi\cconfig(M))$
is a fibrant simplicial set (a.k.a. Kan simplicial set) by construction, and $\underline{Z}_\vep$ is also fibrant
by construction.

\begin{thm} \label{thm-grbound} The inclusion $\underline{Z}_\infty
\hookrightarrow \grmap_{\cfin}(\cconfig_\delta(L;\alpha), \varphi\cconfig(M))$
admits a homotopy right inverse.
\end{thm}
The $\delta$ turns out to have little influence in the proof of this, and so it will be written out for
$\cconfig(L;\alpha)$. (That case could be simulated by taking for $\delta$ a number larger than the diameter of $L$.)
Apart from that, the proof requires a great deal of preparation.

\begin{lem} \label{lem-funnycov} The smooth compact manifold $L$ admits an open $\alpha$-cover $\sV=(V_j)_{j\in J}$ with finite indexing set $J$
such that, for each nonempty subset $S\subset J$, the set $V_S=\bigcap_{j\in S} V_j$ is homeomorphic to a disjoint
union of finitely many copies of $\RR^\ell$.
\end{lem}
\proof Begin with a smooth triangulation $\sT=\sT^0$ of $L$. Let $\sT^k$ be the $k$-fold barycentric subdivision of $\sT^0$.
Let $\sW^{(k)}$ be the open cover of $L$ defined as follows: $V\in \sW^{(k)}$ iff and only if $V$ is a disjoint union of
finitely many subsets each of which is an open star (of a vertex) in one of the triangulations $\sT^0,\sT^1,\dots,\sT^k$. Let us show
that for sufficiently large $k$, the open
cover $\sW^{(k)}$ is an open $\alpha$-cover of $L$. To this end let
$C_k\subset \mapt(\{1,2,\dots,\alpha\},L)$ consist of all $f\co\{1,2,\dots,\alpha\}\to L$ such that $\im(f)$ is not contained
in any of the open sets
which make up $\sW^{(k)}$. Then $C_k$ is compact for all $k$, and we have $C_0\supset C_1\supset C_2\supset ...$, and moreover
the intersection of all the $C_k$ is empty. It follows that there is some $k\ge 0$ for which $C_k$ is empty.
Set $\sV:=\sW^{(k)}$ for such a $k$. \qed

\medskip
Choose an open cover $\sV=(V_j)_{j\in J}$ as in lemma~\ref{lem-funnycov}. The open cover determines a poset $\sK$
as follows. The elements of $\sK$ are the pairs $(S,T)$ where $S$ is a nonempty subset of $J$ and $T$ is a subset of $\pi_0(V_S)$ of cardinality
$\le\alpha$. (We allow $T=\emptyset$.) The ordering has $(S_0,T_0)\le (S_1,T_1)$ if and only if $S_1\subset S_0$ and
the inclusion
\[   V_{S_0}\to V_{S_1} \]
takes every connected component selected by $T_0$ to a component selected by $T_1$. (The resulting map $T_0\to T_1$
does not have to be injective.)

For $(S,T)\in \sK$ let $V_{S,T}\subset V_S$ be the union of the connected components of $V_S$ selected by $T$.
Then $(S,T)\mapsto V_{S,T}$ is a covariant functor.
The inclusions $\cconfig(V_{S,T};\alpha)\to \cconfig(L;\alpha)$ determine a map
\[  \hocolimsub{(S,T)\textup{ in }\sK} \cconfig(V_{S,T};\alpha) \lra \cconfig_\sV(L;\alpha). \]

\begin{lem} \label{lem-pedantic} There exists a \emph{(strict)} factorization
\begin{equation} \label{eqn-dersecagain}
\begin{aligned}
\xymatrix@R=23pt@C=10pt@M=6pt{
 &   {\hocolimsub{(S,T)\textup{ in }\sK} \cconfig(V_{S,T};\alpha)} \ar[d]^-\simeq \\
{\varphi^\lad(\cconfig_\sV(L;\alpha))} \ar@{..>}[ur] \ar[r] & \cconfig_\sV(L;\alpha)
}
\end{aligned}
\end{equation}
\end{lem}
One important message of the lemma is that there is a homotopy colimit decomposition of the Segal space $\cconfig_\sV(L;\alpha)$,
indexed by a finite poset, in which every piece $\cconfig(V_{S,T};\alpha)$ is a Segal space
which admits a homotopy terminal object.

\proof The following commutative diagram is our guide.
\begin{equation} \label{eqn-stairs}
\begin{aligned}
\xymatrix@R=20pt@C=8pt@M=5pt{
&  & {\varphi\!\hocolimsub{(S,T)\textup{ in }\sK} \displaystyle{\cconfig(V_{S,T};\alpha)}} \ar[d] \\
  &   {\hocolimsub{S\subset J} \displaystyle{\cconfig(V_S;\alpha)}} \ar[d] \ar@{..>}[ur]^-H \ar[r] &
  {\varphi\,\hocolimsub{S\subset J} \displaystyle{\cconfig(V_S;\alpha)}} \ar[d] \\
{\varphi^\lad(\cconfig_\sV(L;\alpha))} \ar[r] \ar[ur]^-{u}  & \cconfig_\sV(L;\alpha) \ar[r]  & \varphi\cconfig_\sV(L;\alpha)
}
\end{aligned}
\end{equation}
The three vertical arrows are induced by inclusions. Two of the horizontal arrows are preferred
inclusions and the other one is a preferred projection. The map $u$ comes from proposition~\ref{prop-dersec}. The map $H$
is something that we have to design now,
and we use remark~\ref{rem-imagine} for that.
Fix some $k\ge 0$. For every nonempty $S\subset J$ there is a preferred embedding
\[  \cconfig(V_S;\alpha)_k \lra \coprod_T \cconfig(V_{S,T};\alpha)_k   \]
(whose image is a summand, a.k.a~union of connected components) because for every element $z$ of $\cconfig(V_S;\alpha)_k$
there is a minimal $T$ such that $z\in \cconfig(V_{S,T};\alpha)_k$. These maps are natural in $S$ (for fixed $k$), and so they determine a map
\[  \hocolimsub{S\subset J} \cconfig(V_S;\alpha)_k \lra \hocolimsub{(S,T)\textup{ in }\sK} \cconfig(V_{S,T};\alpha)_k\,. \]
Let this be $H_\sD$ in the case where $\sD$ is the diagram in $\Delta_\inj$ consisting of $[k]$ only. The maps
$H_\sD=H_{[k]}$ are not claimed to be compatible with the simplicial operators in
$X:=\hocolim_{S\subset J} \cconfig(V_S;\alpha)$, respectively
$Y:=\hocolim_{(S,T)} \cconfig(V_{S,T};\alpha)$.
But for $g\co [j]\to [k]$ in $\Delta_\inj$ there is a preferred homotopy
relating $g^*H_{[k]}$ to $H_{[j]}\,g^*$ because for every $z\in \cconfig(V_S;\alpha)_k$ the minimal $(S,T_1)$ in
$\sK$ for which $z\in \cconfig(V_{S,T_1};\alpha)$ is $\ge$ the minimal $(S,T_0)$ for which
$g^*z\in \cconfig(V_{S,T_0};\alpha)$.
Let this homotopy be $H_\sD$ in the case where $\sD$ is
$g\co [j]\to [k]$. And so on; a diagram
\[ \sD = (\,[k_0] \xleftarrow{g_1} [k_1]\xleftarrow{g_2} \cdots \xleftarrow{g_{r-1}} [k_{r-1}]\xleftarrow{g_r} [k_r]\,) \]
in $\Delta_\inj$ determines a map or higher homotopy
\[  H_\sD \co X_{k_0}\times \Delta^r \lra Y_{k_r}\,. \]
These maps $H_\sD$ satisfy the naturality properties listed in remark~\ref{rem-imagine}. Moreover they are maps over $X$.
Therefore broken arrow $H$ has come to life and it does make diagram~\eqref{eqn-stairs} commutative.
--- Now $Hu$ in diagram~\eqref{eqn-stairs} is a map which we can write in the form
\[  \varphi^\lad\varphi^\lad(\cconfig_\sV(L;\alpha)) \lra \hocolimsub{(S,T)\textup{ in }\sK} \cconfig(V_{S,T};\alpha).   \]
Then we may pre-compose with $\varphi^\lad\to \varphi^\lad\varphi^\lad$ of corollary~\ref{cor-phiadcomonad}
to obtain the broken arrow in~\eqref{eqn-dersecagain}. Commutativity in~\eqref{eqn-dersecagain} follows easily from
commutativity in~\eqref{eqn-stairs}.

It remains to be shown that the vertical arrow in diagram~\eqref{eqn-dersecagain} is a weak equiva\-lence.
Observe that the rule taking a finite
nonempty $S\subset J$ to the poset of all $(R,T)\in\sK$ with $R=S$ can be made
into a contravariant functor (from the poset of nonempty finite subsets of $J$ to the category of finite
posets). On the basis of that observation we can set up a diagram
\[
\xymatrix@R=16pt@C=20pt@M=6pt{
{\hocolimsub{S\subset J} \hocolimsub{T:~ (S,T)\in \sK} \cconfig(V_{S,T};\alpha)} \ar[d] \\
{\hocolimsub{(S,T)\in \sK} \cconfig(V_{S,T};\alpha)} \ar[r] & {\hocolimsub{S\subset \sK} \cconfig(V_S;\alpha)}
}
\]
in which the vertical arrow is a weak equivalence by the general theory
of homotopy colimits \cite[\S9]{DwyerKan}. Then it suffices to show that
the composition of these two arrows is a weak equivalence. The composition is the map of
homotopy colimits induced by a natural transformation
\[
\hocolimsub{T:~ (S,T)\in \sK} \cconfig(V_{S,T};\alpha)  \lra \cconfig(V_S;\alpha) \]
of contravariant functors in the variable $S$. It suffices to show that this gives a weak equivalence
for each $S$. But this is clear since for a fixed $S$ the open sets $V_{S,T}$ form an $\alpha$-cover of $V_S$. \qed

\begin{lem} \label{lem-surprise} For fixed $(S,T)\in \sK$, the inclusions $W\hookrightarrow M$ for
$W\in \cconfig(M)_0$
determine a weak equivalence
\[
\xymatrix@R=20pt{
{\hocolimsub{W} \grmap_\cfin(\cconfig(V_{S,T};\alpha),\varphi\cconfig(W))} \ar[d] \\
{\grmap_\cfin(\cconfig(V_{S,T};\alpha),\varphi\cconfig(M))}
}
\]
\emph{(The homotopy colimit is taken over the \emph{discrete} poset of multipatches in $M$.)}
\end{lem}
\proof Let $U$ be a weakly terminal object for $\cconfig(V_{S,T};\alpha)$. That is to say, $U$ is a multipatch in $V_{S,T}$
such that the inclusion $U\to V_{S,T}$ induces a bijection in $\pi_0$. Let $\star_U\subset \cconfig(V_{S,T};\alpha)$
be the simplicial subspace generated by $U$ as an element of $\cconfig(V_{S,T};\alpha)$ in degree 0, so that $\star_U$ has
exactly one element in each degree. Then for every multipatch $W$ in $M$
there is a homotopy cartesian square
\[
\xymatrix{
{\grmap_\cfin(\cconfig(V_{S,T};\alpha),\varphi\cconfig(W))} \ar[d]
\ar[r] & \ar[d] {\grmap_\cfin(\cconfig(V_{S,T};\alpha),\varphi\cconfig(M))} \\
{\map_\cfin(\star_U,\varphi\cconfig(W))}
\ar[r] & {\map_\cfin(\star_U,\varphi\cconfig(M))}
}
\]
where the vertical arrows are restriction maps.
This follows easily from the (derived) universal property of $U$ and the locality statement \cite[Cor.3.6]{BoavidaWeiss2018}.
Therefore it suffices to show that the map
\[
{\hocolimsub{W} \map_\cfin(\star_U,\varphi\cconfig(W))} \lra
{\map_\cfin(\star_U,\varphi\cconfig(M))}
\]
determined by the inclusions $W\hookrightarrow M$ is a weak equivalence. Evidently this can be simplified, using the
weak homotopy invariance of $\hocolim$. We can choose a bijection $\uli k\to \pi_0 U$. Then the
simplified map is
\[
\hocolimsub{W} \emb(\uli k,W) \lra \emb(\uli k,M)\,.
\]
It is a weak homotopy equivalence because it is a Serre microfibration with contractible fibers \cite[Lem. 2.2]{WeissWhatclass}. (The fiber
over some $f\in \emb(\uli k,M)$ is the classifying space of the poset of all multipatches $W$ in $M$ which contain
the image of $f$.) \qed

\medskip
\proof[Proof of theorem~\ref{thm-grbound}] The most important observation here is that for
an element $(S,T)\in \sK$ and a multipatch $W$ in $M$,
any grounded map from $\cconfig(V_{S,T};\alpha)$ to $\varphi\cconfig(W)$ is automatically strongly $\vep$-bounded
for some $\vep$, if we view it as a map to $\varphi\cconfig(M)$. Another useful observation, or a pair of useful
observations, is that
lemma~\ref{lem-pedantic} gives us a weak equivalence of simplicial sets
\[
\xymatrix@M=6pt@R=15pt{
{\begin{array}{cl}
& \holimsub{(S,T)} \map_\cfin(\cconfig(V_{S,T};\alpha),\varphi\cconfig(M)) \\
\cong & \map_\cfin(\hocolimsub{(S,T)}\cconfig(V_{S,T};\alpha),\varphi\cconfig(M))
\end{array}} \ar[d] \\
{\map_\cfin(\cconfig(L;\alpha),\varphi \cconfig(M)).}
}
\]
By a mechanical refinement of that same lemma, this lifts to a weak equivalence
\[
\xymatrix@M=6pt@R=15pt{
{\holimsub{(S,T)} \grmap_\cfin(\cconfig(V_{S,T};\alpha),\varphi\cconfig(M))}
\ar[d] \\
{\grmap_\cfin(\cconfig(L;\alpha),\varphi \cconfig(M)).}
}
\]
which, by inspection, takes $\holim_{(S,T)}\,\underline{Z}_\infty(V_{S,T},M)$ to $\underline{Z}_\infty(L,M)$.
(The inspection should use the full strength of lemma~\ref{lem-pedantic},
including the commutativity of that triangle. Let $x$ be a simplex in $\holim_{(S,T)}\,\underline{Z}_\infty(V_{S,T},M)$. Then $x$ has finitely many ``coordinates''
$x_\sD$ corresponding to nondegenerate simplices $\sD$ in the nerve of the finite poset $\sK$.
For each $\sD$ we can choose a strong bound $\vep_\sD$
satisfied by the coordinate $x_\sD$. The maximum of the $\vep_\sD$ is a strong bound for the image of $x$
in $\underline{Z}_\infty(L,M)$.) Therefore the composition
\[
\xymatrix@R=20pt@M=5pt{
{\holimsub{(S,T)}\hocolimsub{W} \grmap_\cfin(\cconfig(V_{S,T};\alpha),\varphi\cconfig(W))} \ar[d] \\
{\holimsub{(S,T)} \grmap_\cfin(\cconfig(V_{S,T};\alpha),\varphi\cconfig(M))} \ar[d]^-\simeq \\
{\grmap_\cfin(\cconfig(L;\alpha),\varphi \cconfig(M))}
}
\]
lands in $\underline{Z}_\infty(L,M)$, and it is a weak equivalence by lemma~\ref{lem-surprise}.
\qed

\section{Bounded lifting}  \label{sec-bolift}
Now we need to make matters more complicated by introducing a finite index subgroup $\Gamma_1$
of the group $\Gamma$ of covering translations
of $\Pi_M\co A_M\to M$. The subgroup $\Gamma_1$ acts freely, tautologically and isometrically on $A_M$ and $A_L$. We write
$E_{M}$ and $E_{L}$ for the orbit spaces, respectively, and
\[  \pi_{M}\co E_{M} \lra M\,, \qquad \Pi_{M,\Gamma_1}\co A_M\lra E_M\,, \]
\[  \pi_{L}\co E_{L} \lra L\,, \qquad \Pi_{L,\Gamma_1}\co A_L\to E_L  \]
for the resulting covering spaces, so that $\pi_{M}\Pi_{M,\Gamma_1}=\Pi_M$ and $\pi_{L}\Pi_{L,\Gamma_1}=\Pi_L$.
All these maps are local isometries, by definition. The finite index assumption allows us to say
that $E_L$ is still compact, because $L$ is compact.

\smallskip
Let $\underline{\emb}^s(L,M)$ be the pullback (limit) of
\[  \emb^s(L,M) \to \map(L,M) \leftarrow \map^\Gamma(A_L,A_M) \]
where the first arrow is forgetful and the second one is obtained by passage to $\Gamma$-orbits.
In \cite{confcover} we constructed the broken arrow in a homotopy commutative diagram
\[
\xymatrix@M=5pt@R=20pt{
\underline{\emb}^s(E_{L},E_{M}) \ar[r] & \rgrmap_\fin(\config(E_{L}),\config(E_{M})) \\
\underline{\emb}^s(L,M) \ar[u] \ar[r] & \rgrmap_\fin(\config(L),\config(M)) \ar@{..>}[u]
}
\]
where the $\RR$ in $\rgrmap$ indicates passage to an unspecified derived setting. (Section 7 of
\cite{confcover} is closest to what we need here.)
This used the ordered configuration categories, and for these the particle models, which makes the horizontal arrows obvious.
Switching to the Riemannian patch models is inconvenient from this point of view, but we have to do it.

The ``abstract'' description of the broken arrow does not suffer much in the translation. In the language
used here it is given by a soon-to-be-described map
\begin{equation} \label{eqn-soontobe} 
\grmap_\cfin(\cconfig(L),\varphi\cconfig(M)) \lra \map_\cfin(\cconfig(E_{L}),\varphi\cconfig(E_{M}))
\end{equation} 
together with a preferred lift across
\[
\grmap_\cfin(\cconfig(E_{L}),\varphi\cconfig(E_{M})) \to \map_\cfin(\cconfig(E_{L}),\varphi\cconfig(E_{M})) \;.
\]
(To specify that lift we use the pullback square in remark~\ref{rem-othergr}.)

As for the map~\eqref{eqn-soontobe}, it is a composition 
\begin{equation} \label{eqn-nosuffer}
\begin{aligned}
\xymatrix@R=15pt{
\grmap_\cfin(\cconfig(L),\varphi\cconfig(M)) \ar[d] \\
\map_\cepifin(\cconfig(\pi_{L}),\varphi\cconfig(\pi_{M})) \ar[d] \\
\map_\cfin(\cconfig(\pi_{L}),\varphi\cconfig(E_{M})) \ar[d] \\
\map_\cfin(\cconfig(E_{L}),\varphi\cconfig(E_{M}))
}
\end{aligned}
\end{equation}
in which the first arrow is an instance of naturality, the second is induced by the
forgetful map $\cconfig(\pi_{M})\to \cconfig(E_{M})$, and the third one is a homotopy inverse for the
map
\begin{equation} \label{eqn-wonderwhy}
\begin{aligned}
\xymatrix@R=15pt{
\map_\cfin(\cconfig(E_{L}),\varphi\cconfig(E_{M})) \ar[d] \\
\map_\cfin(\cconfig(\pi_{L}),\varphi\cconfig(E_{M}))
}
\end{aligned}
\end{equation}
induced by the forgetful map $\cconfig(\pi_{L})\to \cconfig(E_{L})$. Of course,
\eqref{eqn-wonderwhy} is a homotopy equivalence because $\cconfig(\pi_{L})\to \cconfig(E_{L})$ is a
conservatization map over $\cfin$.

\smallskip
For us it is important
to have some form of metric understanding of ``the'' homotopy inverse for the map~\eqref{eqn-wonderwhy}. This is not provided by \cite{confcover},
but it is provided by the next theorem.

\begin{thm} \label{thm-step2} For positive $\vep$ and $\delta$, let $\rho:=(2\alpha+1)\vep+(\alpha+1)\delta$.
There exist posi\-tive $\delta_1<\delta$ and a filler making the following diagram homotopy commutative:
\begin{equation} \label{eqn-struggle}
\begin{aligned}
\xymatrix@R=14pt@M=4pt{
\underline{Z}_\vep(L,M) \ar@{..>}[dd] \ar[r]^-{\textup{incl.}} & \grmap_{\cfin}(\cconfig_\delta(L;\alpha), \varphi\cconfig(M)) \ar[d]  \\
 &  \map_{\cfin}(\cconfig_{\delta_1}(\pi_L;\alpha),\varphi \cconfig(E_M)) \\
{Z}_\rho (E_L,E_M) \ar[r]^-{\textup{incl.}}  &  \map_{\cfin}(\cconfig_{\delta_1}(E_L;\alpha), \varphi\cconfig(E_M)) \ar[u]^\simeq
}
\end{aligned}
\end{equation}
\end{thm}
\emph{Alert.} There is no underlining in ${Z}_\rho (E_L,E_M)$.
\proof
Choose a finite open $\alpha$-cover $\sW=(W_j)_{j\in J}$ of $L$ as in lemma~\ref{lem-funnycov},
and such that every connected component of every $W\in \sW$
has diameter $<\delta$. Then use lemma~\ref{lem-alphacover} to find $\delta_1>0$ such that
\[  \cconfig_{\delta_1}(L;\alpha) \subset \cconfig_\sW(L;\alpha). \]
Let $\sV$ be the open $\alpha$-cover of $E_L$ consisting of the $V_j:=\pi_L^{-1}(W_j)$ for $j\in J$. For a nonempty
$S\subset J$ let
\[ V_S=\bigcap_{j\in S}V_j=\pi_L^{-1}(W_S). \]
Let $\sQ$ be the set of all pairs
$(S,T)$ where $S$ is a nonempty subset of $J$ and $T\subset\pi_0V_S$ is a subset of cardinality $\le\alpha$.
For $(S,T)$ in $\sQ$ let $V_{S,T}\subset V_S$ be the union of the components selected by $T$. Make $\sQ$ into a poset
in such a way that $(S,T)\mapsto V_{S,T}$ is a covariant functor. As in lemma~\ref{lem-pedantic} we can make a map
\begin{equation} \label{eqn-brisklike2} \varphi^\lad\cconfig_\sV(E_L;\alpha) \lra \hocolimsub{(S,T)\in\sQ} \cconfig(V_{S,T};\alpha)
\end{equation}
over $\cconfig(E_L;\alpha)$.
On the other hand, by \cite{confcover} there is a preferred map
\[  |\underline{Z}_\vep(L,M)|\times \cconfig_\delta(\pi_L;\alpha) \lra \varphi\cconfig(E_M;\alpha). \]
By the definition of $\underline{Z}_\vep(L,M)$, this map can be refined mechanically to a
natural transformation of functors in the variable $(S,T)\in\sQ$,
\begin{equation} \label{eqn-underived}
|\underline{Z}_\vep(L,M)|\times \cconfig(V_{S,T}\to W_S;\alpha) \lra \varphi\cconfig(\Theta_\vep V_{S,T};\alpha).
 \end{equation}
Since the forgetful maps $\cconfig(V_{S,T}\to W_S;\alpha) \to \cconfig(V_{S,T};\alpha)$
are conservatization maps over $\cfin$, this factors up to derived homotopy
through a \emph{derived} natural transformation
\begin{equation} \label{eqn-bewarederived}
\xymatrix{|\underline{Z}_\vep(L,M)|\times \cconfig(V_{S,T};\alpha) \ar@{..>}[r] & \varphi\cconfig(\Theta_\vep V_{S,T};\alpha) }
\end{equation}
over $\cfin$. Here we have to be precise. Let $E_0$, $E_1$ and $F$ be the functors on $\sQ$ sending $(S,T)$ to
\[
\cconfig(V_{S,T} \to W_{S}; \alpha) \;\; , \; \; \cconfig(V_{S,T}; \alpha) \;\; \textup{and} \; \; \varphi \cconfig(\Theta_\varepsilon V_{S,T}; \alpha)
\]
respectively. Then \eqref{eqn-underived} amounts to a map
\begin{equation}\label{eqn-bddtoE0}
\underline{Z}_\vep(L,M) \to \map(E_0, F) \hookrightarrow \rmap(E_0, F)
\end{equation}
where $\map$ (and $\rmap$) refers to the space of (derived) natural transformations from $E_0$ to $F$. We have a forgetful map $E_0 \to E_1$ which is a conservatization, so that the induced map
\begin{equation}\label{eqn-E1toE0}
\rmap(E_1, F) \to \rmap(E_0, F)
\end{equation}
is a weak equivalence. Therefore, \eqref{eqn-bddtoE0} has a lift across \eqref{eqn-E1toE0} up to homotopy. This gives us \eqref{eqn-bewarederived}. Next,~\eqref{eqn-bewarederived} determines an honest map
\begin{equation} \label{eqn-bewarestill}
\xymatrix@R=6pt@C=20pt@M=6pt{  {|\underline{Z}_\vep(L,M)|\times \hocolimsub{(S,T)\in \sQ} \cconfig(V_{S,T};\alpha)} \ar[r] &
{\colimsub{(S,T)\in \sQ} \varphi\cconfig(\Theta_\vep V_{S,T};\alpha)} \\
& \subset \quad\varphi\cconfig(E_M;\alpha).
}
\end{equation}
Pre-composing in the second input variable with~\eqref{eqn-brisklike2} gives
\[  |\underline{Z}_\vep(L,M)|\times \varphi^\lad\cconfig_\sV(E_L;\alpha) \lra \varphi\cconfig(E_M;\alpha) \]
with an adjoint which we can write in the form
\[  \underline{Z}_\vep(L,M) \lra \map_{\cfin}(\varphi^\lad\cconfig_\sV(E_L;\alpha), \varphi\cconfig(E_M)). \]
In order to simplify this we trade the $\varphi^\lad$ prefix for a $\varphi$ prefix attached to $\varphi\cconfig(E_M)$,
and apply the monadic transformation $\varphi\varphi\to \varphi$ of lemma~\ref{lem-phimonad}. We may also replace the
$\sV$-condition (subscript to $\cconfig$) the stronger $\delta_1$-condition, and so we obtain at last
\begin{equation} \label{eqn-endofstruggleless}
\underline{Z}_\vep(L,M) \lra \map_{\cfin}(\cconfig_{\delta_1}(E_L;\alpha), \varphi\cconfig(E_M)).
\end{equation}
It remains to show two things:
\begin{itemize}
\item[(i)] the map~\eqref{eqn-endofstruggleless} makes diagram~\eqref{eqn-struggle} homotopy commutative, and
\item[(ii)] its image is contained in $Z_\rho(E_L,E_M)$.
\end{itemize}
We begin with (i). We are comparing two maps out of $\underline{Z}_\vep(L,M)$ with the same target
$\map_{\cfin}(\cconfig_{\delta_1}(\pi_L;\alpha),\varphi \cconfig(E_M))$. It is allowed
to post-compose with the map
\[
\xymatrix@R=10pt{
\map_{\cfin}(\cconfig_{\delta_1}(\pi_L;\alpha),\varphi \cconfig(E_M)) \ar[d] \\
\map_{\cfin}(\hocolimsub{(S,T)\in\sQ}\cconfig_{\delta_1}(V_{S,T}\to W_{S,T});\alpha) ,\varphi \cconfig(E_M))
}
\]
induced by the projection $\hocolim\to \colim$, because this map is a weak equivalence. Now we have two maps
out of $\underline{Z}_\vep(L,M)$ with the same target
\begin{equation} \label{eqn-fattarget}
\map_{\cfin}(\hocolimsub{(S,T)\in\sQ}\cconfig_{\delta_1}(V_{S,T}\to W_{S,T});\alpha) ,\varphi \cconfig(E_M)) 
\end{equation}
Both of these maps factor through $\RR \map(E_0, F)$ and the two maps from $\underline{Z}_\vep(L,M)$ to $\RR \map(E_0, F)$ in these factorizations are homotopic by construction.

The proof of (ii) is much harder. We fix an element of $|\underline{Z}_\vep(L,M)|$
and think of that as a grounded map $(f,f^\infty)$ from $\cconfig_\delta(L;\alpha)$
to $\varphi\cconfig(M;\alpha)$, over $\cfin$. Constructions~\eqref{eqn-underived} and~\eqref{eqn-bewarederived}
should be specialized accordingly. Let $g$ be the map from $\cconfig_{\delta_1}(E_L;\alpha)$
to $\varphi\cconfig(E_M;\alpha)$, over $\cfin$, which we obtain by evaluating~\eqref{eqn-endofstruggleless} on $(f,f^\infty)$.
Here are some guiding principles.
\begin{itemize}
\item[(a)] The bound
$\vep$ is still intact in~\eqref{eqn-underived}.
\item[(b)] In making the step from~\eqref{eqn-underived} to~\eqref{eqn-bewarederived}
and~\eqref{eqn-bewarestill}, we may lose information about metric bounds, due to ``scattering''. The reason is that \eqref{eqn-bewarederived} is a \emph{derived} natural transformation.
\item[(c)] No further damage is done by the passage from~\eqref{eqn-bewarestill} to~\eqref{eqn-endofstruggleless}.
This step relies mainly on pre-composition with the map~\eqref{eqn-brisklike2}.
This is a map over $\cconfig(E_L;\alpha)$,
so that it is metrically beyond reproach.
\end{itemize}
Therefore our main business is to understand the metric aspects of the step from~\eqref{eqn-underived} to~\eqref{eqn-bewarederived}
and~\eqref{eqn-bewarestill}, as in (b). Here is another guiding principle.
\begin{itemize}
\item[(d)] Suppose that $U\subset E_L$ and $V\subset E_M$ are open sets, and we have a map $h$ from $\cconfig_{\delta_1}(U;\alpha)$ to $\varphi\cconfig(V;\alpha)$ over $\cfin$. Then $h$ induces a map
from $\pi_0(U)$ to $\pi_0(V)$. Let $V_1\subset V$ be the union of the connected components which are in the
image of that map. Then $h$ takes $\cconfig_{\delta_1}(U;\alpha)$ to
$\varphi\cconfig(V_1;\alpha)\subset \varphi\cconfig(V;\alpha)$.
\end{itemize}
Fix some $(S,T)\in \sQ$. It may look as if
$\cconfig(V_{S,T};\alpha)$ is taken to $\cconfig(\Theta_\vep V_{S,T})$ under the map~\eqref{eqn-bewarederived}. This would suggest that $g$ takes $\cconfig_{\delta_1}(V_{S,T};\alpha)$ to $\cconfig(\Theta_\vep V_{S,T})$.
But we must beware of scattering as in guiding principle (b).
The stark truth is that~\eqref{eqn-bewarederived} is a collection of (compatible) maps
\begin{equation} \label{eqn-laststand}
|\underline{Z}_\vep(L,M)|\times \Delta^r \times \cconfig(V_{S_0,T_0};\alpha)
\lra \varphi\cconfig(\Theta_\vep V_{S_r,T_r};\alpha).
\end{equation}
There is one such map for every string $(S_0,T_0)\le (S_1,T_1)\le \cdots \le (S_r,T_r)$
in $\sQ$. Therefore we must be more cautious. For fixed $(S,T)$ in $\sQ$, let $\sN(V_{S,T})\subset E_M$ be the union of all connected
components $C$ of all $\Theta_\vep V_{S',T'}$ such that $(S,T)\le (S',T')$ in $\sQ$
and $C\cap \Theta_\vep V_{S,T}$ is nonempty.
By guiding principle (d), it is safe to say that $g$ takes $\cconfig_{\delta_1}(V_{S,T};\alpha)$ to $\varphi \cconfig(\sN(V_{S,T}))$.
We can also estimate the size of $\sN(V_{S,T})$. Any connected component $C$ of $\Theta_\vep V_{S',T'}$
must have diameter $<\alpha(2\vep+\delta)$. (The $\varepsilon$-neighborhood of each connected component of $V_{S^\prime, T^\prime}$ has diameter $<2\varepsilon + \delta$, but there are $\leq \alpha$ of these and they may intersect.) If $C$ has nonempty intersection with $\Theta_\vep V_{S,T}$,
then it is contained in $\Theta_{\alpha(2\vep+\delta)+\vep}V_{S,T}$. Therefore, $\sN(V_{S,T})$ is contained in $\Theta_{\alpha(2\vep+\delta)+\vep}V_{S,T}$.

Now, suppose that $U$ is a multipatch in $E_L$ with not more than $\alpha$ connected components, all
of which have diameter less than $\delta_1$. Then we can find some $V_{S,T}$ containing $U$,
and we can take $T$ minimal, so that every component of $V_{S,T}$ contains some component of $U$.
So $V_{S,T}$ is contained in a $\delta$-neighborhood of $U$. We conclude that
$g$ takes $\cconfig(U;\alpha)$ to $\cconfig(\Theta_{\rho} U)$ where
\[ \rho= \alpha(2\vep+\delta)+\vep+\delta=(2\alpha+1)\vep+(\alpha+1)\delta.   \]
Therefore $g$ is $\rho$-bounded. \qed

\begin{rem}
There is a preferred lift
\[
\xymatrix@R=14pt@M=4pt{
\underline{Z}_\vep(L,M) \ar[d]^-{} \ar@{..>}[r] & \grmap_{\cfin}(\cconfig_{\delta_1}(E_L;\alpha), \varphi\cconfig(E_M))  \ar[d]  \\
{Z}_\rho (E_L,E_M) \ar[r]^-{\textup{incl.}}  &  \map_{\cfin}(\cconfig_{\delta_1}(E_L;\alpha), \varphi\cconfig(E_M))
}
\]
where the left-hand map is the dotted arrow in theorem \ref{thm-step2}.
Namely, we may compose on the right with the pullback square
\[
\xymatrix@R=14pt@M=4pt{
\grmap_{\cfin}(\cconfig_{\delta_1}(E_L;\alpha), \varphi\cconfig(E_M)) \ar[d] \ar[r] & \map^{\Gamma_1}(A_L, A_M)  \ar[d]  \\
 \map_{\cfin}(\cconfig_{\delta_1}(E_L;\alpha), \varphi\cconfig(E_M)) \ar[r]  &  \map(E_L,E_M)
}
\]
and then it suffices to produce a lift as in 
\[
\xymatrix@R=14pt@M=4pt{
\underline{Z}_\vep(L,M) \ar[d]^-{} \ar@{..>}[r] &  \ar[d] \map^{\Gamma_1}(A_L, A_M) \\
{Z}_\rho (E_L,E_M) \ar[r]^-{\textup{incl.}}  &  \map(E_L,E_M).
}
\]
Such a lift exists by the homotopy commutativity of \eqref{eqn-struggle}.
 
Note that this does not automatically give us a map from $\underline{Z}_\vep(L,M)$ to $\underline{Z}_\rho (E_L,E_M)$.
We did not establish a strong bound for the lifted map.
\end{rem}

\begin{appendices}
\section{Choice of geometric models} \label{sec-models}
Which models are we going to choose for the configuration categories~?
The question is important because we use (geo)metric arguments
in order to achieve the goals set out in section~\ref{sec-goal}. Here we make some
decisions.

\begin{defn} \label{defn-models1} For several reasons we will need unordered configuration categories more than
ordered configuration categories. (We authors have made some attempts to avoid them, but these attempts have come to nothing and we have
learned that such attempts must be viewed with unrelenting suspicion.)

Let $L$ be a smooth manifold with a complete Riemannian metric.
The preferred model for $\cconfig(L)$ is the Riemannian multipatch model. It is the nerve of a \emph{topological} poset
whose elements/objects are the multipatches. Since we will use
this consistently, there is no special notation for this, no embellishments. The multipatches are
finite disjoint unions of \emph{open} balls $B(x,s)$ (for $x\in L$ and $s>0$) defined using the geodesic metric. We always assume
that $s$ is less than the injectivity radius of the Riemannian manifold $L$ at the point $x$. In such a case
the center point of $B(x,s)$ is uniquely determined by $B(x,s)$; it is the element of $B(x,s)$
which has maximal distance from the complement of $B(x,s)$ in $L$.

One more condition:
we always assume that the patches do not touch, i.e., their closures are still disjoint.
This has the consequence that the object space of $L$ is a disjoint union of manifolds (without boundary).
The cardinality $k$ contribution is a manifold of dimension $k(\dim(L)+1)$.

Often we will impose stronger conditions on the multipatches. Typically these are meaningful when we
have an upper bound $\alpha$ on the \emph{number} of patches in a multipatch. We write $\cconfig(L;\alpha)$ as usual
for the corresponding full sub-poset. Suppose that $\sW$ is an open $\alpha$-cover of $L$. (This means that every
finite subset of $L$ having cardinality $\le\alpha$ is contained in some open set of $\sW$.)
We write $\cconfig_\sW(L;\alpha)$ for the full topological sub-poset of $\cconfig(L;\alpha)$ consisting of the
multipatches $V$ whose closure in $L$ is contained in some open set of $\sW$.

Let $\delta$ be a positive real number. We may write $\cconfig_\delta(L;\alpha)$
for the full topological sub-poset of $\cconfig(L;\alpha)$ consisting of the
multipatches $V$ in which every patch has diameter $<\delta$. This could also be written
in the form $\cconfig_\sW(L;\alpha)$, for example if we let $\sW$ consist of all the open subsets
$W$ of $L$ such that every connected component of $W$ has diameter $<\delta$.

Let $V\subset L$ be an open set. With the Riemannian metric induced from $L$, this may not be complete. In this context, we may nevertheless
write $\cconfig(V)$ for the nerve of the topological sub-poset of $\cconfig(L)$ consisting of all multipatches whose closure
is contained in $V$. Similarly, $\cconfig(V;\alpha)$ means $\cconfig(V)\cap \cconfig(L;\alpha)$.

We use Riemannian patch models for the ordered configuration categories, too. An object of $\config(L)$ is an object
$V$ of $\cconfig(L)$ with the extra datum of a bijection $\uli k\to \pi_0V$ for the appropriate $k\ge 0$. The morphisms
are inclusions as in $\cconfig(L)$. The extra data allow us to make a forgetful functor $\config(L)\to \fin$. The price for that is that
$\config(L)$ is not the nerve of a topological poset.
\end{defn}

\begin{expl} \label{expl-alphacover}
Let $L_1,L_2$ be closed topological manifolds and let $\pi\co L_1\to L_2$ be a
covering space, aka fiber bundle with discrete fibers. Let
$\sW$ be an open $\alpha$-cover of $L_2$. (Meaning: every finite subset of $L_2$
of cardinality $\le \alpha$ is contained in some $W\in \sW$.) Let $\sV$ be the open cover of $L_1$ consisting
of the sets $\pi^{-1}(W)$ for $W\in \sW$. Then $\sV$ is an open $\alpha$-cover of $L_1$.
\end{expl}

Let $\pi\co E\to L$ be a covering space. In \cite{confcover} we introduced $\config(\pi)$ and a forgetful map
from $\config(\pi)$ to $\config(E)$ over $N\fin$. (These definitions relied on the ``particle'' models for $\config(E)$ and
$\config(L)$,
but it is easy to adapt them to the Riemannian patch models in the case where $L$ has a complete Riemannian metric.)
Here we need similar definitions for unordered configuration categories.

\begin{defn} \label{defn-configpi}
Let $\pi \co E \to L$ be a covering space where $L$ is a (smooth, complete) Riemannian manifold. Then $E$ inherits a Riemannian metric from $L$
which makes $\pi$ into a local isometry. The topological poset $\cconfig(\pi)$ has as objects pairs $(V, W)$ where $W$ is an object of
$\cconfig(L)$, i.e., a multipatch in $L$, and $V$ is a union of finitely many connected components of $\pi^{-1}(W)$,
therefore a multipatch in $E$. \emph{Condition}: every connected component of $W$ is the image under $\pi$ of some connected component of $V$.
A morphism in $\cconfig(\pi)$ is a commutative square
\[
\xymatrix@R=16pt{
V_0 \ar[d]^\pi & V_1 \ar[d]^\pi \ar[l] \\
W_0 & W_1 \ar[l]
}
\]
There are forgetful maps $\cconfig(\pi) \to \cconfig(L)$ and $\cconfig(\pi) \to \cconfig(E)$.
Both of these are maps over $\cfin$, but two distinct forgetful maps
from $\cconfig(\pi)$ to $\cfin$ are involved. One of these takes an object $(V,W)$ to $\pi_0W$, and the other takes
the same object to $\pi_0V$. From the definition it is rather clear that $\cconfig(\pi)$ is a complete Segal space.
\end{defn}

\begin{defn}
The definition of $\config(\pi)$, Riemannian patch model, is similar to definition~\ref{defn-configpi} but we spell it out (nearly) in order
to make a few useful observations. An object of $\config(\pi)$ is an object $(V,W)$ of $\cconfig(\pi)$ together
with bijections $\uli k\to \pi_0V$ and $\uli\ell\to \pi_0W$ for apropriate $k$ and $\ell$. \emph{Condition}: the resulting map $\uli k \to
\uli\ell$ is selfic. (See \cite[\S2]{confcover}.) The precise meaning of \emph{selfic} is not important. The important point
is that it makes the bijection $\uli\ell\to \pi_0W$ redundant; it is determined by the bijection $\uli k\to \pi_0V$.
As a consequence, the commutative square
\[
\xymatrix@R=16pt{
\config(\pi) \ar[d] \ar[r] &  \cconfig(\pi) \ar[d] \\
\fin  \ar[r] &  \cfin
}
\]
is homotopy cartesian. (The vertical arrows are the forgetful maps taking an object $(V,W)$ to
the set $\pi_0(V)$, with resp.~without total ordering.) This confirms something which we already observed in \cite{confcover},
that $\config(\pi)$ is fiberwise complete over $\fin$ wrt the vertical arrow of the left-hand column. ---  There are forgetful maps
$\config(\pi)\to \config(L)$ and $\config(\pi)\to \config(E)$, defined much like the forgetful
maps $\cconfig(\pi)\to \cconfig(L)$ and $\cconfig(\pi)\to \cconfig(E)$
\end{defn}

\begin{defn} \label{defn-models5} We describe a model $X$ for $\cfin$, a complete Segal space. (To be consistent
with \cite{confpres} we should go for something like  $\cconfig(\RR^\infty)$.
This is alright as far as the object space is concerned, but then the morphism spaces contain a lot of useless information.)
Let $\sA$ be the category whose objects are the standard finite sets $\uli k$
for $k\ge 0$, with all maps between them as morphisms. For $r\ge 0$ let $\sV_r$ be the groupoid whose objects are the
functors from $[r]^\op$ to $\sA$. The morphisms in $\sV_r$ are the natural isomorphisms between such functors.
Then $[r] \mapsto B\sV_r$ is a simplicial space $X$, and it is a model for $\cfin$.
\end{defn}

\medskip
Often we can avoid specific models for $\cfin$ altogether, e.g., by relying on definition~\ref{defn-leftfib}
and lemma~\ref{lem-classy} just below.
\begin{defn} \label{defn-leftfib} A left fibration with \emph{discrete} fibers on a simplicial space $X$
is a map of simplicial spaces $p\co E\to X$ with the following properties:
\begin{itemize}
\item[-] the map $p_0\co E_0\to X_0$ is a fibration with discrete fibers (a covering projection);
\item[-] for every $r>0$, the commutative square
\[
\xymatrix@C=60pt{ E_r \ar[r]^{(d_0)^r} \ar[d]^-{p_r} & E_0 \ar[d]^-{p_0}  \\
X_r \ar[r]^{(d_0)^r}  & X_0
}
\]
is a (strict) pullback square of spaces.
\end{itemize}
Here $(d_0)^r$ is the ``ultimate source'' operator.

This can also be rephrased as follows.
A left fibration with discrete fibers on a simplicial space $X$ is a
covering projection $\pi\co E_0\to X_0$ together with:
\begin{itemize}
\item[-] for every $y\in X_1$, a map of sets $\psi(y)\co \pi^{-1}(d_0y)\to \pi^{-1}(d_1y)$
which depends continuously on $y$;
\item[-] such that, for every $z\in X_2$ we have $\psi(d_1z)=\psi(d_2z)\circ \psi(d_0z)$.
\end{itemize}
\end{defn}

\begin{expl} \label{expl-magicleft} Let $X$ be the simplicial space of definition~\ref{defn-models5}.
Now we can make a left fibration $E\to X$ with discrete fibers as follows. Define $\bar\sA$ like $\sA$
in definition~\ref{defn-models5}, but for the objects take the standard finite sets with a choice of
base element, and as morphisms allow all maps which respect the choice of base element. Let $\bar\sV_r$
be the groupoid whose objects are the functors from $[r]^\op$ to $\bar\sA$. The morphisms
in $\bar\sV_r$ are the natural isomorphisms between such functors. Then
\[ [r] \mapsto B\bar\sV_r  \]
is a simplicial space $E$. There is a forgetful simplicial map $E\to X$. It is a left fibration with discrete fibers.
\end{expl}

\begin{lem} \label{lem-classy} Among simplicial spaces, $\cfin$ is a classifying object for left
fibrations with finite fibers. \qed
\end{lem}
The lemma speaks of left fibrations with fibers which are not only discrete, but in addition
finite. A more elaborate formulation is that for cofibrant
simplicial spaces $X$, there is a zigzag of natural weak equivalences relating the mapping space
$\rmap(X,\cfin)$ to the classifying space of the (discrete) groupoid whose objects are the left fibrations
with finite fibers on $X$. (The ``cofibrant'' condition on $X$ can be weakened. Degreewise
paracompact is enough.) Of course the lemma also wants to say that $\cfin$ carries a universal
left fibration with finite fibers. We have already seen that in example~\ref{expl-magicleft}.

\medskip
In the case of the simplicial space $\cconfig(L)$, the standard reference map from $\cconfig(L)$ to $\cfin$
corresponds to a left fibration $\cconfig^b(L)\to \cconfig(L)$ which we can imagine as follows.
In the (Riemannian) multipatch model, $\cconfig(L)$ is the nerve of a certain topological poset of
multipatches in $L$. Define $\cconfig^b(L)$ similarly, but use multipatches with one distinguished
component (and morphisms, i.e. inclusions, respecting this). There is a forgetful map
$\cconfig^b(L)\to \cconfig(L)$.

\medskip
Let $X$ be a simplicial space. We want a practical method for finding isomorphisms between left fibrations $p\co E\to X$
and $q\co F\to X$ with discrete fibers.

\begin{lem} \label{lem-itsenough} To construct an isomorphism $E\to F$ over $X$, it suffices to
construct an isomorphism $h\co  E_0 \to F_0$ of covering spaces over $X_0$ and to verify that for
every $g\in X_1$ the diagram
\[
\xymatrix{ p^{-1}(d_0g)\cong p^{-1}(g) \ar[r]^-{d_1} \ar[d]^-h & p^{-1}(d_1g) \ar[d]^-h \\
q^{-1}(d_0g) \cong q^{-1}(g) \ar[r]^-{d_1} & q^{-1}(d_1g)
}
\]
commutes. \qed
\end{lem}
The method of lemma~\ref{lem-itsenough} can be simplified under a condition on $p\co E\to X$. Suppose that for every
$v\in E_0$ there exists $w\in E_1$ such that $d_1w$ is in the path component of $v$ and $p^{-1}(p(w))$
has only one element, inevitably $w$. Then it is enough to verify commutativity of the diagram in
lemma~\ref{lem-itsenough} in all cases where $g\in X_1$ satisfies $|p^{-1}(g)|=1$, or equivalently,
$|p^{-1}(d_0g)|=1$.

\section{Homotopical decisions}  \label{sec-homodels}
\begin{notn} \label{notn-mapt} Where we encounter mapping spaces, more often than not they come as simplicial sets because it is
tiresome to give them any other status. We will write $\map(-,-)$ for these. There are a few cases where
mapping spaces as topological spaces are more appropriate. For such cases we have the
notation $\mapt(-,-)$. (It is unlikely that we can be consistent in using it.)
\emph{Example:} let $Y$ be a simplicial space (each $Y_n$ is a space, a.k.a. object of $\topcat$).
Then $\mapt(\Delta[n],Y)$ is a space homeomorphic to $Y_n$, whereas $\map(\Delta[n],Y)$ is a
simplicial set isomorphic to the singular simplicial set of $Y_n$.
\end{notn}

\begin{defn} \label{defn-models4} The current preference for a model category structure on $\topcat$
is the mixed model structure \cite{Cole2006} in which the objects of $\topcat$ are the compactly generated weak Hausdorff
spaces \cite{Vogt}, the categorical weak equivalences are the weak homotopy equivalences and the
categorical fibrations are the Hurewicz fibrations. This was also used in \cite{confcover}.
The current preference for a model category structure
on the category of simplicial spaces is the corresponding Reedy model structure. 
One of the main reasons for this is that $\config(L)$ and $\cconfig(L)$ are already Reedy cofibrant;
this is a statement about the degeneracy operators in $\config(L)$ and $\cconfig(L)$ only.
It is clear what \emph{Reedy fibrant} means for a simplicial space $X$. It means that
for every $n\ge 1$, the map
\[   X_n \quad \lra \quad \match_mX~~:=\!\!\!\! \limsub{\twosub{f\co [m]\to [n] \textup{ in }\Delta_\inj}{m<n}} X_m \]
induced by the various $f$ is a fibration. More generally, a map
of simplicial spaces $X\to Y$ is a Reedy fibration if for every $m\ge 0$ the map
from $X_m$ to the limit of $(Y_m\to \match_mY \leftarrow \match_mX)$ determined by
the commutative square
\[
\xymatrix@R=15pt{
X_m \ar[r] \ar[d] & \match_mX \ar[d] \\
Y_m \ar[r] & \match_mY
}
\]
is a fibration. --- It is less clear what
could be meant by the \emph{standard} (functorial) Reedy fibrant replacement $\varphi X$ of a simplicial space $X$.
But there is such a thing. The formula is
\[ (\varphi X)_n := \holimsub{[m]\to [n] \textup{ in }\Delta_\inj} X_m  \]
where $[m]\to [n]$ runs over the comma category $(\Delta_\inj\!\downarrow\![n])$.
The simplicial operators are defined as follows. (See also \cite[\S3]{confcover}.)
Let $\vep_n$ be the forgetful functor
from the comma category $(\Delta_\inj\!\downarrow\![n])$ to $\Delta$.
A morphism $f\co [m]\to [n]$ in $\Delta$ induces a functor
\[  f_!\co (\Delta_\inj\!\downarrow\![m])\to (\Delta_\inj\!\downarrow\![n])  \]
which is the composition of $f\circ \co (\Delta_\inj\!\downarrow\![m])\to (\Delta\!\downarrow\![n])$
with the left adjoint of the inclusion $(\Delta_\inj\!\downarrow\![n])\to (\Delta\!\downarrow\![n])$.
There is a preferred and obvious natural transformation
$\vep_{m} \to \vep_{n} f_!$
which induces contravariantly $X \vep_{n} f_! \to X \vep_{m}$ and then
\[ \rule{9mm}{0mm} (\varphi X)_{n}= \holim (X\vep_{n}) \to \holim (X\vep_{n}f_!) \to \holim (X\vep_{m}) =(\varphi X)_{m} \,. \]
\end{defn}

\begin{prop} \label{prop-phifibration} If $f\co X\to Y$ is a map of simplicial spaces which is a fibration in the projective structure,
i.e., a levelwise fibration, then the induced map $\varphi f$ from $\varphi X$ to $\varphi Y$ is a Reedy fibration.
\end{prop}
\proof In this proof, a commutative square of spaces
\[
\xymatrix@R=15pt{
A \ar[r] \ar[d] & \ar[d] B \\
C \ar[r] & D
}
\]
will be called \emph{$\tau$-fibrant} if the tautological map from $A$ to the homotopy
limit of $C\to D\leftarrow B$ is a fibration in $\topcat$. ---
We need to show that for every $n\ge 0$ the commutative square
\[
\xymatrix@R=20pt{
{\holimsub{[m]\to [n]} X_m} \ar[r] \ar[d] & {\limsub{\twosub{[m]\to [n]}{m<n}} \holimsub{[\ell]\to [m]} X_m} \ar[d] \\
{\holimsub{[m]\to [n]} Y_m} \ar[r] & {\limsub{\twosub{[m]\to [n]}{m<n}} \holimsub{[\ell]\to [m]} Y_m}
}
\]
is $\tau$-fibrant. The square simplifies immediately to
\begin{equation} \label{eqn-reedysimplif}
\begin{aligned}
\xymatrix@R=20pt{
{\holimsub{[m]\to [n]} X_m} \ar[r] \ar[d] & {\holimsub{\twosub{[m]\to [n]}{m<n}} X_m} \ar[d] \\
{\holimsub{[m]\to [n]} Y_m} \ar[r] & {\holimsub{\twosub{[m]\to [n]}{m<n}} Y_m}
}
\end{aligned}
\end{equation}
We begin with the following observation. If in a commutative diagram of spaces
\[
\xymatrix@R=15pt{
A_0 \ar[r] \ar[d] & B_0 \ar[d]  & \ar[l] C_0 \ar[d] \\
A_1 \ar[r]  & B_1 & \ar[l] C_1
}
\]
the vertical arrows are fibrations in $\topcat$, then the induced map between the homotopy
limits of the rows is again a fibration in $\topcat$. This observation applies to the diagram
\[
\xymatrix{
X_m \ar[r] \ar[d]_-{f_m} &  \ar[d]^-{f_*}
{\holimsub{\twosub{[m]\to [n]}{m<n}} X_m} & \ar[d]^-{\id} \ar[l]_-{\id} {\holimsub{\twosub{[m]\to [n]}{m<n}} X_m} \\
Y_m \ar[r] &  {\holimsub{\twosub{[m]\to [n]}{m<n}} Y_m} & \ar[l]_-{f_*} {\holimsub{\twosub{[m]\to [n]}{m<n}} X_m}
}
\]
so that the map between the homotopy limits of the rows is a fibration. But this is exactly the tautological map
from the initial term in~\eqref{eqn-reedysimplif} to the homotopy limit of the three-term diagram obtained by deleting the
initial term. \qed

\begin{cor} \label{cor-EZfibrant} Let $g\co X\to Y$ be a map of simplicial spaces. Let $\varphi_Y X$ be the levelwise
homotopy limit of
\[   \varphi X \xrightarrow{\varphi(g)} \varphi Y \xleftarrow{\textup{incl.}} Y\,.\]
Then the projection $\varphi_Y X\to Y$ is a Reedy fibration.
\end{cor}
\proof Let $X^h$ be the homotopy limit of $X \xrightarrow{g} Y \xleftarrow{\id} Y$ (degreewise)
and let
\[ g^h\co X^h\to Y \]
be the projection. This is a fibration in the projective model structure,
so that $\varphi(g^h)\co \varphi(X^h) \to \varphi Y$ is a Reedy fibration by proposition~\ref{prop-phifibration}.
Pulling this back along the inclusion $Y\to \varphi Y$ gives $\varphi_Y X\to Y$, which is therefore
also a Reedy fibration. \qed

\begin{lem} \label{lem-reedyone}
The Reedy fibrant replacement $\varphi$ has a left adjoint $\varphi^\lad$.
\end{lem}
\proof This follows
from the adjoint functor theorem, but we give an explicit description.
For $n\ge 0$ and a (constant) space $Z$, the simplicial space $Z\times \Delta[n]$ is a co-representing object
for the functor $X\mapsto \mor_\topcat(Z,X_n)$. It follows that $\varphi^\lad(Z\times \Delta[n])$ is a co-representing
object for $X\mapsto \mor_\topcat(Z,(\varphi X)_n)$. Therefore
\begin{equation} \label{eqn-adphidef}  \varphi^\lad(Z\times\Delta[n]) = \hocolimsub{[m]\to [n] \textup{ in }\Delta_\inj} Z\times\Delta[m]
~~=~~Z\times \varphi^\lad(\Delta[n])
\end{equation}
where the homotopy colimit is taken (degreewise) in the category of simplicial spaces.
It is still easy to understand how a morphism $u\co [n_0]\to [n_1]$ in $\Delta$ induces a
simplicial map
\[ \varphi^\lad(\Delta[n_0])\to \varphi^\lad(\Delta[n_1]). \]
(If $v\co [m]\to [n_0]$ is a morphism in $\Delta_\inj$, then $uv$ need not be in $\Delta_\inj$
but it has a unique factorization $[m]\to [m']\to [n_1]$ in $\Delta$ where $[m]\to [m']$ is surjective
and $[m']\to [n_1]$ is injective.) An arbitrary simplicial space $X$ can be written in the coend form
\[  \Big(\coprod_n X_n\times \Delta[n]\Big)\Big/\textup{relations}.  \]
Therefore
\begin{equation} \label{eqn-adphi} \varphi^\lad X=
\Big(\coprod_n X_n\times \varphi^\lad(\Delta[n])\Big)\Big/\textup{relations}. \end{equation}
The relations are the familiar ones,
\[ (f^*a,b)\simeq (a,f_*b) \]
for $(a,b)\in X_n\times \varphi^\lad(\Delta[m])$ and $f\co [m]\to [n]$ in $\Delta$. \qed

\medskip
\emph{Remark.} It is correct to say that $\varphi$ has little effect in low degrees, e.g., no effect in degree 0,
but it is wrong to think that $\varphi^\lad$ has little effect in low degrees. For example, $\Delta[1]$ in degree 0
is discrete (with two elements), but $\varphi^\lad\Delta[1]$ in degree 0 is the disjoint union of two intervals.
Related to this observation: it is always true that $\varphi(\skel_k X)\cong \skel_k(\varphi X)$, but for most $k$
and simplicial spaces $X$ it would be wrong to claim $\varphi^\lad(\skel_k X)\cong \skel_k(\varphi^\lad X)$.

\begin{lem} \label{lem-reedyfour} If $X$ is a Reedy cofibrant simplicial space, then
the map
\[ v\co \varphi^\lad X\to X \]
adjoint to the inclusion $u\co X\to \varphi X$
is a weak equivalence.

\end{lem}
\proof Let $Y$ be a Reedy fibrant simplicial space. There is a commutative triangle
\[
\xymatrix@C=0pt{
\map(X,Y) \ar[rr]^-{\circ v} \ar[dr]^-{u\circ} && \map(\varphi^\lad X,Y)  \\
& \map(X,\varphi Y) \ar[ur]^-\cong
}
\]
in which the arrow labeled $u\circ$ is a weak equivalence (because $u$ is a weak equivalence of Reedy fibrant
objects and $X$ is Reedy cofibrant). It follows that the horizontal arrow is always a weak equivalence, and
it is a map between spaces (here simplicial sets) which can call themselves derived mapping spaces. \qed

\begin{prop} \label{prop-reedyfive} Let $X$ and $Y$ be simplicial spaces and let
$f\co X\to Y$ be a Reedy cofibration. Then $\varphi^\lad f\co \varphi^\lad X\to \varphi^\lad Y$ is a
projective cofibration. In particular, if $Y$ is a Reedy cofibrant simplicial space, then $\varphi^\lad Y$ is projectively cofibrant.
\end{prop}
\proof Let $g\co P\to Q$ be a map of simpiicial
spaces which is an acyclic fibration in the projective structure, i.e., a levelwise fibration.
We are supposed to show that $g$ has the appropriate lifting property with respect to $\varphi^\lad f$.
By adjunction, this is equivalent to showing that $\varphi g$ has the appropriate lifting
property with respect to $f$. But it does have that property because $\varphi g$ is a Reedy
fibration by proposition~\ref{prop-phifibration} and $f$ is a Reedy cofibration by assumption. \qed

\medskip
Frequently we are confronted with a pile-up of fibrant replacements. In such a case
the following lemma can provide relief.
\begin{lem} \label{lem-phimonad} The fibrant replacement $X\mapsto \varphi X$ and the preferred natural inclusion
$\iota_X\co X\to \varphi X$ together admit the structure of a monad.
\end{lem}

\proof We have to find a natural transformation $\mu\co \varphi\circ\varphi\to \varphi$
satisfying associativity and having $\iota$ as a two-sided unit. Let $X$ be a simplicial space and choose $n\ge 0$. Then
\[ (\varphi\varphi X)_n \cong \mapt(\Delta[n],\varphi\varphi X)\cong
\mapt(\varphi^\lad\Delta[n],\varphi X) \cong \mapt(\varphi^\lad\varphi^\lad\Delta[n],X). \]
Therefore we turn our attention to $\varphi^\lad\varphi^\lad\Delta[n]$. Using~\eqref{eqn-adphi} we obtain
\[
\begin{array}{rcl}
\varphi^\lad\varphi^\lad\Delta[n]  & \cong  &
\coprod_k \left((\varphi^\lad\Delta[n])_k\times\varphi^\lad\Delta[k]\right)\Big/\textup{relations} \\
& = & \rule{0mm}{4.5mm}
\coprod_k \Big(\hocolimsub{[m]\to [n]} (\Delta[m])_k~\times \hocolimsub{[j]\to[k]} \Delta[j]\Big)\Big/ \textup{relations} \\
& = &
\coprod_k \Big(\hocolimsub{[m]\to [n]} \mor_\Delta([k],[m])~\times \hocolimsub{[j]\to[k]} \Delta[j]\Big)\Big/ \textup{relations} \\
& \cong  & \hocolimsub{[m]\to [n]}
\Big(\coprod_k \big( \mor_\Delta([k],[m])~\times \hocolimsub{[j]\to[k]} \Delta[j]\big)\Big/ \textup{relations}\Big) \\
& \cong & \hocolimsub{[m]\to [n]} \hocolimsub{[j]\to[m]} \Delta[j]
\end{array}
\]
where $[m]\to[n]$, $[j]\to [k]$ and $[j]\to [m]$ denote morphisms in $\Delta_\inj$. This simplifies some more,
so that we have
\begin{equation}  \label{eqn-illumin}
\varphi^\lad\varphi^\lad\Delta[n] \quad \cong   \hocolimsub{[j]\to [m]\to [n]} \Delta[j].
\end{equation}
In more detail, there is a category $\sB_n$
whose objects are diagrams $[j]\to [m]\to [n]$ in $\Delta_\inj$, with fixed $[n]$. A morphism is a commutative diagram
\[
\xymatrix@R=12pt{
[j_0]\ar[r] \ar[d] & [m_0] \ar[r] \ar[d] & [n] \ar[d]^-=  \\
[j_1] \ar[r] & [m_1] \ar[r] & [n]
}
\]
Then there is a functor $Q_n$ from $\sB_n$ to simplicial sets (simplicial spaces) taking $[j]\to [m]\to [n]$ to $\Delta[j]$.
We have calculated $\varphi^\lad\varphi^\lad\Delta[n]\cong \hocolim~Q_n$.

This calculation is natural in $[n]$. Namely, a morphism $f\co [m]\to [n]$ in $\Delta$, not necessarily injective,
induces a functor $f_* \co\sB_m\to \sB_n$
(details left to the reader) and a natural transformation $Q_m\to Q_nf_*$. These in turn determine a map
\[   \hocolim~Q_m \to \hocolim~Q_n. \]
Using~\eqref{eqn-illumin} and~\eqref{eqn-adphidef}, it is easy to produce a natural transformation
\begin{equation} \label{eqn-relief}  \varphi^\lad\Delta[n]=\!\!\hocolimsub{[j]\to [n]}\Delta[j] \quad\lra~~
\hocolimsub{[j]\to[m]\to[n]}\Delta[j]= \varphi^\lad\varphi^\lad\Delta[n]  \end{equation}
(natural in $[n]$ as an object of $\Delta$). Namely, let $\sA_n$ be the category whose objects are diagrams
$[j]\to [n]$ in $\Delta_\inj$, with fixed $[n]$, etc., so that $\hocolim_{[j]\to [n]}\Delta[j]$
is the hocolim of the functor $P_n$ on $\sA_n$ given by $P_n([j]\to[n])=\Delta[j]$. We need functors
$u_n\co \sA_n\to \sB_n$ and we define them by
\[  \big([j]\to [n]\big) \mapsto \big([j]\xrightarrow{=} [j]\to [n]\big). \]
(This definition of $u_n$ has the desirable consequence $u_n f_*=f_*u_m$ for morphisms $f\co [m]\to [n]$
in $\Delta$.) Then $Q_n u_n=P_n$ and this leads to a map from $\hocolim~P_n$ to $\hocolim~Q_n$ which implements~\eqref{eqn-relief}.
The resulting map
\[   (\varphi\varphi X)_n \lra (\varphi X)_n  \]
(for a simplicial space $X$) is covariantly natural in $X$ and contravariantly natural in $[n]$, object of $\Delta$,
by inspection. Therefore we have constructed $\mu\co \varphi\circ\varphi \to \varphi$. The associativity
property is easily verified on the basis of
\[  \varphi^\lad\varphi^\lad\varphi^\lad\Delta[n] \quad \cong   \hocolimsub{[i]\to [j]\to [m]\to [n]} \Delta[j] \]
(where $[i]\to [j]$, $[j]\to [m]$ and $[m]\to [n]$ denote morphisms in $\Delta_\inj$). This is a calculation
similar to~\eqref{eqn-illumin}, the details of which are left to the reader. The unit properties are more obvious.
Keep in mind that $\iota\co \id\to \varphi$ is adjoint to a natural transformation
$\varphi^\lad \to \id$ which is induced by certain maps $\varphi^\lad\Delta[n]\to \Delta[n]$, for $n\ge 0$.
These maps are just the usual maps from a homotopy colimit to the corresponding colimit,
$\hocolim_{[j]\to[n]}\Delta[j]\to \colim_{[j]\to[n]}\Delta[j]\cong \Delta[n]$. \qed

\begin{cor} \label{cor-phiadcomonad} The functor $X\mapsto \varphi^\lad X$ and the natural projection
$\varphi^\lad X\to X$ together admit the structure of a comonad. \qed
\end{cor}

\begin{rem} \label{rem-imagine} Let $X$ and $Y$ be simplicial spaces. How should we imagine a map from $X$ to $\varphi Y$,
or equivalently, a map from $\varphi^\lad X$ to $Y$~? Understanding this can be more useful than having a good
idea of what $\varphi Y$ is, or what $\varphi^\lad X$ is. A map from $X$ to $\varphi Y$ is a package which provides, for every $r\ge 0$
and every diagram
\[  \sD~:=~~[k_0] \xleftarrow{g_1} [k_1]\xleftarrow{g_2} \cdots \xleftarrow{g_{r-1}} [k_{r-1}]\xleftarrow{g_r} [k_r] \]
in $\Delta_\inj$\,, a map $H_\sD\co X_{k_0} \times \Delta^r \lra Y_{k_r}$.
These maps $H_\sD$ are subject to two naturality conditions. Firstly, for any morphism $u\co [q]\to [r]$ in $\Delta_\inj$ the diagram
\[
\xymatrix@C=35pt{
X_{k_0} \times \Delta^r \ar[r]^-{H_\sD} & Y_{k_r} \\
X_{k_0}\times \Delta^q \ar[u]^-{\id\times u_*} \ar[d]_-{(g_1g_2\cdots g_{u(0)})^*\times\id}  & \\
X_{k_{u(0)}}\times \Delta^q \ar[r]^-{H_{u^*\sD}} & \ar[uu]_-{(g_{u(r)+1}g_{u(r)}\cdots g_{r-1}g_r)^*} Y_{k_{u(r)}}
}
\]
commutes. Secondly, if in a commutative diagram
\[
\xymatrix@M=5pt@C=15pt@R=15pt{
[k_0] \ar[d]^-{p_0} & \ar[l] [k_1] \ar[d]^-{p_1} & \ar[l] \cdots & \ar[l] [k_{r-1}] \ar[d]^-{p_{r-1}} & \ar[l] [k_r] \ar[d]^-{p_r}  \\
[\ell_0] & \ar[l] [\ell_1] & \ar[l] \cdots & \ar[l] [\ell_{r-1}]  & \ar[l] [\ell_r]
}
\]
in $\Delta$, all horizontal arrows are injective and all vertical arrows are surjective (top row $\sD$, bottom row $\sE$),
then
\[
\xymatrix{
X_{k_0} \times \Delta^r \ar[r]^-{H_\sD} & Y_{k_r} \\
X_{\ell_0} \times \Delta^r \ar[u]_-{p_0^*\times\id} \ar[r]^-{H_\sE} & \ar[u]_-{p_r^*} Y_{\ell_r}
}
\]
commutes.
\end{rem}

\section{Joins in a simplicial setting} \label{sec-joins}
\begin{defn} \label{defn-joingeo} The (geometric) join $Y*Z$ of two spaces $Y$ and $Z$ is
the homotopy colimit, aka homotopy pushout, of
\[  Y \xleftarrow{\textup{proj.}} Y\times Z \xrightarrow{\textup{proj.}} Z\,.  \]
It is often thought of as a quotient of the disjoint union
$Y\sqcup(Y\times Z\times \Delta^1)\sqcup Z$.
If $Y$ and $Z$ are both nonempty, it can of course be
viewed as a quotient of $Y\times Z\times \Delta^1$. In that case elements can be labeled $(sy,tz)$
where $y\in Y$, $z\in Z$ and $(s,t)\in \Delta^1$. Then the relations are $(sy_1,tz_1)\sim (sy_2,tz_2)$
iff $(s,t)=(1,0)$ and $y_1=y_2$ or $(s,t)=(0,1)$ and $z_1=z_2$.

The geometric join can pass for a monoidal product in $\topcat$.
The space $\emptyset$ is a two-sided unit for the join.
\end{defn}

\begin{expl} For integers $p,q,n\ge 0$ such that $p+q+1=n$,
there is a preferred homeomorphism $\Delta^p*\Delta^q\to \Delta^n$ given by
\[  ((y_0,y_1,...,y_p),(z_0,z_1,\dots,z_q),(s,t)) \mapsto (sy_0,\dots,sy_p,tz_0,\dots,tz_q). \]
\end{expl}

\begin{defn} \label{defn-joinsi}
The \emph{join} of two simplicial spaces $A$ and $B$ is a simplicial space $A*B$. The definition begins
with a few auxiliary conventions, $[-1]:=\emptyset$ and $A_{-1}:=\star=:B_{-1}$. Then
\[  (A*B)_n:= \coprod_{\twosub{p,q\ge -1}{p+q+1=n}} A_p\times B_q \]
for $n\ge 0$. The simplicial operator corresponding to a morphism $f\co [m]\to [n]$ in $\Delta$ acts on
$A_p\times B_q$ (where $p+q+1=n$) as
follows. Determine $m_1\ge -1$ and $m_2\ge -1$ in such a way that $m_1+m_2+1=m$ and $f(x)\le p$ if and only if
$x\le m_1$. Then there are unique monotone maps $f_1\co [m_1]\to [p]$ and $f_2\co [m_2]\to [q]$ such that
$f_1(x)=f(x)$ for $x\in [m_1]$ and $f_2(x)+p+1=f(x+m_1+1)$ for $x\in [m_2]$. For $(z_1,z_2)\in A_p\times B_q$ let
$f^*(z_1,z_2):= (f_1^*(z_1),f_2^*(z_2))\in A_{m_1}\times B_{m_2} \subset (A*B)_m$.

There are preferred inclusions $A\to A*B$ and $B\to A*B$. (The first of these in degree $n$ is given by
$A_n\cong A_n\times B_{-1}\subset (A*B)_n$\,.) There is a preferred simplical map
$A*B \lra \Delta[1]$ taking the summand $A_p\times B_q$ of $(A*B)_n$ to the monotone map $g\co [n]\to [1]$
which has $g(x)=0$ for $x\le p$ and $g(x)=1$ for $x>p$.
\end{defn}

\begin{expl} The join of $\Delta[p]$ and $\Delta[q]$ is isomorphic as a simplicial set to
$\Delta[p+q+1]$. This is an easy exercise.
\end{expl}
The join as in definition~\ref{defn-joinsi} is important to us (particularly in section~\ref{sec-pseudomet}) because
of the following observation.
\begin{prop} If the simplicial spaces $A$ and $B$ are Segal spaces, then $A*B$ is also a Segal space. \qed
\end{prop}
Looking at some special cases, we can be a little more precise. Suppose for example
that $\sC$ and $\sD$
are (discrete) small categories. Let $N\sC$ and $N\sD$ be their nerves. Then $N\sC*N\sD$ is isomorphic to
the nerve of a third category $\sK$, defined as follows. The set of objects $\ob(\sK)$ is the
disjoint union of $\ob(\sC)$ and $\ob(\sD)$. For $x,y\in \ob(\sK)$ we let
\[   \mor_\sK(x,y) := \left\{\begin{array}{cl} \mor_\sC(x,y) & \textup{ if }x,y\in \ob(\sC) \\
\mor_\sD(x,y) & \textup{ if }x,y\in \ob(\sD) \\
\star & \textup{ if $x\in\ob(\sD)$ and $y\in \ob(\sC)$} \\
\emptyset & \textup{ if $x\in\ob(\sC)$ and $y\in \ob(\sD).$}
\end{array} \right.
\]
A similar example: suppose that $\sP$ and $\sQ$ are topological posets. Then $N\sP*N\sQ$ is isomorphic
(as a simplicial space) to the nerve of another topological poset. As a space, this is $\sP\sqcup \sQ$,
with the order relation which has $x\le y$ if and only if
\begin{itemize}
\item[] $x,y\in \sP$ and $x\le y$ in $\sP$, or $x,y\in \sQ$ and $x\le y$ in $\sQ$, or $x\in\sQ$ and $y\in \sP$.
\end{itemize}

\begin{defn} \label{defn-fracdi} Let $C$ be a simplicial space. There is a map of simplicial spaces
\[  \kappa\co C\times \Delta[1]\lra C*C  \]
as follows. Given $x\in C_n$ and monotone $f\co [n]\to [1]$, let $p$ be the maximum of the $x\in [n]$ such that
$f(x)=0$ (and let $p=-1$ if there is no such $x$). Let $q:=n-p-1$. Let $g\co [p]\to [n]$ be the
inclusion and let $h\co [q]\to [n]$ be defined by $h(x)=x+p+1$. Then
\[  \kappa(x,f):=(g^*(x),h^*(x))\in C_p\times C_q\subset (C*C)_n\,. \]
We may call $\kappa$ the \emph{separation diagonal} if a name is needed.
\end{defn}
If $C$ is the nerve of a topological poset, $C=N\sP$, then $\kappa$ in definition~\ref{defn-fracdi}
is an embedding of simplicial spaces with closed image. (We assume that the order relation
is a closed subset of $\sP\times \sP$.)

\medskip
The fibrant replacement $\varphi$ of section~\ref{sec-homodels} interacts in an interesting way with the join
of definition~\ref{defn-joinsi}. To explain this we begin with the left adjoint $\varphi^\lad$.
Suppose that $p,q,n$ are integers such that $p+q+1=n$, where $n\ge 0$ and $p,q\ge -1$. Then there is a preferred map
of simplicial spaces
\begin{equation} \label{eqn-phijoin}
   \varphi^\lad\Delta[n] \lra \Delta[p]*\varphi^\lad(\Delta[q]) \end{equation}
as follows.
\[
\xymatrix@C=20pt@R=10pt@M=5pt{
\varphi^\lad\Delta[n] \ar[d]^-= &   \Delta[p]\,*\,\varphi^\lad\Delta[q]     \\
 {\hocolimsub{[m]\to[n]} \Delta[m]} \ar[d]^-\cong  &
{\Delta[p]*\big(\!\!\hocolimsub{[m_2]\to[q]} \Delta[m_2]\big)} \ar[u]^-= \\
{\hocolimsub{\twosub{[m_1]\to[p]}{[m_2]\to[q]}} \Delta[m_1]~*~\Delta[m_2]} \ar[r] & {\hocolimsub{[m_2]\to[q]} \Delta[p]~*~\Delta[m_2]} \ar[u]
}
\]
(All maps which appear here under hocolim signs are monotone injective. Beware that $p,q,m_1,m_2$
can take the value $-1$. As before, $[-1]$ means $\emptyset$, but we must also agree that
$\Delta[-1]=\emptyset$.)

\begin{defn} \label{defn-phijoin} Let $A$ and $B$ be simplicial spaces. There is a preferred map
\[  A*\varphi(B) \lra \varphi(A*B) \]
as follows. An element of $A*\varphi(B)$ in degree $n$ is a pair of simplicial maps
$\Delta[p]\to A$, $\varphi^\lad\Delta[q]\to B$ for some $p,q\ge -1$ such that $p+q+1=n$.
These determine a simplicial map
\[  \Delta[p]~*~\varphi^\lad\Delta[q]\lra A*B \]
which we may pre-compose with the map~\eqref{eqn-phijoin}. The composition is a map
from $\varphi^\lad\Delta[n]$ to $A*B$. This is tantamount to an element of
$\varphi(A*B)$ in degree $n$.
\end{defn}

\begin{lem} \label{lem-asinus} Let $\sK$ be a topological poset. For $w\in \sK$ let $\sK(w)\subset \sK$ be the
topological poset consisting of all $v\in \sK$ which satisfy $v\le w$. There is a strict pullback square
of simplicial spaces
\[
\xymatrix{
A  \ar[d]_-{\textup{incl.}} \ar[r] &   \varphi(N\sK\times\Delta[1]) \ar[d]^-{\varphi(\textup{separation diag.})} \\
N\sK*\varphi N\sK \ar[r] &  \varphi(N\sK*N\sK)
}
\]
where $A$ is the simplicial subspace of $N\sK*\varphi N\sK$ defined as follows.
\[  A_n:= \coprod_{\twosub{p,q\ge -1}{p+q+1=n}}\!\!\!\!\! A(p,q)\quad\subset
\coprod_{\twosub{p,q\ge -1}{p+q+1=n}}\!\!\!\!\! N\sK_p \times (\varphi N\sK)_q   \]
and $A(p,q)\subset N\sK_p \times (\varphi N\sK)_q$ consists of the pairs $(x,y)$ such that $y\in (\varphi N\sK(x_p))_q$,
where $x_p\in N\sK_0=\sK$ is the ultimate source of $x\in N\sK_p$. \emph{(If $p=-1$ or $q=-1$, then there is no such condition,
i.e., $A(p,q)=N\sK_p \times (\varphi N\sK)_q$.)} \qed
\end{lem}

\section{Another view of conservatization}  \label{sec-fibreal}
Let $\sC$ be a (discrete) small category. Let $X$ be simplicial space and let $X\to N\sC$ be a map
of simplicial spaces. Recall from \cite{GalKoTon} or \cite[Def.8.1]{BoavidaWeiss2018}
that $X\to N\sC$ is \emph{conservative} if (and only if) for every surjective $f\co [m]\to [n]$ in $\Delta$,
the square
\[
\xymatrix@R=16pt{
X_n \ar[r]^-{f^*} \ar[d] & X_m \ar[d] \\
(N\sC)_n \ar[r]^-{f^*} &  (N\sC)_m
}
\]
is homotopy cartesian.   
In \cite{BoavidaWeiss2018} we constructed a functor $\Lambda$ from the category
of simplicial spaces over $N\sC$ to itself which is derived left adjoint to the inclusion of the full subcategory consisting
of the conservative objects. A map of simplicial spaces over $N\sC$ as in
\[
\xymatrix@C=8pt@R=12pt{
W \ar[dr] \ar[rr] && \ar[dl] X \\
& N\sC
}
\]
is a \emph{conservatization map} if $X$ is conservative over $N\sC$ and the map is derived initial among maps
over $N\sC$ from $W$ to a conservative object; equivalently, if the induced map
$\Lambda W \lra \Lambda X$ is a degreewise weak equivalence of simplicial spaces.

\smallskip
Now we want to describe and use another
variant $\Lambda^\fre$ of $\Lambda$, definition~\ref{defn-freLam}. This is inspired by \cite[\S3.3]{Rezk01}.
It is defined under more general circumstances. Briefly, we can replace $N\sC$ by any simplicial space $Y$.
A map $X\to Y$ of simplicial spaces
is \emph{conservative} if (and only if) for every surjective $f\co [m]\to [n]$ in $\Delta$,
the square
\[
\xymatrix@R=16pt{
X_n \ar[r]^-{f^*} \ar[d] & X_m \ar[d] \\
Y_n \ar[r]^-{f^*} &  Y_m
}
\]
is homotopy cartesian.

\medskip
\emph{In this section and from now on, we write \emph{Space} for simplicial set. We
make an effort to avoid ordinary topological spaces. The \emph{realization} of a simplicial Space is a Space
\emph{(coend construction; see \cite{Quillen73}, proof of Thm A)}. By $\Delta[n]$ we mean a certain simplicial Space
which is discrete in every degree.}

\begin{defn} \label{defn-freLam}
Let $X$ be a simplicial Space and let $u\co X\to Y$ be a
map of simplicial Spaces, where $Y$ is degreewise fibrant.
Let $X\langle u\rangle$ be the bisimplicial Space defined by
\[  (\,[m],[n]\,) \mapsto \,\holim\left(\begin{aligned}
\xymatrix@C=20pt@R=15pt{
& \map(\Delta[m],Y) \ar[d] \\
\map(\Delta[m]\times\Delta[n],\varphi X) \ar[r]^-{u_*} & \map(\Delta[m]\times\Delta[n],\varphi Y)
}\end{aligned}\right).
\]
Let $\Lambda^\fre X$ be the simplicial Space obtained from the bisimplicial Space $X\langle u\rangle$
by realization in the second variable $[n]$.
There is a preferred inclusion of $\varphi_YX$ in $\Lambda^\fre Y$ (try $n=0$) and there is a
preferred forgetful map $\Lambda^\fre X\to Y$. The composition of these is
the standard fibrant replacement of $u\co X\to Y$.
\end{defn}
\emph{Remarks.} The Space $\map(\Delta[m],Y)$ is isomorphic to $Y_m$ and we can say that it is contained in $(\varphi Y)_m$.
The vertical arrow within the large parentheses is
induced by the projection from $\Delta[m]\times\Delta[n]$ to $\Delta[m]$.
By definition $\Delta[n]$ is the simplicial set represented by $[n]$,
object of $\Delta$. With our conventions for \emph{nerve}, it is more obviously identified with
the nerve of $[n]^\op$ than with the nerve of $[n]$.

\medskip
Using exponential notation for derived internal hom objects, we see that there is a
homotopy cartesian square of simplicial Spaces
\begin{equation} \label{eqn-expoLam}
\begin{split}
\xymatrix@R=15pt{
X\langle u\rangle_{m,-} \ar[r] \ar[d] & Y_m  \ar[d] \\
X^{\Delta[m]} \ar[r] &  Y^{\Delta[m]}
}
\end{split}
\end{equation}
where $Y_m$ in the upper row is meant as a constant simplicial Space.

\medskip
In the rest of this section we write $\Lambda$ to mean $\Lambda^\fre$.

\begin{prop} \label{prop-Lambdadoesit} $\Lambda X$ is conservative over $Y$.
\end{prop}
\proof Let $f\co [m+1]\to [m]$ be a surjective morphism in $\Delta$; let $t\in[m]\subset [m+1]$ be the unique element
such that $f(t)=f(t+1)=t$. We have to show that
\begin{equation} \label{eqn-conservativezero}
\begin{aligned}
\xymatrix@R=15pt{
(\Lambda X)_{m} \ar[r]^-{f^*} \ar[d] & (\Lambda X)_{m+1} \ar[d] \\
Y_m  \ar[r]^-{f^*} & Y_{m+1}
}
\end{aligned}
\end{equation}
is homotopy cartesian, or equivalently, that the map
\begin{equation} \label{eqn-conservativeone}
(\Lambda X)_{m} \lra \holim\left(
\begin{aligned}
\xymatrix@R=15pt@C=20pt{
 & (\Lambda X)_{m+1} \ar[d] \\
Y_m  \ar[r]^-{f^*} & Y_{m+1}
}
\end{aligned}
\right)
\end{equation}
determined by~\eqref{eqn-conservativezero} is a weak equivalence.
Now~\eqref{eqn-conservativeone} can be written, up to weak equivalences, in the
form of a map of Spaces $|P| \to |Q|$
induced by a map of simplicial Spaces $f^{(*)}\co P\to Q$. Here $P$ is
the simplicial Space $[n]\mapsto X\langle u\rangle_{m,n}$ and $Q$ is the
simplicial space
\[  [n]~\mapsto~ \holim\left(
\begin{aligned}
\xymatrix@R=15pt@C=20pt{
 & X\langle u\rangle_{m+1,n} \ar[d] \\
Y_m  \ar[r]^-{f^*} & Y_{m+1}
}
\end{aligned}
\right).
\]
Let $g\co [m]\to [m+1]$ be the injective morphism in $\Delta$
which has $t\notin \im(g)$. Then $g^{(*)}f^{(*)}$ is the identity of $P$.
To complete the proof we only need to construct a simplicial homotopy
$Q\times\Delta[1] \lra Q$
from the identity to $f^{(*)}g^{(*)}$. This means that for every
monotone map $c$ from $[n]$ to $[1]$ we need an ``induced'' map $Q_n\to Q_n$. For that we take the map which is
determined by the identity on $Y_m$, the map $(gf)^*\co Y_{m+1}\to Y_{m+1}$ and the map
from $X\langle u \rangle_{m+1,n}$ to $X\langle u \rangle_{m+1,n}$
which is induced by the order-preserving endomorphism of $[m+1]\times[n]$ given by
$(s,t)\mapsto (s,t)$ if $c(t)=0$, and $(s,t)\mapsto(gf(s),t)$ if $c(t)=1$. \qed

\begin{lem}\label{lem-inthomcons}
If $X\xrightarrow{\,u\,} Y$ is a conservative map of simplicial spaces, then the induced map
$X^{\Delta[m]} \to Y^{\Delta[m]}$ is also conservative.
\end{lem}
Once again this uses exponential notation for derived internal hom objects; so $X^{\Delta[m]}$ is the
simplicial space which has $\map(\Delta[m]\times\Delta[n],\varphi X)$ in degree $n$.

\proof Fix $[m]$ and an epimorphism $q\co [n]\to [n-1]$ in $\Delta$. We have to show that
\begin{equation} \label{eqn-havetoshow}
\begin{split}
\xymatrix{  {\map(\Delta[m]\times\Delta[n-1],Y)}   \ar[r]^-u \ar[d] & {\map(\Delta[m]\times\Delta[n-1],Y)}  \ar[d] \\
\map(\Delta[m]\times\Delta[n],\varphi X) \ar[r]^-{u_*} & \map(\Delta[m]\times\Delta[n],\varphi Y)
}
\end{split}
\end{equation}
is homotopy cartesian. Write
\[  \Delta[m]\times\Delta[n]= \colimsub{g\co [k]\to [m]\times[n]} \Delta[k]  \]
where $g\co [k]\to [m]\times[n]$ runs over the nondegenerate simplices in $\Delta[m]\times\Delta[n]$.
(The homotopy limit is taken over $\sQ(m,n)$, the poset of nondegenerate simplices of $\Delta[m]\times\Delta[n]$.)
In the Reedy model structure, this colim is a hocolim, and so
the lower row of~\eqref{eqn-havetoshow} can be re-written (up to weak equivalences)
\[ \holimsub{g\co [k]\to [m]\times[n]} X_k \quad \lra  \holimsub{g\co [k]\to [m]\times[n]} Y_k \]
where we use $X_k\cong\map(\Delta[k],X)$ and $Y_k\cong \map(\Delta[k],Y)$.
Now we propose to write the upper row in a similar manner:
\[ \holimsub{g\co [k]\to [m]\times[n]} X_{k(g)} \quad \lra  \holimsub{g\co [k]\to [m]\times[n]} Y_{k(g)}  \]
where $k(g)$ is determined by the commutative diagram
\begin{equation} \label{eqn-havetoshow2}
\begin{split}
\xymatrix@M=8pt@R=18pt{
[k] \ar@{->>}[d] \ar[r]^-{g} & [m]\times[n] \ar[d]^-{[m]\times q} \\
[k(g)] \ar@{>->}[r] & [m]\times[n-1]
}
\end{split}
\end{equation}
To justify this proposal it suffices to show that the functor $F$ taking $g\in \sQ(m,n)$ to
the injective monotone map $[k(g)]\to [m]$ as in~\eqref{eqn-havetoshow2}, an element of $\sQ(m,n-1)$, is homotopy terminal.
This amounts to showing that certain subposets of $\sQ(m,n)$ have a (weakly) contractible classifying Space.
We can take the view that the objects of $\sQ(m,n)$ are certain nonempty subsets of $[m]\times[n]$,
and similarly for $\sQ(m,n-1)$. Now we have to show that for every $T\in \sQ(m,n-1)$, the full sub-poset
of $\sQ(m,n)$ consisting of the $S\in \sQ(m,n)$ such that $F(S)\supset T$ has a contractible classifying Space.
Here $F(S)$ is simply the image of $S$ under the right-hand vertical arrow in~\eqref{eqn-havetoshow2}.
By an adjunction argument, this is equivalent to showing that for every $T\in \sQ(m,n-1)$, the full sub-poset
$\sQ(m,n,T)$ of $\sQ(m,n)$ consisting of the $S\in \sQ(m,n)$ such that $F(S)=T$ has a contractible classifying Space.

Here we should make a case distinction. Let $x\in [n-1]\subset [n]$ be the unique element such that $q(x)=q(x+1)$.
If $T$ has empty intersection with $[m]\times\{x\}$, then $\sQ(m,n,T)$ has exactly one element, and the classifying Space is contractible.
If that intersection has $r+1$ elements, $r\ge 0$, then $\sQ(m,n,T)$ is isomorphic to $\sQ(r,1,U)$ where
$U$ is all of $[r]\times[0]$. The classifying Space is a concatenation of $2r$ copies of $\Delta[1]$, and so it is contractible.

Having settled that, we can recast~\eqref{eqn-havetoshow} in the form
of a commutative square
\begin{equation} \label{eqn-havetoshow3}
\begin{split}
\xymatrix@C=20pt@R=15pt@M=6pt{
{\holimsub{g\co [k]\to [m]\times[n]} X_{k(g)}}  \ar[r]^-u \ar[d] & {\holimsub{g\co [k]\to [m]\times[n]} Y_{k(g)}}  \ar[d] \\
{\holimsub{g\co [k]\to [m]\times[n]} X_k} \ar[r]^-{u_*} & {\holimsub{g\co [k]\to [m]\times[n]} Y_k}
}
\end{split}
\end{equation}
determined by the collection of commutative squares
\begin{equation} \label{eqn-havetoshow4}
\begin{split}
\xymatrix@C=20pt@R=15pt{
X_{k(g)}  \ar[r]^-u \ar[d] &  Y_{k(g)}  \ar[d] \\
X_k \ar[r]^-{u_*} & Y_k
}
\end{split}
\end{equation}
where $g\in \sQ$, and the vertical arrows are induced by preferred monotone surjections $[k]\to [k(g)]$.
Since the squares~\eqref{eqn-havetoshow4} are all homotopy cartesian by our assumption on $u$,
it follows that~\eqref{eqn-havetoshow3} is homotopy cartesian.
\qed

\begin{prop} If $u\co X\to Y$ is already conservative, then
the second-variable simplicial operators $X\langle u\rangle_{m,0}\to X\langle u\rangle_{m,n}$
in the bisimplicial Space $X\langle u\rangle$ are all weak equivalences. Hence the inclusion
$X\to \Lambda X$ is a weak equivalence.
\end{prop}
\proof By lemma~\ref{lem-inthomcons}, the square
\[
\xymatrix{
X_m \ar[r] \ar[d] & Y_m \ar[d] \\
X^{\Delta[m]} \ar[r] & Y^{\Delta[m]}
}
\]
is homotopy cartesian. Comparison with~\eqref{eqn-expoLam} finishes the proof. \qed

\medskip
Let $\sL(Y)$ be the category of all simplicial Spaces over $Y$ and let
$\sK(Y)\subset \sL(Y)$ be the full subcategory consisting of the objects which are
conservative.

\begin{prop} \label{prop-goodnewLambda} The functor $\Lambda$, viewed as a functor $\sL(Y)\to \sK(Y)$,
is a derived left adjoint for the inclusion $\sK(Y)\to \sL(Y)$. The natural inclusion
$\id\Rightarrow \Lambda$ is the unit of the adjunction.
\end{prop}
\proof Let $X$ be an object of $\sL(Y)$. Write $\iota_X\co X\to \Lambda X$ for the
inclusion. The general theory of derived adjoints
says that we need to establish two properties of $\Lambda$ and $\iota$.
\begin{itemize}
\item[(i)] If $X$ is in $\sK(Y)$, then $\iota_X\co X\to \Lambda X$ is a weak equivalence.
\item[(ii)] For every $X$ in $\sL(Y)$, the map $\Lambda(\iota_X)\co \Lambda X \lra \Lambda(\Lambda X)$
induced by $\iota_X$ is a weak equivalence.
\end{itemize}
Of these, (i) has already been established. Statement (ii) follows from (i) and
\begin{itemize}
\item[(iii)] for every $X$ in $\sL(Y)$ and integer $m\ge 0$, the two maps $\iota_{\Lambda X}$ and
$\Lambda(\iota_X)$ from $\Lambda X$ to $\Lambda(\Lambda X))$ are (weakly) homotopic in degree $m$.
\end{itemize}
In the proof of (iii) we are going to make heavy use of naturality arguments.
Therefore we introduce $\sL$, the arrow category of simplicial Spaces. (An object is a map of simplicial spaces $X\to Y$, and a morphism
is a commutative square. This has a full subcategory $\sK$ spanned by the objects $u\co X\to Y$
where $u$ is conservative. We will write objects of $\sL$ in the form $X\cvar Y$ in order to mark them as objects.
For $\Lambda X$ as in (iii) we will now write $\Lambda(X\cvar Y)$ to show the dependence
on $Y$. The functor $\Lambda$ can be regarded as an endofunctor of $\sL$, if that is convenient.)

For the proof of (iii) we introduce a small subcategory $\sU_m$ of $\sL$. This is the image of the
functor from $\Delta$ to $\sL$ taking $[n]$ to
\[ U(m,n):=\big(\Delta[m]\times\Delta[n]~\cvar~\Delta[m]\big) \]
where the (curved) arrow is the projection.

Note that $\Lambda(U(m,n))$ is $U(m,0)$ up to weak equivalence, by
direct computation. The map $\iota_{U(m,n)}\co U(m,n)\to \Lambda(U(m,n))$ can be identified with the
unique morphism $U(m,n)\to U(m,0)$ in $\sU_m$. Therefore we have a preferred derived natural transformation
\begin{equation} \label{eqn-byedegen}
\begin{split}
\xymatrix@R=12pt@M=4pt{  {\big(U(m,n)\mapsto \star\big)} \ar@{=>}[d] \\
{\big(U(m,n) \mapsto (\Lambda(U(m,n)))_m\big)}
}
\end{split}
\end{equation}
by writing $(\Lambda(U(m,n)))_m\simeq (\Delta[m])_m$ and using the inclusion of the unique nondegenerate element
in the finite set $(\Delta[m])_m$. This helps us in making a crucial observation.
\begin{itemize}
\item[(obs)] For fixed $m\ge 0$, the functor $(X\cvar Y)\mapsto (\Lambda(X\cvar Y))_m$ from $\sL$ to Spaces is the
derived left Kan extension of the constant one-point functor on $\sU_m$.
The derived natural transformation~\eqref{eqn-byedegen} is the unit transformation implied by this
relationship.
\end{itemize}
(Justification of (obs): the Space $(\Lambda(X\cvar Y))_m$ is known to us as the realization of the simplicial Space
$[n] \mapsto \rmap_{\sL}(U(m,n), X \cvar Y)$.
The realization is weakly equivalent to the homotopy colimit of the functor $[n] ~\mapsto~\rmap_{\sL}(U(m,n),X \cvar Y)$
on $\Delta^\op$. But that is also
the preferred formula for the value at $(X\cvar Y)$ of the derived left Kan extension of the
constant functor $\star$ on $\sU_m$ along $\sU_m\hookrightarrow \sL$.)

Because of (obs), statement~(iii) can be deduced from a statement about $\Lambda$ restricted to $\sU_m$.
 \begin{itemize}
\item[(iv)] For $X\cvar Y$ in $\sU_m\subset \sL$, the two maps $\iota_{\Lambda X}$ and
$\Lambda(\iota_X)$ from $\Lambda(X\cvar Y)$ to $\Lambda(\Lambda(X\cvar Y))$ are (weakly) homotopic in degree $m$
by a \emph{derived natural} homotopy.
\end{itemize}
But (iv) is obvious because $\Lambda(\Lambda U(m,n))\simeq \Lambda U(m,n)\simeq U(m,0)$, so that the $m$-th space
of $\Lambda(\Lambda U(m,n))$ is homotopy discrete. \qed

\medskip
Now we are in a good position to prove that conservatization is compatible with products.
In more detail, suppose that we have commutative diagrams of simplicial Spaces
\[
\xymatrix@R=15pt@C=10pt{W \ar[dr] \ar[rr] && X \ar[dl] \\
& Y
}
\quad \quad
\xymatrix@R=15pt@C=10pt{W' \ar[dr] \ar[rr] && X' \ar[dl] \\
& Y'
}
\]
\begin{prop} \label{prop-prodconserv} If $W\to X$ is a conservatization map over $Y$, and $W'\to X'$ is a conservatization map over $Y'$,
then the induced map $W\times W'\to X\times X'$ is a conservatization map over $Y\times Y'$.
\end{prop}
\proof We are told that $X$ is conservative over $Y$, that $X'$ is conservative over $Y'$, and that the maps
$\Lambda W \to \Lambda X$\,, $\Lambda W' \to \Lambda X'$
are weak equivalences. We need to deduce that $\Lambda(W\times W') \to \Lambda(X\times X')$
is a weak equivalence. Fibrant replacement $\varphi$ and the contravariant functor
$\map(\Delta[m]\times\Delta[n], -)$
respect products. The ungeometric realization (applied to Reedy cofibrant simplicial Spaces)
also respects products up to weak equivalence. \qed

\medskip
Let $Y$ be a Segal Space. Let $X \to Y$ be a map of simplicial
Spaces making $X$ into a fiberwise complete Segal Space over $Y$. Write
$\kappa X \to \kappa Y$ for the induced map between their Rezk completions.
\begin{prop}
 The square of simplicial Spaces
\[
\xymatrix@R=16pt{
\Lambda(X) \ar[r] \ar[d] & \ar[d] \Lambda(\kappa X) \\
	Y \ar[r] & \kappa Y
}
\]
is a (degreewise) homotopy pullback square.
\end{prop}
We will use homotopical shorthand for the proof.
In degree $n$, the map from $\Lambda(X)$ to $\Lambda(\kappa X)$
 in the above square is the (ungeometric) realization of the map of simplicial Spaces
\begin{equation}\label{eq:LL}
X^{\Delta[n]} \times^h_{Y^{\Delta[n]}} Y_n \lra (\kappa X)^{\Delta[n]} \times^h_{(\kappa Y)^{\Delta[n]}} \kappa Y_n
\end{equation}
(The exponential means derived inner hom; so $Y_n \subset Y^{\Delta[n]}$ and $\kappa Y_n \subset \kappa Y^{\Delta[n]}$.)

\proof
Since $X \to Y$ is fiberwise complete, the map $X \to \kappa X \times^h_{\kappa Y} Y$ which it induces
is a (degreewise) weak equivalence. Therefore the induced map
\[
X^{\Delta[n]} \to (\kappa X)^{\Delta[n]} \times^h_{(\kappa Y)^{\Delta[n]}} Y^{\Delta[n]}
\]
is also a degreewise weak equivalence for each $n$.
It follows that (\ref{eq:LL}) is the upper row in a homotopy cartesian square whose lower row is
\[  Y_n \to \kappa Y_n \]
viewed as a map between constant simplicial Spaces. Therefore (\ref{eq:LL}) is a lucky case of a map between simplicial Spaces
where we are allowed to interchange realization with taking homotopy fibers. That is, every homotopy fiber of $\Lambda X_n \to \Lambda(\kappa X)_n$ maps by a weak equivalence
to the appropriate homotopy fiber of $Y_n \to \kappa Y_n$, settling the claim.
\qed

\begin{cor} \label{cor-piconserv} \emph{(Notation of definition~\ref{defn-configpi}.)} The forgetful map from
$\cconfig(\pi)$ to $\cconfig(E)$ is a conservatization map over $N\cfin$. \qed
\end{cor}

\end{appendices}

\end{document}